\documentclass[a4paper,12pt]{article}
\usepackage[dvips]{epsfig}
\usepackage{amsfonts}
\usepackage{amsmath}
\usepackage{amssymb}
\makeatletter
\renewcommand{\section}{\@startsection{section}{1}{0.0mm}%
{-\baselineskip}{0.5\baselineskip}{\normalfont\large\bfseries}}
\makeatother
\begin{document}

\title{\textbf{Contributions of Issai Schur to  Analysis}}
\author{\textbf{Harry Dym and Victor Katsnelson\footnote{HD thanks Renee and
Jay Weiss and VK thanks Ruth and Sylvia Shogam for endowing the chairs that
support their respective research. Both authors thank the Minerva Foundation 
for partial support.}}}
\maketitle The name Schur is associated with many terms and
concepts that are widely used in a number of diverse fields of mathematics and
engineering. This survey article focuses on Schur's work in analysis. Here
too, Schur's name is commonplace:  The Schur test and
Schur-Hadamard multipliers (in the study of estimates for
Hermitian forms), Schur convexity, Schur complements, Schur's
results in summation theory for sequences (in particular,  the
fundamental Kojima-Schur theorem), the Schur-Cohn test, the Schur
algorithm, Schur parameters and the Schur interpolation problem for 
functions that are holomorphic and bounded by one in the unit disk. 
In this survey, we shall discuss all of the above mentioned topics and then 
some, as well as  
some of the generalizations that they inspired.  There are nine sections of
 text, each of which is devoted to a separate theme based on Schur's work. 
Each  of these sections has an independent bibliography. There is very little 
overlap. A tenth section presents a  list  of the papers of Schur
 that focus on topics that are commonly considered to be analysis.  
We shall begin
 with a review of Schur's less familiar papers on the
theory of commuting differential operators. 

{\bf Acknowledgement}: The authors extend their thanks to Bernd Kirstein 
for carefully reading the manuscript and spotting a number of misprints.

\setcounter{section}{0} 
\begin{minipage}{15.0cm}
\section{\hspace{-0.5cm}.\hspace{0.2cm}%
{\large Permutable differential operators and fractional powers of
differential operators.\label{PDifOp}}}
\end{minipage}\\[0.18cm]
\setcounter{equation}{0}

Let
\begin{equation}
\label{DifOp1}
 P(y)=p_n(x)\frac{d^{n}y}{dx^{n}}+
p_{n-1}(x)\frac{d^{n-1}y}{dx^{n-1}}+\,\cdots\,+p_{0}(x)y
\end{equation}
%
and
\begin{equation}
\label{DifOp2}
 Q(y)=q_{m}(x)\frac{d^{m}y}{dx^{m}}+
q_{m-1}(x)\frac{d^{m-1}y}{dx^{m-1}}+\,\cdots\,+q_{0}(x)y,
\end{equation}
be formal differential operators, where \(n\geq 0\) and \(m\geq 0\)
are  integers, and
 \(p_k(x)\) and \(q_k(x)\) are complex valued functions.
Then \(Q\)
commutes with \(P\) if \((PQ)(y)=(QP)(y)\). ( It is assumed that the
coefficients \(p_k, q_l\) are are smooth enough, say infinitely
differentiable, so that the product  of the two differential expressions
is defined according to the usual rule for differentiating a
product. The commutativity \(PQ-QP=0\) means that the appropriate
differential expressions, that are constructed from the coefficients
\(p_k, q_l\) according to the usual rules for  differentiating a product,
vanish.)

In \cite{Sch1a} Schur proved the following result:
 \textsl{Let \(P\), \(Q_1\) and \(Q_2\) be differential
 operators of the form \textup{(\ref{DifOp1})} and
 \textup{(\ref{DifOp2})}.
  Assume that each of the operators \(Q_1\) and \(Q_2\)
 commutes with \(P\): \(PQ_1=Q_1P\) and \(PQ_2=Q_2P\).
 Then the operators \(Q_1\) and \(Q_2\) commute with each other}:
 \(Q_1Q_2=Q_2Q_1\).

This result of Schur was forgotten and was rediscovered by S.\,Amitsur  
(\cite{Ami},
Theorem 1) and by I.M.\,Krichever (\cite{Kri1}, Corollary 1 of Theorem 1.2).
(Amitsur does not mention the result of Schur, and Krichever does not
mention either the result of Schur, or the result of Amitsur in \cite{Kri1}, 
but does refer to Amitsur in a subsequent paper \cite{Kri2}.

The method used by Schur to obtain this result is not less
interesting than the result itself. In modern language,  Schur
developed the calculus of formal pseudodifferential operators in
\cite{Sch1a}: for every integer  \(n\) (positive, negative or
zero), Schur considers the formal differential ``Laurent" series
of the form
\begin{equation}
\label{LaurF} F=\sum\limits_{\genfrac{}{}{0pt}{}{-\infty<k\leq n
}{k\ \mathrm{an\ integer}}}f_{k}(x)\,D^{\,k},
\end{equation}
where the coefficients \(f_k(x),\ -\infty<k\leq n\), are smooth
complex-valued functions of \(x\) and \(D=\frac{d}{dx}\).
 (He does not discuss the  existence of an operator
in a  space of functions that corresponds  to this formal
series.)  The sum of two formal ``Laurent" series and the product of
such a series and a complex constant are defined in the usual
way. To define the product \(F\circ G\) of two such  series \(F\)
and
\begin{equation}
\label{LaurG} G=\sum\limits_{\genfrac{}{}{0pt}{}{-\infty<l\leq
m}{l\ \mathrm{an\ integer}}}g_{l}(x)\,D^{\,l},
\end{equation}
one needs a rule for commuting powers of the operator \(D\) with
 powers of the operator of multiplication by the
function \(a(x)\). This
rule is defined by the formulas
\[Da=a(x)D+a^{\prime}(x)I,\]
and
\[
D^{-1}a=a(x)D^{-1}-a^{\prime}(x)D^{-2}+a^{\prime\prime}(x)D^{-3}
+\,\cdots \,+(-1)^{k-1}a^{(k-1)}(x)D^{-k}+\,\cdots\,,
\]
where \(a^{\prime}(x),\,a^{\prime\prime}(x),\,\dots
a^{(k-1)}(x),\,\dots\,\) are the derivatives of the function
\(a(x)\) of the indicated order. The set of all formal
differential  ``Laurent" series provided with such operations
becomes an associative (but not commutative) ring over the field
of complex numbers.
 If  the function \(f_n(x)\) is invertible (in
which case we can and will assume that \(f_n(x)\equiv
1\)), then the formal Laurent series \textup{(\ref{LaurF})} is
invertible, and its inverse is of the form
\begin{equation}
\label{Inv} H=\sum\limits_{-\infty<l\leq -n}h_l(x)D^{l},
\end{equation}
where \(h_{-n}(x)=1\), and the coefficients \(h_k(x)\) are
polynomials in the functions \(f_k(x),\, k<n,\ \) and their
derivatives.

 In particular, a differential operator \(P\) of
the form (\ref{DifOp1}) may be considered as a formal ``Laurent"
series (\ref{LaurF}) whose ``positive" part
 \(F_{+}=\sum\limits_{0\leq k \leq n}f_{k}(x)D^{k}\)
coincides with \(P\) and whose ``negative"
 part
\(F_{-}=\sum\limits_{-\infty< k< 0}f_{k}(x)D^{k}\)
 vanishes. In
\cite{Sch1a}, Schur proved that \textsl{if each of two formal
differential Laurent series \(F_1\) and \(F_2\) commutes with a
differential operator \(P\) of the form \textup{(\ref{DifOp1})\,:}
\(P\circ F_1=F_1\circ P\) and \(P\circ F_2=F_2\circ P\), then
\(F_1\) and \(F_2\) commute with each other\,: \(F_1\circ
F_2=F_2\circ F_1\).} In particular, this result is applicable to
 polynomial differential operators \(Q_1\) and \(Q_2\) of the form
(\ref{DifOp2}) commuting with \(P\). (\(Q_1\) and \(Q_2\) are
considered as differential formal Laurent series whose ``negative"
parts are equal to zero.) Schur gives an explicit description
of the commutant of the differential operator \(P\) (and, even
more generally, the description of the commutant of any formal
differential Laurent series). The notion of \textsf{the
fractional power \(P^{1/n}\) of the differential operator} \(P\)
is involved in this description.

Let \(n\geq 0\) and let \(F\) be a formal differential Laurent series
of the form (\ref{LaurF}). The formal differential Laurent series
\(F^{1/n}\) is defined as the formal differential series \(R\) for which the
equality
\begin{equation}
\label{DefinRoot}%
\underbrace{ R\circ R\circ\,\cdots\,\circ R}\limits_{n\
\mathrm{times}}=F
\end{equation}
holds. In \cite{Sch1a} it is proved that if the function
\((f_n(x))^{1/n}\) exists (in which case we can and  will assume
that \(f_n(x)\equiv 1\)), then such a series \(R=F^{1/n}\) exists
and is of the form
\begin{equation}
\label{DefRoot}%
R=\sum\limits_{{\genfrac{}{}{0pt}{}{-\infty<\rho\leq 1} {\rho\
\mathrm{an\ integer}}}} r_{\rho}(x)D^{\rho}\,,
\end{equation}
where \(r_1(x)\equiv 1\). The coefficients \(r_0(x)\),
\(r_{-1}(x), r_{-2}(x),\,\dots \) can be determined in a
recursive manner as  polynomials of the functions \(
f_{n-1}(x), f_{n-2}(x),\,\dots,\, f_0(x),\,\dots\ \) and their
derivatives. The differential Laurent series \(F^{k/n}\) (\(k\)
an integer) is defined as
\(F^{k/n}\stackrel{\mathrm{def}}{=}(F^{1/n})^k.\)

For example, if
\begin{equation}
\label{StLiouv} L=D^2+q(x)I
\end{equation}
is a Sturm-Liouville differential operator of  second order,
then
\begin{equation}
\label{SqRoot}%
 L^{1/2}=D+s_0(x)I+s_{-1}(x)D^{-1}+s_{-2}(x)D^{-2}
 +s_{-3}(x)D^{-3}+s_{-4}(x)D^{-4}+\,\cdots\ \,,
\end{equation}
where
\begin{equation}
\label{CoefSqR}%
\begin{array}{l}
s_0(x)=0,\, s_{-1}(x)=\frac{1}{2}\,q(x),\
s_{-2}(x)=-\frac{1}{4}\,q^{\,\prime},\
s_{-3}(x)=\frac{1}{8}\,q^{\,\prime\prime}(x)-\frac{1}{8}\,q^2(x),\\[1.2ex]
s_{-4}(x)=-\frac{1}{16}\,q^{\,\prime\prime\prime}(x)+\frac{3}{8}\,q(x)
q^{\,\prime}(x),\,\dots\ .
\end{array}
\end{equation}
Furthermore, \(L^{3/2}=(L^{1/2})^{3}=L\cdot L^{1/2}=L^{1/2}\cdot L\),
and we can calculate
\begin{equation}
\label{3over2}%
L^{3/2}=D^3+t_2(x)D^2+t_1(x)D+t_0(x)I+t_{-1}(x)D^{-1}+\,\cdots\, ,
\end{equation}
where
\begin{equation}
\label{Coef3over2}%
t_2(x)=0,\ t_{1}(x)=\frac{3}{2}\,q(x),\
t_{0}(x)=\frac{3}{4}\,q^{\,\prime}(x),\,
t_{-1}(x)=\frac{1}{8}\,q^{\,\prime\prime}(x)+\frac{3}{8}\,q^{2}(x),\,\dots
\end{equation}

In \cite{Sch1a} it is proved that \textsl{the formal differential
Laurent series \(F\) commutes with a differential operator \(P\)
of the form \textup{(\ref{DifOp1})}  (of order \(n\)) if and only
if \(F\) is of the form
\begin{equation}
\label{LaurRootF} F=\sum\limits_{-\infty<k\leq n}c_{k} P^{k/n}\,,
\end{equation}
where \(c_k\) are complex \textsf{constants} (\(c_k\) do not
depend on \(x\)).} In particular, some series of the form
(\ref{LaurRootF}) can in fact be differential  operators (if by
some \textsl{very special} choice of the \(c_k\) the negative part
\(F_{-}\) of the series (\ref{LaurRootF}) vanishes). These and
only these differential operators  commute with \(P\). Moreover, it is
clear that series of the form (\ref{LaurRootF}) commute with each
other.

The results of Schur on fractional powers of differential
operators were forgotten. The resurgence of interest in this
topic  is related to the
inverse scattering method for solving non-linear evolution
equations. The inverse scattering method was discovered and
applied to the Korteweg - de Vries equation  by C.\,Gardner,
J.\,Green, M.\,Kruskal and R.\,Miura in their famous paper
\cite{GGKM}. This method was then extended to some other
important equations towards the end of the sixties. P.\,Lax
\cite{Lax} developed some machinery (that is now commonly known
as the method of \(L\,\textrm{-}\,A \) pairs, or Lax pairs) that
allows one to use the inverse scattering formalism in a more
organized way. The first step of the Lax method is to express the
given evolution equation in the form
\begin{equation}
\label{LaxRepr}%
\frac{\partial L}{\partial t}=[A,\,L],
\end{equation}
where \(L\) is a differential operator (with respect to \(x\)),
some of whose coefficients depends on \(t\), \(A\) is a
differential operator with respect to \(x\) that does not depend
on \(t\) and (the commutator)
\([A,\,L]=AL-LA\). In subsequent developments, the evolution
equation (\ref{LaxRepr}) was investigated using various analytic
methods drawn from the theory of inverse spectral and scattering problems  
and the
Riemann-Hilbert problem, among others. In their article \cite{GeDi},
I.M.\,Gel'fand and L.A.\,Dikii (=L.\,Dickey) observed that 
fractional powers of differential operators can help in  a
systematic search for pairs \(L\) and \(A\) whose commutator
\([A,\,L]\) is related to a nonlinear evolution equation.  The
idea of Gel'fand and Dikii is to consider the ``positive" part
\((L^{\alpha})_+\) of some fractional power \(L^{\alpha}\) as
such an operator \(A\). Let us explain how the fractional powers
of the Sturm-Liouville operator \(L\) of the form (\ref{StLiouv})
can be applied to construct  the \(L\,\textrm{-}\,A\) pair for the
Korteveg - de Vries equation. Since an operator \(L\) of the form
(\ref{StLiouv}) is  of second order, it suffices to  consider
only integer and half-integer powers of \(L\). Integer powers
do not lead to anything useful: the appropriate \(A\) just commutes
with \(L\). Half-integer powers are more interesting.
According to (\ref{SqRoot})-(\ref{CoefSqR}), \((L^{1/2})_{+}=D\).
The direct computation of the commutator gives: \(
[A,\,L]=q^{\prime}I\) for \(A=D\). The evolution equation
(\ref{LaxRepr}) is of the form \(\frac{\partial q}{\partial t}=
\frac{\partial q}{\partial x}\) in this case. The case
\(A=(L^{3/2})_{+}\) is much more interesting. From
(\ref{3over2})-(\ref{Coef3over2}) it follows that
\begin{equation}
\label{A3over2}%
 A=D^{3}+\frac{3}{2}\,q(x)D+\frac{3}{4}\,q^{\,\prime}(x)I\,.
\end{equation}
The direct calculation of the commutator of the differential
expressions \(A\) and \(L\) of the forms (\ref{A3over2}) and
(\ref{StLiouv}), respectively, gives
\begin{equation}
\label{CalcCom}%
[A,\,L]=\frac{1}{4}\,q^{\,\prime\prime\prime}(x)+\frac{3}{2}\,q(x)\,q^{\prime}(x).
\end{equation}
Thus, the evolution equation (\ref{LaxRepr}) takes the form
\begin{equation}
\label{KdV} \frac{\partial q}{\partial t}=
\frac{1}{4}\,\frac{{\partial}^{3}q}{\partial
x^{3}}+\frac{3}{2}\,q\, \frac{\partial q}{\partial x}.
\end{equation}
This is the Korteweg - de Vries equation. In the paper \cite{GeDi}
a symplectic structure was introduced and a
Hamiltonian formalism was developed. The approach of Gel'fand
and Dikii was further developed by M.\,Adler \cite{Adl}
and by B.M.\,Lebedev and Yu.I.\,Manin \cite{LebMa}. However, the
 results of Schur on permutable differential expressions and on
fractional powers of differential expressions are not mentioned
either in \cite{Adl}, or in \cite{LebMa}, nor are they mentioned
in the well-known surveys \cite{Man}, \cite{Tsuj}, dedicated to
algebraic aspects of non-linear differential equations. The fact
that these results of Schur were largely forgotten may be due to
the lack of a natural area of application for a long time. We
found only one modern source where this aspect of Schur's work is
mentioned:  \textsl{Tata Lectures} by D.\,Mumford. Mumford cites
the paper \cite{Sch1a} in Chapter IIIa\,, \S
 11 of  \cite{Mum} (Proposition 11.7).

The paper \cite{Sch1a} does not discuss the structure of the set
of differential expressions which commute with a given operator
\(P\). The answer ``the differential expressions which commute
with \(P\) are those formal Laurent series in \(P^{1/n}\) for
which ``negative part" vanishes is not satisfactory because it
just replaces the original question by the question ``what is the
structure of formal Laurent series in \(P^{1/n}\) for which
``negative part" vanishes.
 Of course, if \(P\) is a given differential
operator and \(b\) is a polynomial with constant coefficients then
the differential operator \(Q=b(P)\) commutes with \(P\). More
generally, if \(Z\) is any differential operator and \(a\) and
\(b\) are polynomials with constant coefficients then the
operators \(P=a(Z)\) and \(Q=b(Z)\) commute each with other.
However there exist pairs of commuting differential operators
\(P,\,Q\) which are not representable in the form \(P=a(Z)\),
\(Q=b(Z)\). (See formula (1) in \cite{BuCh1}.) The problem of
describing  pairs of commuting differential operators was
essentially solved by J.L.\,Burchnall and T.W.\,Chaundy \cite{BuCh1},
\cite{BuCh2}, \cite{BuCh3} in the twenties. See also
\cite{Bak1}. (The complete answer was obtained for those pairs
\(P,\,Q\) whose orders are coprime.) The answer was expressed in
terms of Abelian functions. In particular, it was proved that the
commuting pair \(P,\,Q\) satisfy the equation
\begin{equation}
 \label{HamCa}%
 r(P,\,Q)=0\,,
\end{equation}
where \(r(\lambda,\,\mu)\) is a (non-zero) polynomial of two
variables with constant coefficients. (This result is known as the 
Burchnall-Chaundy lemma.) The remarkable papers \cite{BuCh1},
\cite{BuCh2}, \cite{BuCh3}, \cite{Bak1} were forgotten. Their
results were rediscovered  and further developed by
I.M.\,Krichever, \cite{Kri1}, \cite{Kri2}, \cite{Kri3},
\cite{Kri4} in the seventies.
(When Krichever started his investigations in this direction,      
he was not aware of the results of Burchnall and Chaundy. In     
his paper \cite{Kri1} he mentioned only the relevant recent works of a       
 group of Moscow  mathematicians. However, in his subsequent papers     
he  referred to  \cite{BuCh1}, \cite{BuCh2},           
\cite{BuCh3} and \cite{Bak1}; see the ``Note in Proof" at the end of       
\cite{Kri2} and references [2-4] in \cite{Kri3}.)              

Thus,  the history of commuting differential expressions,  
which began with the work of Schur \cite{Sch1a},  
 is rich
in forgotten and rediscovered results. 


\begin{minipage}{15.0cm}
\section{\hspace{-0.5cm}.\hspace{0.2cm}%
{\large Generalized limits of infinite sequences \newline%
\hspace*{-0.6ex}and their matrix
transformations.\label{GenLimits}}}
\end{minipage}\\[0.18cm]
\setcounter{equation}{0}

One of the basic notions of  mathematical analysis is the notion
of \textsl{the limit of a sequence of real or complex numbers.} A
sequence \(\{x_k\}_{1\leq k<\infty}\) of complex numbers for which
the \(\lim_{k\to\infty}x_k\) exists is said to be
\textsl{convergent}. A sequence \(\{x_k\}_{1\leq k<\infty}\) of
complex numbers for which \(\sup_{k}|x_k|<\infty\) is said to be
\textsl{bounded}. Let the set of all convergent sequences be
denoted by \(\mathbf{c}\), the set of all bounded sequences be
denoted by \(\mathbf{m}\), and the set of all sequences be denoted
by \(\mathbf{s}\).

 It is clear that each of the sets %
\(\mathbf{ c}\), \(\mathbf {m}\) and \(\mathbf {s}\) is a vector
space, and that \(\mathbf{c\subset m \subset s}\).

Sometimes one needs to define  \textsf{a generalized limit}
\(\textsl{\textsf{R\,-}}\lim_{k\to\infty}x_k\) (according to some
rule \(\textsl{\textsf{R}}\))
 for some sequences
for which the ``usual" limit \(\lim_{k\to\infty}x_k\) may not
 exist. Let
\(\mathbf{c}_{\textsl{\textsf{R}}}\) denote the set of all
sequences \(\{x_k\}_{1\leq k<\infty}\) for which the
\(\textsl{\textsf{R\,-}}\lim_{k\to\infty}x_k\) is defined (in
other words,``exists").
Usually some natural requirements are imposed on such a rule
\(\textsl{\textsf{R}}\). Thus, for example, it is often  required
that the set \(\mathbf{c}_{\textsl{\textsf{R}}}\)  be a vector
space. In this case, if the condition
\[\mathbf{c}\subset \mathbf{c}_{\textsl{\textsf{R}}}\ \ \  %
\mathrm{and}\ \ \ 
\textsl{\textsf{R-\!\!\!\!\!}}\lim_{k\to\infty\ \ \ }\!\!\!\!x_k=
\lim\limits_{k\to\infty}x_k\ \ \ %
\mathrm{for\ all}\ \ \{x_k\}_{1\leq k<\infty}\in \mathbf{c}\,,
\]
is satisfied, then the generalized limit
\(\textsl{\textsf{R-}}\lim\) is said to be \textsl{regular}.

A familiar example of  a generalized limit is the well known
 (C\`esaro)
\(\textsf{\textsl{C}}\)-limit:   Given a sequence
\(\{x_k\}_{1\leq k<\infty}\), the sequence %
\(\{y_k\}_{1\leq k<\infty}\) is defined as%
\[y_n=\frac{x_1+x_2+\,\cdots\,+x_n}{n}
\quad (n=1, \,2,\,3\,,\,\dots\,)\] By definition, the
\(\textsf{\textsl{C}}\)-limit of the sequence \(\{x_k\}\) exists,
if usual limit of the sequence \(\{y_k\}\) exists, and
\[
\textsl{\textsf{C}}\textrm{\ -\!\!\!\!}\lim\limits_{k\to\infty\ \
\ }\!\!\!\!x_k
\stackrel{\mathrm{def}}{=}\lim\limits_{k\to\infty}y_k\,.
\]
It is not difficult to prove that the \(\textsl{\textsf{C}}\)\
-\,limit is regular. There are sequences for which the
\(\textsl{\textsf{C}}\)\,-\,limit exists, but the ``usual" limit
does not exist. For example,  the sequence %
\(x_k\stackrel{\mathrm{def}}{=}\frac{1+(-1)^k}{2}\) does not tend
to a limit as $k \rightarrow \infty $
, but the C\`esaro limit exists, and %
\(\textsl{\textsf{C}}\textrm{\
-}\lim\limits_{k\to\infty}x_k=\frac{1}{2}\).

Matrix transformations can be used to define generalized limits.
Let \(A\) be an infinite matrix:
\begin{equation}
\label{MatrTransf}%
 A= \left[
\begin{array}{ccccc}
a_{11}&a_{12}&\ \cdots\ &a_{1k}\ & \cdots\ \\[1.5ex]
a_{21}&a_{22}&\ \cdots\ &a_{2k}\ &\cdots\ \\[1.5ex]
\cdot   & \cdot & \cdot & \cdot & \cdot \\[1.5ex]
a_{n1}  & a_{n2} &  \cdots & a_{nk} & \cdots \\[1.5ex]
\cdot   & \cdot & \cdot & \cdot & \cdot
\end{array}
\right]\, ,
\end{equation}
where the matrix entries \(a_{jk}\) are real or complex numbers.
The matrix transformation \(x\rightarrow y=Ax\) is defined for
those sequences \(x=\{x_k\}_{1\leq k <\infty}\) for which all the
series \(\sum\limits_{1\leq k<\infty}a_{nk}x_{k},\
n=1,\,2,\,3,\,\dots\,,\) converge. The resulting sequence
\(y=Ax\), \(y=\{y_k\}_{1\leq k<\infty}\) is defined by
\(y_n=\sum\limits_{1\leq k<\infty}a_{nk}x_{k},\
(n=1,\,2,\,3,\,\dots\,).\) It is clear that the domain of
definition \({\cal{D}}_A\) of the matrix transformation generated
by the matrix \(A\) is a vector space, \({\cal{D}}_A\subset
\mathbf{s}.\) Moreover, there is a natural generalized limit
associated with each such
infinite matrix \(A\)
(that we denote as \(\textsl{\textsf{A}}\)-limit and which we
shall refer to as \textsl{the matrix generalized limit generated
by the matrix \(A\)} ). Namely, by definition, the
\textsl{\(\textsl{\textsf{A}}\)-limit of a sequence
\(x=\{x_k\}_{1\leq k<\infty}\) exists}, if \(x\in{\cal{D}}_A\)
(i.e. the matrix transformation \(Ax\) is defined), and the
sequence \(y=Ax\) is convergent: \(y\in\mathbf{c}\). By
definition,
\begin{equation}
\label{DefAlim}%
\textsl{\textsf{A}}\,\textrm{-}\!\!\!\!\lim\limits_{k\to\infty\ \
\ }\!\!\!\!x_k=\lim\limits_{k\to\infty}y_k\,.
\end{equation}
The C\`esaro generalized limit (\(\textsl{\textsf{C}}\)\,-limit)
can be considered as the matrix generalized limit generated by the
lower triangular matrix \(A\) for which \(\displaystyle
a_{nk}=\frac{1}{n}\ \ \mathrm{for}\ \ k=1,\,2,\,\dots\,,\,n\) and
\(a_{nk}=0\ \ \mathrm{for}\ \ k>n\,\). The systematic
investigation of matrix generalized limits was initiated by
O.Toeplitz, \cite{Toep}. A fundamental contribution to the theory
of matrix generalized limits was made by Schur. In \cite{Sch16a}
he introduced three classes of matrix transformations:
\textsl{convergence preserving, convergence generating} and
\textsl{regular}.

  A matrix transformation \(x\rightarrow Ax\) is said to
be \vspace{-1.2ex}
\begin{enumerate}
\item
\textsl{convergence preserving}, if it is defined for every
sequence \(x\in\textbf{c}\), and for \(x\in\textbf{c}\) the
sequence \(y=Ax\) belongs to \(\textbf{c}\) as
well.\vspace{-1.2ex}
\item
\textsl{convergence generating}, if it is defined  for every
sequence \(x\in\textbf{m}\), and for \(x\in\textbf{m}\) the
sequence \(y=Ax\) belongs to \(\textbf{c}\).\vspace{-1.2ex}
\item
\textsl{regular}, if it is convergence preserving and, moreover,
if \(x\in\mathbf{c}\) and \(y=Ax\), then the equality
\(\lim\limits_{k\to\infty}y_k=\lim\limits_{k\to\infty}x_k\) holds.
\end{enumerate}

Schur obtained necessary and sufficient conditions for a matrix
transformation \(x\rightarrow Ax\) to belong to each of these
three classes. These conditions are  presented in the following
three theorems that are taken from \cite{Sch16a}. They are
formulated in terms of the numbers
\begin{equation}
\label{RowValue}%
{\sigma}_n=\sum\limits_{1\leq k<\infty}a_{nk} \quad {\rm and}
\quad {\zeta}_n=\sum\limits_{1\leq k<\infty}|a_{nk}|\,
\end{equation}
for those \(n,\ 1\leq n<\infty\), for which these values exist.
The values \({\sigma}_n\) are said to be the \textsl{row sums};
the values \({\zeta}_n\) are said to be the \textsl{row norms}.

\textsf{THEOREM I}. \textsl{The matrix transformation \(A\) is
convergence preserving if and only if the following three
conditions are satisfied}:
\begin{enumerate}
\item
\textsl{For every \(k\) the following limit exists}
\begin{equation}
\label{ColLim} a_{k}\stackrel{\mathrm{def}}{=}
\lim\limits_{n\to\infty}{a_{nk}}\,.
\end{equation}
\item
\textsl{The row sums \({\sigma}_k\) tend to the finite limit}
\(\sigma\):
\begin{equation}
\label{LimRoS} \sigma=\lim\limits_{k\to\infty}{\sigma}_k\,.
\end{equation}
\item
\textsl{The sequence of  row norms is bounded}:
\begin{equation}
\label{RowNorBound} \sup\limits_{1\leq
n<\infty}{\zeta}_{\,n}<\infty.
\end{equation}
\end{enumerate}
\textsl{If these conditions are satisfied, then the series
\(\sum\limits_{1\leq k<\infty}a_k\) converges absolutely and, if 
\begin{equation}
\label{SumColLim}
\alpha\stackrel{\mathrm{def}}{=}\sum\limits_{1\leq k<\infty}a_k\,,
\end{equation}
then for every convergent sequence \(\{x_k\}\)}
\begin{equation}
\label{LimY}%
\lim\limits_{n\to\infty}\sum\limits_{1\leq k<\infty}a_{nk}x_k=
 (\sigma - \alpha)\lim\limits_{k\to\infty}x_k+\sum\limits_{1\leq
 k<\infty}a_{k}x_{k}\,.
\end{equation}

\textsf{THEOREM II}. \textsl{The convergence preserving matrix
transformation \(A\) is regular if and only if all the column
limits \(a_k\) defined in  (\ref{ColLim}) are equal to zero:
\begin{equation}
\label{ColLimZero}%
a_k=0\qquad (k=1,\,2,\,3,\,\dots\,)\,,
\end{equation}
and the limit \(\sigma\) of the row sums \({\sigma}_n\) defined in
(\ref{LimRoS}) is  equal to 1:}
\begin{equation}
\label{LimRoSOne}%
\sigma=1\,.
\end{equation}

\textsf{THEOREM III.} \textsl{The matrix transformation \(A\) is
convergence generating if and only the three assumptions of
Theorem I  are satisfied and
 \textit{the series \(\sum\limits_{1\leq
k<\infty}|a_{nk}|,\,n=1,\,2,\,3,\,\dots\,\), converge uniformly
with respect to \(n\).} In this case}
\begin{equation}
\label{LimVal}%
\lim\limits_{n\to\infty}\sum\limits_{1\leq k<\infty}a_{nk}x_k=
\sum\limits_{1\leq k<\infty}a_{k}x_{k}\,.
\end{equation}

Theorem \textsf{II}  was formulated and proved by O.\,Toeplitz in
(\cite{Toep}). However, Toeplitz considered only lower-triangular
matrices \(A\) . Theorem \textsf{II} is commonly known as the
Toeplitz theorem or as the Silverman-Toeplitz theorem, since  part
of Theorem \textsf{II} was obtained also by L.L.\,Silverman in his
PhD thesis, \cite{Silv}. Theorem \textsf{I} is known as the
Schur-Kojima theorem. (Part of Theorem \textsf{I} was also
obtained by T.Kojima for lower-triangular matrices.) The paper
\cite{Koj} by Kojima was published earlier than the paper
\cite{Sch16a} by Schur. However, in a footnote on the last page of
\cite{Sch16a}, Schur remarks that he only became aware of the
paper \cite{Koj} while reading the proofs of his own paper. The
matrices \(A\) which correspond to convergent generated
transformations are called \textsl{Schur matrices} in  \cite{Pet}.
There is a rich literature dedicated to matrix generalized limits
and to matrix summation methods. (If a considered sequence is a
sequence of partial sums of a series,  then the terminology
``generalized summation method" is used instead of ``generalized
limit" or
 ``generalized limitation method".) We mention only the books \cite{Bo}, 
\cite{Coo}, \cite{Har},  \cite{Pet}, \cite{Pey} and \cite{Zel}. In
all these books, the sections that deal with the basic theory of
generalized limits and generalized summation methods
 cite the results of Schur
and refer to him as one of the founders of this theory.

  In a footnote near the beginning of his
paper \cite{Sch16a}, Schur notes that his considerations have many
points in common with the considerations of H.\,Lebesgue and
H.\,Hahn, dedicated to the sequence of integral transformations of
the form
\[y_n(r)=\int\limits_{a}^{b}A_n(r,\,s)\,ds\,.\]
He also considers  some applications of his Theorems \textsf{I -
III} to the multiplication of series and to Tauberian theorems. In
particular, he derives the Tauberian therem by Tauber (about power
series) from Theorem \textsf{II}.

In  his other paper \cite{Sch6a}, Schur consider H\"older and
C\`esaro limit methods of \(r\)-th order and proves that these
limit methods are equivalent.

Given a sequence \(x_1,\,x_2,\,x_3,\,\dots\,\) of real or complex
numbers, we form the sequences
\[
h^{(1)}_n=\frac{x_1+x_2+\,\cdots\,+x_n}{n}\,,\qquad h^{(2)}_n=
\frac{h^{(1)}_1+h^{(1)}_2+\,\cdots\,+h^{(1)}_n}{n}\,,
\]
\[h^{(3)}_n=
\frac{h^{(2)}_1+h^{(2)}_2+\,\cdots\,+h^{(2)}_n}{n}\,,\quad
\dots\quad,\, h^{(r)}_n=
\frac{h^{(r-1)}_1+h^{(r-1)}_2+\,\cdots\,,\,+h^{(r-1)}_n}{n}\,\cdot
\]
The sequence \(h^{(r)}_1,\,h^{(r)}_2,\,\dots\,,h^{(r)}_k,\,\dots\)
is said to be \textsl{the sequence of H\"older means of order
\(r\)} (constructed from the initial sequence
\(x_1,\,x_2,\,x_3,\,\dots\,\)). Another class of sequences can be
constructed as follows. Let
\[
s^{(1)}_n=x_1+x_2+\,\cdots\,+x_n,\quad
s^{(2)}_n=s^{(1)}_1+s^{(1)}_2+\,\cdots\,+s^{(1)}_n\,,
\]
\vspace{1.5ex}
\[
s^{(3)}_n=s^{(2)}_1+s^{(2)}_2+\,\cdots\,+s^{(2)}_n\,,\quad\cdots\quad
,\, s^{(r)}_n=s^{(r-1)}_1+s^{(r-1)}_2+\,\cdots\,+s^{(r-1)}_n\,,
\]
and set
\[c^{(r)}=\frac{s^{(r)}}{\left({n+r-1}\atop{r}\right)}
\]
The sequence \(c^{(r)}_1,\,c^{(r)}_2,\,\dots\,,c^{(r)}_k,\,\dots\)
is said to be \textsl{the sequence of C\`esaro means of order
\(r\)} (constructed from the initial sequence
\(x_1,\,x_2,\,x_3,\,\dots\,\)). The transformations
\[\{x_1,\,x_2,\,\dots\,,\,x_k,\,\dots\,\}\rightarrow
\{h^{(r)}_1,\,h^{(r)}_2,\,\dots\,,\,h^{(r)}_k,\,\dots\,\}
\]
and
\[\{x_1,\,x_2,\,\dots\,,\,x_k,\,\dots\,\}\rightarrow
\{c^{(r)}_1,\,c^{(r)}_2,\,\dots\,,\,c^{(r)}_k,\,\dots\,\}
\]
can be considered as matrix transformations based on appropriately
defined matrices that we denote by \(\textsf{H}^{(r)}\) and
\(\textsf{C}^{(r)}\), respectively. These matrices are
lower-triangular. Both generalized limits
\(\textsf{\textsl{H}}^{(r)}\)-limit and
\(\textsf{\textsl{C}}^{(r)}\)-limit are regular. In \cite{Sch6a}
it is shown that these generalized limits are equivalent in the
following sense:

\textsl{Let a sequence \(x_1,\,x_2,\,x_3,\,\dots\,\) and a natural
number \(r\) be given. Then the sequence of C\`esaro means
\(\{c^{(r)}_1,\,c^{(r)}_2,\,\dots\,,\,c^{(r)}_k,\,\dots\,\}\)
tends to a finite limit if and only if the sequence of H\"older
means \(\{h^{(r)}_1,\,h^{(r)}_2,\,\dots\,,\,h^{(r)}_k,\,\dots\,\}
\) tends to a finite limit.  Moreover, in this case, the two
limits must agree.}

Schur obtained this result by showing
that both the matrices
\((\textsf{H}^{(r)})^{-1}\cdot\textsf{C}^{(r)}\) and
\((\textsf{C}^{(r)})^{-1}\cdot\textsf{H}^{(r)}\) satisfy the
assumptions of Theorem \textsf{II} (the Toeplitz regularity
criterion). Thus, the appropriate matrix transformations are
regular.

This result by Schur was not new. At the time that the paper
\cite{Sch6a} was published  proofs of the equivalency of
H\"older's and C\`esaro's methods had already been obtained by
K.\,Knopp, by W.\,Schnee and by W.B.\,Ford. However, these proofs
were very computational, very involved and not very transparent.

\setcounter{equation}{0}
\begin{minipage}{15.0cm}
\section{\hspace{-0.5cm}.\hspace{0.2cm}%
{\large Estimates for matrix and integral operators, \newline%
\hspace*{-0.6ex}bilinear forms and related
inequalities.\label{SchurTes}}}
\end{minipage}\\[0.18cm]

The terms \textsl{Schur test, Schur (or Hadamard-Schur)
multiplication of matrices, Schur (or Hadamard-Schur) multipliers}
are all related to Schur's contributions to the estimates of
operators and bilinear forms, see \cite{Sch4a}. In this section we
consider the  Schur test. Results related to the Schur (or
Schur-Hadamard) product will be considered in the next section.

Let \(A=\|a_{jk}\|\) be a matrix, finite or infinite, with real or
complex entries. This matrix generates the bilinear form
\begin{equation}  
\label{BilF}  
A(x,y)=\sum\limits_{j,k}a_{jk}x_ky_j\,,
\end{equation}    
where \(x\) and \(y\) are vectors with entries \(\{x_j\}\) and
\(\{y_k\}\) that are real or complex. The matrix \(A\) also
generates the linear operator
\begin{equation}  
\label{LinOp}  
x\rightarrow Ax, \quad \mathrm{where} \quad
(Ax)_j=\sum\limits_{k}a_{jk}x_k\,.
\end{equation}  
If the matrix \(A\) is not finite, we consider only
\textsl{finite} vectors \(x\) and \(y\), i.e.,  vectors with only
finitely many nonzero entries. This allows us to avoid troubles
related to the convergence of infinite sums.
 If the
sets of vectors \(x\) and \(y\) are provided with  norms, then a
problem of interest is to estimate the bilinear form (\ref{BilF})
in terms of the norms of the vectors \(x\) and \(y\). In
particular, the sets of vectors \(x\) and \(y\) can be provided
with \(l_2\) norms:
\begin{equation}  
\label{LN} 
\|x\|_{l_2}=\Big\{\sum\limits_{j}|x_j|^2\Big\}^{1/2}\,,\quad
\|y\|_{l_2}=\Big\{\sum\limits_{k}|y_k|^2\Big\}^{1/2}\,.
\end{equation}  
If the estimate
\begin{equation}  
\label{BE}  
|A(x,y)|\le C\|x\|_{l^2}\|y\|_{l^2}
\end{equation}  
holds for every pair of vectors \(x\) and \(y\) for some constant
\(C<\infty\), then the bilinear form (\ref{BilF}) is said to be
\textsl{bounded}. The smallest constant \(C\), for which the
inequality (\ref{BE}) holds, is denoted by \(C_A\) and is termed
the norm of the bilinear form \(A(x,y)\):
\begin{equation}  
\label{NBF}  
 C_A=\sup\limits_{x\neq 0,\,y\neq 0}\frac{|A(x,y)|}{\|x\|_{l_2}\|y\|_{l_2}}\, .
\end{equation}  
 The norm
\(C_A\) of the bilinear form generated by the matrix \(A\)
coincides with the norm of the linear operator generated by this
matrix,  considered as a linear operator acting from \(l_2\) into
\(l_2\):
\begin{equation}  
\label{NOP}  
 \|A\|_{l_2\to l_2}=\sup\limits_{x\neq 0}\frac{\|Ax\|_{l_2}}{\|x\|_{l_2}}\,.
\end{equation}   
The cases in which it is possible to express the norm \(C_A\) in
terms of the entries of the matrix \(A\) are very rare. Thus, the
problem of estimating the value of  \(C_A\) in terms of the matrix
entries is a very important problem. In particular, if the matrix
\(A\) is infinite, it is important to recognize whether the value
\(C_A\) is finite or not. Schur made important contributions to this
circle of problems.

 In \cite{Sch4a}\,(\S 2, Theorem I), the following estimate was
obtained.\\[3.5ex]
\textsf{THEOREM} (\textsf{The Schur test}).
\begin{slshape}
 Let \(A=\|a_{jk}\|\) be a matrix, and let
\begin{equation}  
\label{EV}  
\zeta(A)=\sup\limits_{j}\sum\limits_{k}|a_{jk}|\,,\quad 
\kappa(A)=\sup\limits_{k}\sum\limits_{j}|a_{jk}|\,.
\end{equation}
Then
\begin{equation} 
\label{ST} 
C_A\le \sqrt{\zeta(A)\kappa(A)}\,.
\end{equation} 
\end{slshape}

It is enough to prove the estimate (\ref{ST}) for finite matrices
\(A\) (of  arbitrary size). The proof of the estimate (\ref{ST})
that was obtained in \cite{Sch4a} is based on the fact that
\begin{equation} 
\label{EE}  
 C_A =\sqrt{ {\lambda}_{\max}},
\end{equation}  
where \({\lambda}_{\max}\) is the largest eigenvalue of the matrix
\(B=A^*A\). Let \(\boldsymbol{\xi}=\{\xi_k\}\) be the eigenvector
which corresponds to the eigenvalue \({\lambda}_{\max}\):
\begin{equation*}
\lambda_{\max}\boldsymbol{\xi}=B\boldsymbol{\xi}\,.
\end{equation*}
Let \(|{\xi}_p|=\max\limits_{k}|\xi_k|\). Then, since %
\(\lambda_{\max}|\xi_p|\leq
\Big(\sum\limits_{k}|b_{pk}|\Big)|\xi_p|\, \), it is easily seen
that
\begin{equation*}
\lambda_{\max}\le \sum\limits_{k}|b_{pk}|\,,
\end{equation*}
where the \(\{b_{jk}\}\) are the entries of the matrix \(B=A^*A\):
\(b_{jk}=\sum\limits_{r}\overline{a}_{rj}a_{rk}\,.\) Thus,
\begin{equation*}
\sum\limits_k|b_{pk}|\le
\sum\limits_k\sum\limits_r|\overline{a}_{rp}||a_{rk}|\le
\Big(\sum\limits_{r}|a_{rp}|\Big)\cdot
\Big(\max\limits_{r}\sum\limits_k|a_{rk}|\Big)\le %
\kappa(A)\cdot\zeta(A)\,.
\end{equation*}
This completes the proof.

Another proof, which does not use the equality (\ref{EE}), is even
shorter:

\begin{eqnarray}  
\label{CI}
\big|A(x,y)\big| &\le&  
\sum\limits_{j,k}|a_{jk}|\cdot |x_k|\cdot |y_j| \notag\\
&=&  
\sum\limits_{j,k}\big(|a_{jk}|^{1/2}|x_k|\big)\cdot  
\big(|a_{jk}|^{1/2}|y_j|\big)\notag \\ & = &
\Big(\sum\limits_{j,k}|a_{jk}||x_k|^2\Big)^{1/2}\cdot 
\Big(\sum\limits_{j,k}|a_{jk}||y_j|^2\Big)^{1/2} \\ &\le&
\Big(\sup\limits_k\sum\limits_j|a_{jk}|\cdot  
\sum\limits_k|x_k|^2\Big)^{1/2}\cdot \notag
\Big(\sup\limits_j\sum\limits_k|a_{jk}|\cdot  
\sum\limits_j|y_j|^2\Big)^{1/2}\notag \\ &=&
\sqrt{\kappa(A)\zeta(A)}\|x\|_{l_2}\|y\|_{l_2}\,.\notag
\end{eqnarray}  

The estimate (\ref{ST}) can be considered as a special case of an
interpolation theorem that is obtained by introducing the \(l_1\)
and \(l_{\infty}\) norms. If \(x=\{x_k\}\) is a finite sequence of
real or complex numbers, then these norms are defined by the usual
rules:
\begin{equation}  
\label{NM}  
\|x\|_{l_1}=\sum\limits_{k}|x_k|\quad {\rm and} \quad
\|x\|_{l_\infty}=\sup\limits_{k}|x_k|\,,
\end{equation}   
respectively. If \(A\) is a matrix, we can consider the linear
operator generated by this matrix as an operator acting in the
space \(l_1\) as well as an operator acting in the space
\(l_{\infty}\). The corresponding norms \(\|A\|_{l_1\to l_1}\) and
\(\|A\|_{l_{\infty}\to l_{\infty}}\) are defined by the formulas
\begin{equation*}
\|A\|_{l_1\to
l_1}=\sup\limits_{x\neq 0}\frac{\|Ax\|_{l_1}}{\|x\|_{l_1}}\quad  
{\rm and} \quad \|A\|_{l_{\infty}\to
l_{\infty}}=\sup\limits_{x\neq
0}\frac{\|Ax\|_{l_{\infty}}}{\|x\|_{l_{\infty}}}\,\cdot
\end{equation*}
Unlike the norm\(\|A\|_{l_2\to l_2}\), the norms \(\|A\|_{l_1\to
l_1}\) and \(\|A\|_{l_{\infty}\to l_{\infty}}\) can be expressed
explicitly in terms of the matrix entries \(\{a_{jk}\}\):
\begin{equation*}
\|A\|_{l_1\to l_1}=\kappa(A)\quad {\rm and} \quad
\|A\|_{l_{\infty}\to l_{\infty}}=\zeta(A)\,,
\end{equation*}
where the numbers \(\zeta(A)\) and \(\kappa(A)\) are defined in
(\ref{EV}). The estimate (\ref{ST}) takes the form
\begin{equation}  
\label{IT}  
\|A\|_{l_2\to l_2}\le \sqrt{\|A\|_{l_1\to
l_1}\cdot\|A\|_{l_{\infty}\to l_{\infty}}}\,.
\end{equation}  
The inequality (\ref{IT}) is a direct consequence of the
\textsf{M.\,Riesz' Convexity Theorem.} To apply
this theorem, let
 \(\|A\|_{l_p\to l_q}\) denote the norm of the operator, generated
by a matrix \(A\), considered as an operator from \(l_p\) into
\(l_q\) for \(1\le p\le \infty,\,1\le q\le \infty\). Then, Riesz'
theorem states that \(\log\,\|A\|_{l_p\to l_q}\) is a convex
function of the variables \(\alpha=1/p\) and \(\beta=1/q\) in the
square \(0\le \alpha\le 1,\,0\le \beta \le 1\). This theorem
can be found in \cite{HLPa},\,Chapter VIII, sec.\,
\textbf{8.13.}  G.O.\,Thorin, \cite{Tho}, found a very beautiful
and ingenious proof of this theorem using a new method based on
Hadamard's Three Circles Theorem from complex analysis. Therefore
this theorem is also called the Riesz-Thorin Convexity Theorem.
Now this theorem is presented in many sources, and even in
textbooks.  The Riesz-Thorin Convexity Theorem belongs to a
general class of interpolation theorems for linear operators. A
typical interpolation theorem for linear operators deals with a
linear operator that is defined by a certain analytic expression,
for example by a certain matrix or kernel, but is considered not
in a fixed space, but in a whole ``scale" of spaces. A typical
interpolation theorem claims that if the linear operator,
generated by the given expression, is bounded in  two spaces of
the considered ``scale of spaces", then it also is bounded in all
the
  ``intermediate"  spaces. Moreover, the norm
 of the operator in the ``intermediate" spaces is estimated through the norms
 of the operators in the original two spaces. The Riesz-Thorin theorem
states that the spaces \(l_p\) with \(1< p< \infty\) are
``intermediate" for
 the pair of spaces \(l_1\) and \(l_{\infty}\).

The estimate (\ref{ST}) can also be
 considered as a special case of another interpolation theorem for linear
 operators, the so-called interpolation theorem for modular spaces.
 This theorem is based on quite another circle of ideas that are
more geometrical in nature  and was
 partially inspired by Schur's work (\cite{Sch18b}). We will discuss this
 in the next section.

 For practical application, the ``weighted" version of the
 Schur estimate (\ref{ST}) is useful.
In fact, this version was also considered in (\cite{Sch4a}) ( but
not as explicitly, as the ``unweighted" version).
 In the weighted version, a positive
 sequence \(\{r_k\},\,\,r_k>0,\) appears
and the ``weighted" \(l_1\)- and
 \(l_{\infty}\)-norms
\begin{equation}  
\label{WN}  
\|x\|_{l_1,r}=\sum\limits_{k}|x_k|\cdot r_k
\qquad \textrm{and}\qquad \|x\|_{l_{\infty},\,r^{-1}}=%
\sup\limits_k\frac{\|x_k\|}{\,r_k}
\end{equation}  
are considered.

\noindent  
\textsf{THEOREM} (\textsf{The weighted Schur test}).
\begin{slshape}
 Let \(A=[a_{jk}]\) be a matrix and let \(r_k\) be a sequence of strictly
 positive numbers: \(r_k>0\).
 Let
\begin{equation}  
\label{EVM}  
\zeta_r(A)=\sup\limits_{j}\frac{1}{r_j}\cdot\sum\limits_{k}|a_{jk}|\cdot
r_k\,
\quad {\rm and}\quad
\kappa_r(A)=\sup\limits_{k}\frac{1}{r_k}\sum\limits_{j}|a_{jk}|\cdot
r_j\,.
\end{equation}
Then the value \(C_A\), defined in \textup{(\ref{NBF})} is subject
to the bound
\begin{equation} 
\label{STW} 
C_A\le \sqrt{\zeta_r(A)\kappa_r(A)}\, .
\end{equation} 
\end{slshape}

It is easy to see that
\begin{equation*}  
\label{WNR} 
\zeta_r(A)=
\sup\limits_{x\neq 0}\frac{\|Ax\|_{l_1,\,r}}{\|x\|_{l_1,\,r}}\quad 
{\rm and} \quad
\kappa_r(A)=\sup\limits_{x\neq 0}
\frac{\|Ax\|_{l_{\infty},\,r^{-1}}}{\|x\|_{l_{\infty},\,r^{-1}}}\,.
\end{equation*}  
 Thus, the estimate (\ref{STW}) can be presented in
the form
\begin{equation}  
\label{ITW}  
\|A\|_{l_2\to l_2}\le \sqrt{\|A\|_{l_1,\,r\to l_1,\,r}\cdot
\|A\|_{l_{\infty},\,r^{-1}\to l_{\infty},\,r^{-1}}}\, .
\end{equation}   
The inequality (\ref{ITW}) is also an ``interpolation" inequality.
It shows that the space \(l_2\) is an ``intermediate" space,
between the spaces \(l_{1,\,r}\) and \(l_{{\infty},\,r^{-1}}\).

The inequality (\ref{STW}) can be proved in much the same way as
the  special case (\ref{ST}).

As an example, we consider a Toeplitz matrix, i.e., a matrix \(A\)
of the form \(a_{jk}=w_{j-k}\). The Schur test leads to the
estimate
\begin{equation*}
  C_A\le \sum\limits_{l}|w_l|\, .
\end{equation*}

The same bound holds for Hankel matrices, i.e., matrices \(A\) of
the form \(a_{jk}=w_{j+k}\)\, .

As a second example, let us consider the Hilbert matrix
\(H^+=\left[\dfrac{1}{j+k-1}\right]_{j,k=1}^{\infty}\). For this
matrix, \(\sum\limits_{k}|h^+_{jk}|=\infty\), so the ``unweighted"
Schur test does not work. However, if we chose \(r_l=l^{-\alpha}\)
with a fixed \(\alpha\in (0,\,1)\), then \(\sup\limits_{1\le
j<\infty }\left(\,j^{\alpha}\sum\limits_{1\le
k<\infty}\dfrac{k^{-\alpha}}{j+k}\right)=s(\alpha)<\infty\). Thus,
\(\|H^+\|\le s(\alpha).\) Then we can optimize the estimate by
choosing the ``best" \(\alpha\). In the discrete case the precise
value \(s(\alpha)\) is unknown. Nevertheless, it is reasonable to
choose \(\alpha=1/2\), since this is the optimum value for the
continuous analogue of the matrix \(H^+\):
\begin{equation}  
\label{CA}  
x^{\alpha}\,\int\limits_{0}^{\infty}\frac{1}{x+y}y^{-\alpha}dy=
\frac{\pi}{\sin{\pi\alpha}}\,,\quad
\min\limits_{\alpha\in(0,\,1)}\dfrac{\pi}{\sin{\pi\alpha}}=\pi
\text{ is attained at
 the point } \alpha=1/2\,.
\end{equation}  
Some other applications of the Schur test can be found in
\cite{BiSo}, Chapter 2, Section \textbf{10}.

 Schur used the estimate (\ref{ST}) in (\cite{Sch4a}) to study the infinite
Hilbert forms
\begin{equation}
\label{HF}  
H^-=\sum_{\substack{p,q=1\\ p\not=q\
}}^{\infty}\frac{x_py_q}{p-q},\qquad
H^+=\sum\limits_{p,q=1}^{\infty} \frac{x_py_q}{p+q-1},
\end{equation}
and the generalized Hilbert forms
\begin{equation*}  
\label{GHF} 
H_{\lambda}^-=\sum\limits_{p,q=1}^{\infty}\frac{x_py_q}{p-q+\lambda
},\qquad H_{\lambda}^+=\sum\limits_{p,q=1}^{\infty}
\frac{x_py_q}{p+q-1+\lambda}, \qquad (0<\lambda<1).
\end{equation*}  
For the Hermitian matrices  
\[N=(H^+)^*H^++(H^-)^*H^-=[n_{pq}]_{p,q=1}^{\infty} \quad \mbox{and}
 \quad 
 N_{\lambda}=(H_{\lambda}^+)^*H_{\lambda}^++(H_{\lambda}^-)^*H_{\lambda}^-,\]
 the conditions
 \begin{equation*}  
 \sum\limits_{1\le q<\infty}|n_{pq}|<
 3\sum_{\substack{q=-\infty\\q\not=p}}^{\infty}\frac{1}{(p-q)^2}=
 {\pi}^2,\quad
 \sum\limits_{1\le q<\infty}|(n_{\lambda})_{pq}|\le \sum\limits_{-\infty<r<\infty}
 \frac{1}{(r+\lambda)^2}=\frac{\pi^2}{\sin^{\,2}{\pi\lambda}}
 \end{equation*}   
 are satisfied for every \(p\). According to (\ref{ST}), the estimates
 \begin{equation}  
\label{HE} C_{H^+}\le\pi,\quad C_{H^-}\le\pi,\quad
C_{H^+_{\lambda}}\le \frac{\pi}{\sin{\pi\lambda}},\quad
C_{H^-_{\lambda}}\le \frac{\pi}{\sin{\pi\lambda}}
 \end{equation}  
 hold. It turns out that in fact  equality prevails in the first two inequalities
 in (\ref{HE}), i.e.,
 the method of Schur gives the exact values for the norms of
 the matrices \(H^+,\,H^-\). (See \cite{HLPa}, Chapter IX). It should be remarked
 that an essential part of Chapters VIII and IX of \cite{HLPa} is based on results
 of the paper \cite{Sch4a}.

In {\S} 6 of \cite{Sch4a}, the infinite quadratic form
\begin{equation}
\label{QuFo} F(t)=\sum_{\substack{p,q=1\\(p\not=q)}}^{\infty}
\frac{\sin{(p-q)t}}{p-q}x_{p}x_{q}\,,
\end{equation}
is considered, where \(t, -\pi<t<\pi,\) is a parameter. It is
shown that the form \(F(t)\) is bounded and that
\begin{equation}
\label{EF} -t\sum\limits_{p=1}^{\infty}x_p^2\leq F(t) \leq
(\pi-t)\sum\limits_{p=1}^{\infty}x_p^2\,.
\end{equation}
 It is also shown that the quadratic form
\begin{equation*}
\sum_{\substack{p,q=1\\(p\not=q)}}^{\infty}\left|
\frac{\sin{(p-q)t}}{p-q}\right|x_{p}x_{q}
\end{equation*}
is unbounded for \(t\not=0\). The form (\ref{QuFo}) is interesting
because it is an example of a \textsl{symmetric} infinite matrix
\(\big[a_{pq}\big]\) that corresponds to a bounded bilinear
form, whereas the form related to the matrix
\(\big[\,|a_{pq}|\,\big]\) is unbounded.   The Hilbert matrix
\([h^-_{pq}]\), also generates a bounded bilinear form \(H^-\)
(see (\ref{HF})) and the matrix \(\big[\,|h^-_{pq}|\,\big]\)
also corresponds to an unbounded form. However, the Hilbert matrix
\(H^-\) is \textsl{antisymmetric}.

\begin{minipage}{15.0cm}
\section{\hspace{-0.5cm}.\hspace{0.2cm}%
{\large The Schur product and Schur multipliers.\label{SchurPr}}}
\end{minipage}\\[0.18cm]
\setcounter{equation}{0}

Let \(A\) and \(B\) be  matrices of the same size whose entries
are either real or complex numbers (or even belong to some ring
\(\mathfrak{R}\)): \(A=[a_{pq}], B=[b_{pq}].\) The
\textsl{Schur product} \(A\circ  B\) of the matrices \(A\) and
\(B\) is the matrix \(C=[c_{pq}]\) (of the same size as \(A\)
and \(B\)) for which \(c_{pq}=a_{pq}\cdot b_{pq}\).

The
term \textsl{Schur product} is used because  the 
 product \(A\circ B\) was introduced in (\cite{Sch4b}) for matrices, and 
some basic results about this product
  were  obtained by  Schur in that paper.  
 The most basic of these results states that  the cone
of positive semidefinite  matrices is closed under the Schur
product. We recall
 that a
square matrix \(M=[m_{pq}]\) (with complex entries) is said to
be \textsl{positive semidefinite} if  the inequality
\(\sum\limits_{p,q}m_{pq}x_q\overline{x_p}\geq 0\) holds for every
sequence \(\{x_k\}\) of complex numbers. (In the case of an
infinite matrix \(M\), only  sequences \(\{x_k\}\) with finitely
many \(x_k\) different  from zero are considered.)

\noindent %
THEOREM (The Schur product theorem, Theorem VII, \cite{Sch4b}). \textsl{If
\(A\) and \(B\) are positive semidefinite matrices (of the same
size), then their Schur product \(A\circ B\) is a positive
semidefinite matrix as well.}

For self-evident reasons, the Schur product is sometimes called
\textsl{the entrywise product} or \textsl{the elementwise
product}. It is also often referred to as the
 \textsl{Hadamard product}. The term \textsl{Hadamard product} seems to have
 appeared in print for the first time  in the 1948 (first)
edition of \cite{Hal1}. This may be due to the well known paper of
Hadamard \cite{Had}, in which he studied two Maclaurin series
\(f(z)=\sum_na_nz^n\) and \(g(z)=\sum_nb_nz^n\) with positive
radii of convergence and their composition
\(h(z)=\sum_na_nb_nz^n\), which he defined as the
\textsl{coefficientwise} product. Hadamard showed that
\(h(\cdot)\) can be obtained from \(f(\cdot)\) and \(g(\cdot)\) by
an integral convolution. He proved that any singularity \(z_1\) of
\(h(\cdot)\) must be of the form \(z_1=z_2z_3\), where \(z_2\) is
a singularity of \(f(\cdot)\) and \(z_3\) is a singularity of
\(g(\cdot)\). (This result is commonly known as \textsl{the
Hadamard composition theorem}.)
 Even though Hadamard did not study entrywise products of matrices in this 
paper, the
enduring influence of the cited result  as well as his
mathematical eminence seems to have linked   his name firmly  with
term-by-term products of all kinds, at least for analysts.
(Presentations of the Hadamard composition theorem can be found,
for example, in \cite{Bie}, Theorem 1.4.1, and in \cite{Tit},
Section 4.6.

\noindent%
PROOF of the Schur product theorem. It is enough to prove this
theorem for matrices of  arbitrary finite size. First we prove the
theorem for matrices \(A\) and \(B\) of rank one. In this case the
matrices \(A\) and \(B\) must be of the form \(A=a\cdot a^*,\,
B=b\cdot b^*\), where \(a\) and \(b\) are column vectors. It is
evident that the matrix \(C=A\circ B\) is of the form \(C=c\cdot
c^*\) where the column vector \(c\) is just the Schur product of
the column vectors \(a\) and \(b\): \(c=a\circ b\). Hence, the
matrix \(C\) is positive semidefinite. In the general case, we use
the spectral decomposition theorem. This theorem states that every
finite positive semidefinite matrix \(M\) admits a decomposition
of the form \(M=\sum\limits_{\lambda\in\sigma(M)} M(\lambda)\),
where the summation index \(\lambda\) runs over the spectrum
\(\sigma(M)\) of the matrix \(M\), and the matrices \(M(\lambda)\)
are either positive semidefinite matrices of rank one or zero
matrices. Decomposing the given matrices \(A\) and \(B\) in this
way:
\(A=\sum\limits_{\lambda\in\sigma(A)}A(\lambda),\,B=\sum\limits_{\mu\in\sigma(B)}B(\mu)\),
we see that \(\displaystyle A\circ B=\sum_{\substack{\lambda\in\sigma(A)\\
\mu\in\sigma(B)}}A(\lambda)\circ B(\mu)\) is a sum of positive
semidefinite matrices: The Schur product \(A(\lambda)\circ
B(\mu)\) of positive definite matrices of rank one is a positive
semidefinite matrix, whereas, if at least one of the matrices
\(A(\lambda)\) or \(B(\mu)\) is equal to zero,
 then
their Schur product is equal to zero. Thus, the theorem is proved.
 \hfill \framebox[0.45em]{ }


Every matrix \(H\), finite or infinite, generates a linear
operator \(\mathfrak{T}_H\) acting in the space of all matrices of
the same size as
 \(H\):
\begin{equation*}
\label{SchTr}%
\mathfrak{T}_H:A\rightarrow H\circ A,\ \ \text{or}\ \
\mathfrak{T}_HA=H\circ A.
\end{equation*}
The linear operator \(\mathfrak{T}_H\) is said to be \textsl{the
Schur transformator generated by the matrix \(H\)}. (The term
\textsl{transformator} is borrowed from \cite{GoKr}, who used it
to designate a linear operator that acts in a space of matrices
(operators).) If the Schur transformator \(\mathfrak{T}_H\) is a
bounded operator in a space of infinite matrices, equipped with a
norm, then the matrix \(H\) is said to be \textsl{the Schur
multiplier} (with respect to this norm).

The first basic estimate of the norm of the transformator
\(\mathfrak{T}_H\) was obtained by Schur in \cite{Sch4b} :

\noindent%
THEOREM (The Schur estimate for  positive definite Schur
transformators). \textsl{Let \(H=[h_{pq}]\) be a positive
semidefinite matrix for which
\begin{equation}
\label{DN}%
   D_H\stackrel{\mathrm{def}}{=} \sup\limits_{p}h_{pp}<\infty\,.
\end{equation}
Then}
\begin{equation}
\label{SME}%
 \|H\circ A\|_{l^2\to l^2}\leq D_H\|A\|_{l^2\to l^2}\,.
\end{equation}

(Here, as before, \(\|A\|_{l^2\to l^2}\) is the operator norm of
the matrix \(A\) considered in the appropriate space \(l^2\) of
sequences).

\noindent%
PROOF of the estimate (\ref{SME}). We reproduce here the reasoning
of Schur from \cite{Sch4b}. It suffices to consider only finite
matrices.
 The proof is based essentially on the fact that a
positive semidefinite matrix \(H\) admits a factorization of the
form
\begin{equation}
\label{Fact}%
                  H=LL^*,
\end{equation}
where \(L=[l_{pq}]\), i.e.,  
\begin{equation}
\label{EW}%
h_{pq}=\sum\limits_{r}l_{pr}\overline{l_{qr}}\quad
(\forall\,p,q)\,.
\end{equation}
Therefore, the number \(\sum\limits_{p,q}a_{pq}h_{pq}y_q\overline{x_p}\)
can be rewritten as
\begin{equation*}
\label{Rew} \sum\limits_{p,q}a_{pq}h_{pq}y_q\overline{x_p}=
\sum\limits_{p,q}a_{pq}
\Big(\sum\limits_{r}l_{pr}\overline{l_{qr}}\Big)y_q\overline{x_p}=
\sum\limits_{r}\sum\limits_{p,q}a_{pq}(l_{pr}\overline{x_p})\overline{(l_{qr}}y_q)\,.
\end{equation*}
Thus,
\begin{equation*}
\label{CoE}%
\begin{array}{l}
\Big|\sum\limits_{p,q}a_{pq}h_{pq}\overline{x_p}y_q\Big|\leq  %
\sum\limits_r\Big|\sum\limits_{p,q}a_{pq}(l_{pr}\overline{x_p})\overline{(l_{qr}}y_q)\Big|%
\leq \sum\limits_r\|A\|
\Big(\sum\limits_{k}\big|l_{kr}x_k\big|^2\Big)^{1/2}
\Big(\sum\limits_{k}\big|l_{kr}y_k\big|^2\Big)^{1/2}\\
\\
 =\|A\|\sum\limits_r
\big(\sum\limits_{k}\big|l_{kr}x_k\big|^2\big)^{1/2}
\Big(\sum\limits_{k}\big|l_{kr}y_k\big|^2\Big)^{1/2}\leq %
\|A\|\big(\sum\limits_{r}\sum\limits_{k}\big|l_{kr}x_k\big|^2\big)^{1/2}
\big(\sum\limits_{r}\sum\limits_{k}\big|l_{kr}y_k\big|^2\big)^{1/2}
\\
\\
\hspace{7.0ex}\leq\|A\|\Big(\sum\limits_{k}
\big(\sum\limits_{r}\big|l_{kr}\big|^2\big
)\big|x_k\big|^2\Big)^{1/2}
\Big(\sum\limits_{k}\big(\sum\limits_{r}\big|l_{kr}\big|^2\big
)\big|y_k\big|^2\Big)^{1/2}
\\
\\
\hspace{7.0ex}\leq\|A\|\Big(\sum\limits_{k}\big(\max\limits_{k}
\sum\limits_{r}\big|l_{kr}\big|^2\big
)\big|x_k\big|^2\Big)^{1/2}
\Big(\sum\limits_{k}\big(\max\limits_{k}\sum\limits_{r}\big|l_{kr}\big|^2\big
)\big|y_k\big|^2\Big)^{1/2}
\\
\\
\hspace{7.0ex}= \|A\|\,\Big(\max\limits_k\sum\limits_r|l_{kr}|^2\Big)%
\big(\sum\limits_{k}|x_k|^2\big)^{1/2}\big(\sum\limits_{k}|y_k|^2\big)^{1/2}\,.
\end{array}
\end{equation*}
According to (\ref{EW}), \(\sum\limits_{r}|l_{kr}|^2=h_{kk}\).
Thus,
\(\max\limits_{k}\Big(\sum\limits_{r}|l_{kr}|^2\Big)=\max\limits_kh_{kk}=D_H\).
Finally,
\begin{equation}
\label{FE}%
\Big|\sum\limits_{p,q}a_{pq}h_{pq}\overline{x_p}y_q\Big| \leq
\|A\|\cdot D_H\cdot
\big(\sum\limits_{k}|x_k|^2\big)^{1/2}\big(\sum\limits_{k}|y_k|^2\big)^{1/2}\,,
\end{equation}
where \(\{x_k\}\) and \(\{y_k\}\) are \textsl{arbitrary}
sequences.
 This is the estimate (\ref{SME}).  \hfill \framebox[0.45em]{ }

In fact, the reasoning of Schur allows us to prove a slightly more
general result:

\noindent%
THEOREM (The Schur factorization estimate for  Schur
transformators). \textsl{Let \(H=[h_{pq}]\) be a matrix which
admits a factorization of the form
\begin{equation}%
\label{GFF}%
H=L\cdot M^*\,,\quad\text{i.e.,}\quad
h_{pq}=\sum\limits_{r}l_{pr}\overline{m_{qr}} \ \ \  (\forall p,q),
\end{equation}%
where the matrices \(L=[l_{pr}]\) and \(M=[m_{rq}]\) satisfy
the conditions
\begin{equation}%
\label{FD}%
D_L\stackrel{\mathrm{def}}{=}\sup\limits_p\sum\limits_r|l_{pr}|^2
<\infty \ \
 and \ \ D_M\stackrel{\mathrm{def}}{=}\sup\limits_q\sum\limits_r|m_{qr}|^2
<\infty\,.
\end{equation}%
Then for every matrix \(A\) (of the same size as H) the following
inequality holds:}
\begin{equation}%
\label{SMEN}%
\|H\circ A\|_{l^2\to l^2}\leq\sqrt{D_LD_M}\,\|A\|_{l^2\to l^2}\,.
\end{equation}%

\noindent%
REMARK. The matrices \(L\), \(M\) and \(H\) need not be square.
The only restriction is that the matrix multiplication \(L,M\to
L\cdot M^*\) is feasible. In fact, the set over which  the
summation index \(r\) runs  in (\ref{GFF})  need not be a subset
of the set of integers. It can be of a much more general nature.
Thus, for example,  \textsl{let \(\mathfrak{X}\) be a measurable space
carrying a sigma-finite non-negative measure \(dx\). Let
\(\{l_{p}(x)\}\) and \(\{m_{q}(x)\}\) be sequences of
\(\mathfrak{X}\)-measurable functions defined on \(\mathfrak{X}\)
and satisfying the conditions \(D_L<\infty,\,D_M<\infty\), where
now
\begin{equation}%
\label{Now}%
D_L=\sup\limits_k\int\limits_{\mathfrak{X}}|l_{k}(x)|^2\,dx\quad
and\quad
D_M=\sup\limits_k\int\limits_{\mathfrak{X}}|m_{k}(x)|^2\,dx\,.
\end{equation}%
Let \(H\) be a matrix with entries
\begin{equation}%
\label{FNow}%
h_{pq}=\int\limits_{\mathfrak{X}}l_{p}(x)\overline{m_{q}(x)}\,dx\quad
(\forall\,p,q)\,
\end{equation}%
(i.e., the matrix \(H\) admits a factorization of the form
\(H=L\cdot M^*\), where \(L\) and \(M\) are operators acting from
the Hilbert space \(L^2(\mathfrak{X}, dx)\) into  appropriate
spaces of \(l^{\infty}\) sequences).
 Then the inequality (\ref{SMEN}) holds for an arbitrary matrix \(A\)
 (of the appropriate size), where now
\(D_L\) and \(D_M\) are defined in (\ref{Now}).}

The last result (with \(\mathfrak{X}=(a,b)\), a finite or infinite
subinterval of \(\mathbb{R}\), and  Lebesgue measure \(dx\) on
\((a,b)\)) appears as Theorem VI in \cite{Sch4b}.

The matrix %
\[H=\Big[\frac{1}{{\lambda}_p+{\mu}_q}\Big]_{1\leq
p,q<\infty}\,,\] where \({\lambda}_k\) and \({\mu}_k\) are
sequences of positive numbers that are separated from \(zero\):
\(\inf_k{\lambda}_k>0,\,\inf_k{\mu}_k>0\), serves as an example.
Here, 
\[h_{pq}=\int\limits_{0}^{\infty}e^{-{\lambda}_px}\cdot e^{-{\mu}_qx}dx\,,\,\,
i.e.,\,\,\,\,l_p(x)=e^{-{\lambda}_px}\,,m_q(x)=e^{-{\mu}_qx}\,,\quad1\leq
p,q<\infty\,,\]
 and for this \(H\) the inequality (\ref{SME}):
\[\Big\|\Big[\frac{a_{pq}}{{\lambda}_p+{\mu}_q}\Big]\Big\|_{l^2\to l^2}
\leq D_H\big\|\big[a_{pq}\big]\big\|_{l^2\to l^2}\]
 holds with
 \[D_H=\dfrac{1}{2}\cdot\dfrac{1}{\sqrt{\inf_k{\lambda}_k}}\cdot
 \dfrac{1}{\sqrt{\inf_k{\mu}_k}}.\]
 (This example is adopted from \cite{Sch4b};  it appears at the end 
of \S 4.)

 It is remarkable that the existence of a factorization of  the form
 \(H=L\cdot M^*\) for the matrix \(H\)  is not only sufficient but is also
a necessary condition for
 the operator \(A\to H\circ A\) to be a bounded operator in the space of all
matrices \(A\) (equipped  with the operator norm in \(l^2\)). This
\textsl{converse} result was
 proved by G.Bennett in \cite{Ben}.

\noindent%
THEOREM (The inversion of the Schur factorization estimate).
\textsl{Let a given matrix \(H=[h_{pq}]\) (finite or infinite)
satisfy the inequality
\begin{equation}%
\label{IFE}%
  \|H\circ A\|_{l^2\to l^2}\leq D\| A\|_{l^2\to l^2}
\end{equation}%
for all matrices \(A\) of the same size as \(H\), for some finite
constant \(D\) that does not depend on \(A\). Then for every
\(\epsilon>0\), the matrix \(H\) can be factored in the form
\(H=L\cdot M^*\), where the matrices \(L=[l_{pr}]\) and
\(M=[m_{rq}]\) act from \(l^2\) to \(l^{\infty}\) and satisfy
the inequality \(\sqrt{D_L\cdot D_M}<D+\epsilon\), and the values
\(D_L\) and \(D_M\) are defined in (\ref{FD}), i.e.,
\(D_L=\|L\|_{l^2\to l^{\infty}} \ and
 \ D_M=\|M\|_{l^2\to
l^{\infty}}\).}

This theorem appears as Theorem 6.4 in \cite{Ben}. It shows that
the Schur factorization  gives a result which is in some sense
optimal. The proof of this theorem of G.\,Bennett is essentially
based on results obtained by A.\,Pietsch on absolute summing
operators in Banach spaces, see \cite{Pie1} and \cite{Pie2} (which, in 
turn, are based on fundamental
results of A.\,Grothendieck, see the references in \cite{Pie1} and
\cite{Pie2}).

In \cite{Sch4b}, Schur considers a new class of functions of
matrices, namely, the so called Schur (or Schur-Hadamard)
functions of matrices. Let \(A=[a_{pq}]\) be an infinite matrix
whose entries have a common finite bound: \(|a_{pq}|\leq R\ \
(\forall p,q)\), where \(R<\infty\). Let \(f(\,\cdot\,)\) be a
function that is defined in the closed disk \(\{z:|z|\leq R\}\).
The matrix \(f^{\circ}(A)\) is defined ``entrywise" as follows:
\begin{equation*}
\label{SchFunc}%
f^{\circ}(A)\stackrel{\mathrm{def}}{=}[f(a_{pq})] \,.
\end{equation*}
 The following result is proved in \cite{Sch4b}:
\textsl{ Let \(f(z)=\sum\limits_{k=1}^{\infty}c_kz^k\),
 where \(
\sum\limits_{k=1}^{\infty}|c_k\,|R^k<\infty\)\,
 and let the operator generated by
the matrix \(A\) be bounded, i.e., \(\|A\|_{l^2\to l^2}<\infty\).
Then the operator generated by the matrix \(f^{\circ}(A)\) is also
bounded:
\(\big\|f^{\circ}(A)\big\|_{l^2\to l^2}<\infty\).}\\
This result appears as Theorem IV in \cite{Sch4b}.

The concept of the Schur (Schur-Hadamard) product arises in
several different areas of analysis (complex function theory,
Banach spaces, operator theory, multivariate analysis); see the
references in the introduction to \cite{Ben}. The paper \cite{Sty}
contains some applications of the Schur product to multivariate
analysis as well as a rich bibliography of books and articles
related to Schur-Hadamard products. The paper \cite{HorR1}
contains a lot of facts about Schur-Hadamard products and
Schur-Hadamard functions of matrices as well as a rich
bibliography. In particular,
 it discusses  fractional Schur-Hadamard
powers of a positive matrix, infinite Schur-Hadamard divisibility
of a positive matrix and its relation to the conditional
positivity of the \(\text{logarithmic}^{\circ}\) matrix. Chapter 5
of the book \cite{HorR2} (about eighty pages) is dedicated to the
Schur-Hadamard product of matrices.

 A very fruitful generalization of
the Schur transformator is the  \textsl{the Stieltjes
double-integral operator}. This notion seems to have
 appeared first in the papers of Yu.L.\,Daletskii and S.G.\,Krein \cite{DaKr1},
\cite{DaKr2}, \cite{Da1}, \cite{Da2}. Later on, the theory of
double-integral operators
 was elaborated on in great detail by M.S.\,Birman and M.Z.\,Solomyak in
\cite{BiSo1} -- \cite{BiSo4}.

Let \(\Lambda\) and \(\mathrm{M}\) be measurable spaces, i.e.,
sets provided with sigma-algebras of subsets, and let
\(E(d\lambda)\) and \(F(d\mu)\) be two \textsl{orthogonal measures
in a separable Hilbert space \(\mathfrak{H}\)} that are  defined on
\(\Lambda\) and \(\mathrm{M}\), respectively, i.e.,
weakly-countably-additive functions taking their values in the set
of orthogonal projectors in \(\mathfrak{H}\) and satisfying the
condition \(E(\alpha)E(\beta)=0\) if \(\alpha\cap\beta=\emptyset\)
and \(F(\gamma)F(\delta)=0\) if \(\gamma\cap\delta=\emptyset\). We
assume also that the orthogonal  measures \(E(d\lambda)\) and
\(F(d\mu)\) are \textsl{spectral measures}, i.e., they also
satisfy the conditions \(E(\Lambda)=I \ \mbox{and} \  F(\mathrm{M})=I\), where
\(I\) is the identity operator in \(\mathfrak{H}\). If \(A\) is a
bounded linear operator in \(\mathfrak{H}\), then
\begin{equation}
\label{DOI}%
A=\iint\limits_{\mathrm{M}\times\Lambda}F(d\mu)AE(d\lambda)\,,
\end{equation}
where the integral can be understood in any reasonable sense. The
equality (\ref{DOI}) can be considered as a direct generalization
of the matrix representation of an operator in a Hilbert space
with respect to two orthonormal bases. Namely, let the orthogonal
spectral measures \(E(d\lambda)\) and \(F(d\mu)\) be discrete and
let their ``atoms" be one-dimensional orthogonal projectors, i.e.,
the atom of the measure \(E(d\lambda)\), located at the point
\(\lambda\in\Lambda\), is of the form \(E(\{\lambda\})=\langle
\,\,\cdot\,, e_{\lambda}\rangle e_{\lambda}\)  and the atom of the
measure \(F(d\mu)\), located at the point \(\mu\in\mathrm{M}\), is
of the form \(F(\{\mu\})=\langle \,\,\cdot\,, f_{\mu}\rangle
f_{\mu}\), where \(e_{\lambda}\) and \(f_{\mu}\) are normalized
vectors generating the one-dimensional subspaces
\(E(\{\lambda\})\mathfrak{H}\) and \(F(\{\mu\})\mathfrak{H}\),
respectively. The collection of all the vectors \(\{e_{\lambda}\}\)  
corresponding to all the atoms of the measure \(E(d\lambda)\)
forms an orthonormal basis of the space \(\mathfrak{H}\).
Analogously,  the collection of all the vectors \(\{f_{\mu}\}\) 
corresponding to all the atoms of the measure \(F(d\mu)\) also
forms an orthonormal basis of the space \(\mathfrak{H}\). Consequently, 
the representation (\ref{DOI}) of the operator
\(A\) takes the form
\begin{equation}%
\label{DOID}%
A=\sum\limits_{\lambda,\mu} f_{\mu}\,a_{\mu,\lambda}\,
\langle\,\,\cdot\,,e_{\lambda}\,\rangle\,,
\end{equation}%
where \(a_{\mu,\lambda}=\langle Ae_{\lambda}\,,f_{\mu}\rangle\,.\)
 Thus, in the  case of discrete orthogonal spectral measures
with one-dimensional atoms, the representation (\ref{DOI}) turns
into the matrix representation of a given operator with respect to
given orthonormal bases. The matrix \([a_{\mu,\lambda}]\)
corresponds to the operator \(A\). If \(h(\mu, \lambda)\) is a
measurable function defined on \(\mathrm{M}\times \Lambda\), then
the sum
\begin{equation}%
\label{DST}%
\mathfrak{T}_hA\stackrel{\mathrm{def}}{=}\sum\limits_{\lambda,\mu}\,
f_{\mu}\,h_{\mu,\lambda}\cdot a_{\mu,\lambda}\, \langle\,\cdot\,,
e_{\lambda}\,\rangle\,
\end{equation}%
can be pictured as an application of the Schur transformator
corresponding  to the matrix \([h_{\mu,\lambda}]\) to the
operator \(A\): \(A\mapsto \mathfrak{T}_h A\). The sum on the
right hand side of the equality (\ref{DST}) can be formally
written as an integral:
\begin{equation}%
\label{SDOI}%
\mathfrak{T}_hA=\iint\limits_{\mathrm{M}\times\Lambda}h%
(\mu,\lambda)F(d\mu)AE(d\lambda)\,.
\end{equation}
However, one can consider integrals of the form (\ref{SDOI}) for
arbitrary orthogonal spectral measures \(E(d\lambda)\) on
\(\Lambda\) and \(F(d\mu)\)  on \(\mathrm{M}\), and  more or less
arbitrary functions  \(h(\mu,\lambda)\) on \(\mathrm{M}\times
\Lambda\). If the integral (\ref{SDOI}) exists in a reasonable
sense (either as a Lebesgue integral, or  a Riemann-Stieltjes
integral, or  some other integral), it is said to be \textsl{a
Stieltes double-integral operator}. The problem of establishing
the existence of a Stieltes double-integral operator is intimately
associated with estimates for it in various norms. In particular,
the estimates
\begin{equation}%
\label{UN}%
\big\|\mathfrak{T}_hA\big\|_{\mathfrak{R}\to\mathfrak{R}} \leq
C\big\|A\big\|_{\mathfrak{R}\to\mathfrak{R}}
\end{equation}%
and
\begin{equation}%
\label{NN}%
\big\|\mathfrak{T}_hA\big\|_{\mathfrak{S}_1\to\mathfrak{S}_1} \leq
C\big\|A\big\|_{\mathfrak{S}_1\to\mathfrak{S}_1}
\end{equation}%
are extremely important. Here
\(\|\Phi\|_{\mathfrak{R}\to\mathfrak{R}}\) is the ``uniform" norm
of the operator \(\Phi\), acting in \(\mathfrak{H}\):
\(\|\Phi\|_{\mathfrak{R}\to\mathfrak{R}} = \sup_{v\in\mathfrak{H},
v\neq 0}\dfrac{\|\Phi v\|_{\mathfrak{H}}}{\|v\|_{\mathfrak{H}}}\),
and \(\|\Phi\|_{\mathfrak{S}_1\to\mathfrak{S}_1}\) is its ``trace"
norm.

In \cite{BiSo4} the estimate (\ref{UN}) was obtained for functions
\(h(\,\cdot\,,\,\cdot \,)\) which admit a ``factorization" of the
form
\begin{equation}
\label{HSF}%
h(\mu, \lambda)=\int\limits_{\mathfrak{X}}m(\mu,x)\cdot
l(\lambda,x)dx\,,
\end{equation}
where \(\mathfrak{X}\) is a measurable space carrying a
non-negative sigma-finite measure \(dx\),
\begin{equation}
\label{KF}%
C_m=\mathrm{ess}\hspace{-1.7ex}\sup\limits_{\mu\in
\mathrm{M}\hspace{3.1ex}}\hspace{-1.5ex}%
\int\limits_{\mathfrak{X}}\big|m(\mu,x)\big|^2dx\,,
\quad C_l=\mathrm{ess}\hspace{-1.5ex}\sup\limits_{\lambda\in
\Lambda\hspace{3.3ex}}\hspace{-1.5ex}%
\int\limits_{\mathfrak{X}}\big|l(\lambda,x)\big|^2dx\,
\end{equation}
and
\begin{equation}%
\label{SFaE}%
C=\sqrt{C_m\cdot C_l} < \infty\,.
\end{equation}%
The inequality (\ref{UN}) is then obtained (with the same constant
$C$) by invoking
 the duality between the set \(\mathfrak{R}\) of all
bounded operators in \(\mathfrak{H}\) and the set
\(\mathfrak{S}_1\) of all trace class operators. The estimate
(\ref{NN}) holds with the same constant \(C\) (that is given in
(\ref{SFaE})). Unfortunately, the paper \cite{BiSo4} is not
translated into English, but some results of this paper, in
particular, the estimate (\ref{UN}), (\ref{SFaE}), are reproduced
in \cite{ABF}, Section 2.

The estimate (\ref{UN}) is a direct analog of the Schur
factorization estimate (\ref{SMEN}), (\ref{FD}) and is obtained by
the same method that Schur used. However, when Birman and Solomyak
started to develop the theory of  Stieltjes double-integral
operators, they were not aware of the paper \cite{Sch4b} by Schur.
The close relationship between double-integral operators and the
results of Schur was only discovered later. In Section 2 of  
\cite{Pel}, V.\,Peller 
obtained a result that ``inverts" the estimate
(\ref{UN}) by Birman and Solomyak in the same sense that Theorem 6.4 of 
 \cite{Ben} (that was stated earlier) 
inverts the factorization estimate by Schur. Peller proved an even
stronger result, a ``maximal" version of the inverse result.
Namely, he proved that if the function \(h\) is such that the
estimate (\ref{UN}) holds for every bounded operator \(A\) in
\(\mathfrak{H}\) with a finite constant \(C\) that is independent
of \(A\), then the function \(h(\,\cdot\,,\,\cdot\,)\) admits a
factorization of the form (\ref{HSF}), where the functions
\(m(\,\cdot,\,\cdot\,)\) and \((\,\cdot,\,\cdot\,)\) satisfy the
conditions
\begin{equation}%
\int\limits_{\mathfrak{X}}\Big(
\mathrm{ess}\hspace{-1.7ex}\sup\limits_{\mu\in
\mathrm{M}\hspace{3.1ex}}\hspace{-1.5ex}\big|m(\mu,x)\big|\Big)^2dx<\infty
\quad  \mbox{and} \quad \int\limits_{\mathfrak{X}}\Big(
\mathrm{ess}\hspace{-1.7ex}\sup\limits_{\lambda\in
\Lambda\hspace{3.1ex}}\hspace{-1.5ex}\big|l(\lambda,x)\big|\Big)^2dx<\infty.
\label{MC}%
\end{equation}%
The estimate (\ref{UN}), (\ref{SFaE}) is ``semi-effective": given
the function \(h(\mu,\lambda)\), it is not so easy to see when it
admits a factoraization of the form (\ref{HSF}). To overcome this
difficulty, Birman and Solomyak developed another approach that
reduces the study of
 Stieltes double-integral operators to the study of integral operators of the
form %
\begin{equation}%
\label{IntOp}%
 u(\lambda)\to v(\mu)=\int_{\Lambda}h(\mu,\lambda)u(\lambda)\rho(d\lambda)\,.
\end{equation}%
This reduction is explained in \cite{BiSo2}, Theorem 2, and also
in \cite{BiSo4}, Lemma 1.1. The operator (\ref{IntOp}) acts from
the space \(L^2(\Lambda,d\rho(\lambda))\) into the space
\(L^2(\mathrm{M},d\sigma(\mu))\), where
 \(\rho(d\lambda)=\langle
E(d\lambda)\omega,\omega\rangle,\ \sigma(d\mu)=\langle
F(d\mu)\theta,\theta\rangle\) and $\omega, \theta \in
\mathfrak{H}$\,.
 The estimates for the integral
operators (\ref{IntOp}) must be carried out for all vectors
\(\omega, \theta\in\mathfrak{H}\) and must be uniform with respect
to the measures \(\rho(d\lambda)\) and \(\sigma(d\mu)\). To obtain
such estimates, Birman and Solomyak developed a method that is
based on the approximation of functions from the
Sobolev-Slobodetski\u{\i} classes \(W_p^{\alpha}\) by
piecewise-polynomial functions, \cite{BiSo5}, \cite{BiSo6},
\cite{BiSo7}, \S\S\,8 - 9, \cite{BiSo8}, Chapter 3,\,\S\S\,5 - 7.
In the construction of the approximating functions, a partition of
the domain of definition of the approximated function appears. To
achieve the desired uniformity of the approximation with respect
to the measures \(\rho(d\lambda)\) and \(\sigma(d\mu)\), this
partition must be adapted to these measures.

The approach, based on  piecewise-polynomial approximations,
allows one to approximate  the kernels of the integral operators
(\ref{IntOp}) by finite-dimensional kernels, and thus to obtain
the needed estimates for the singular values of the Stieltjes
double-integral operators. The estimates of the double-integral
operators are made not only in the uniform and trace norms, but
also in many other norms. These estimates depend upon 
the smoothness of the function \(h(\,\cdot\,,\,\cdot\,)\)
(assuming that \(\Lambda\) and \(\mathrm{M}\) are smooth
manifolds).

Double-integral operators appear in the formula for
differentiating functions of Hermitian operators  with respect to
a parameter. Namely, let \(\tau\to H(\tau)\) be a function on some
open subinterval of the real axis \(\mathbb{R}\) whose values are
self-adjoint operators in a Hilbert space \(\mathfrak{H}\). Let
\(f: \mathbb{R}\to \mathbb{R}\) be a real-valued function that is defined
and bounded on \(\mathbb{R}\) and let \(E(d\lambda,\tau)\) be the
spectral measure of the operator \(H(\tau)\). Under appropriate 
assumptions, Yu.L.\,Daletskii and S.G.\,Krein, \cite{DaKr1},
obtained the formula
\begin{equation}%
\label{DaKr}%
\frac{\partial
f(H(\tau))}{\partial\tau}=\iint\limits_{\mathbb{R}\times\mathbb{R}}
\frac{f(\lambda)-f(\mu)}{\lambda-\mu}\,E(d\mu,\tau)\,\frac{\partial
H(\tau)}{\partial \tau}\,E(d\lambda,\tau)\,.
\end{equation}%
This formula, which expresses  the derivative \ \(\dfrac{\partial
f(H(\tau))}{\partial\tau}\) \  as a Stieltjes double-integral
operator, seems to be the first recorded application of Stieltjes
double-integral operators. The paper \cite{Da1} contains a version
of  Taylor's formula for operator functions. The paper
\cite{DaKr2} (and, to some extent, the paper \cite{Da2}) contains a
more detailed presentation of the results of the papers
\cite{DaKr1} and \cite{Da1} as well as some  extensions. Later on,
Stieltjes double-integral operators were widely used in
scattering theory.
 M.Sh.\,Birman, \cite{Bi1}, used them to prove the existence of wave operators.
 ( See also \cite{BiSo2}, especially the last paragraph of this paper.)
 Double-integral operators are involved in the study of the so called
 \textsl{spectral shift function} (see \cite{BiSo10} and \cite{BiYa}).
 The paper \cite{BiSo11} is devoted to the application of
double-integral operators  to the estimation of perturbations and
commutators  of functions
 of self-adjoint operators.
 It is worth  noticing that double-integral operators allow one to make
 an abstract and symmetric definition of a pseudodifferential operator
 with prescribed symbol (see item \textbf{3} of the paper \cite{BiSo9}).

Thus, the ideas of Issai Schur on the termwise multiplication of
matrices, partially forgotten and rediscovered, are seen to lead
very far from the original setting.

\setcounter{equation}{0}
\begin{minipage}{15.0cm}
\section{\hspace{-0.5cm}.\hspace{0.2cm}%
{\large The Schur Convexity Theorem. \label{SChurConvTheor}}}
\end{minipage}\\[0.18cm]

The well known \textsl{Hadamard inequality} states that
\begin{equation}%
\label{HadIn}%
\text{det}H\leq\prod\limits_{1\leq k \leq n } h_{kk}
\end{equation}%
for every non-negative definite Hermitian matrix
\(H=[h_{jk}]_{1\leq j,k \leq n}\).  (There are many proofs; see,
for example, \cite{HoJo}, Section 7.8.) In a short but penetrating
paper published in 1923, Issai Schur \cite{Sch18a} gave a highly
effective method for deriving this inequality. However the
importance of the paper \cite{Sch18a} rests primarily on the ideas
which are contained there and by the impact which the paper had on
various areas of mathematics, some of which lie very far from the
original setting. This paper has generated and continues to
generate many fruitful investigations.

Given a Hermitian matrix \(H=[h_{jk}]_{1 \leq j,k \leq n}\), it
can be reduced to the
diagonal form %
\begin{equation}%
\label{DiagF}%
 H=U\,\text{diag}(\omega_1,\,\dots\,,\omega_n)\,U^*,
\end{equation}%
 where
\(\omega_1,\,\dots\,,\omega_n\) are the eigenvalues of the matrix
\(H\), and \(U=[u_{jk}]_{1\leq j,k \leq n}\) is a unitary matrix.
(If the Hermitian matrix \(H\) is real,  then the matrix \(U\) can
be chosen real also, i.e., if $H$ is real and symmetric, then $U$
is orthogonal.) In particular, the equality (\ref{DiagF}) implies
that
\begin{equation}%
\label{OrtAv}%
\left[%
\begin{array}{c}%
h_{11}\\ \vdots\\
 h_{nn}
\end{array}%
\right]%
= \left[%
\begin{array}{ccc}%
|u_{11}|^2&\dots&|u_{1n}|^2\\[3pt]
\dots&\dots&\dots\\[3pt]
|u_{n1}|^2&\dots&|u_{nn}|^2\\
\end{array}%
\right]%
\left[%
\begin{array}{c}%
\omega_1\\
\vdots\\
\omega_n
\end{array}%
\right].%
\end{equation}%
Since the matrix \(U\) in (\ref{DiagF}) is unitary (orthogonal),
the matrix \(M=[m_{jk}]_{1\leq j,k\leq n}\), with
\begin{equation}%
\label{OSDS}%
\ m_{jk}=|u_{jk}|^2,
\end{equation}%
as in (\ref{OrtAv}), possesses the properties
\begin{equation}%
\label{Doub-Stoch}%
\begin{array}{rll}%
\text{i.}  &m_{jk}\geq 0,\quad &1\leq j,k \leq n;\\[5pt]
\text{ii.} &\sum\limits_{1\leq k\leq n}m_{jk}=1, \quad& 1\leq j
\leq n;\\[13pt]
\text{iii.} & \sum\limits_{1\leq j \leq n}m_{jk}=1, \quad& 1\leq
k\leq n\,.
\end{array}%
\end{equation}%

It turns out to be fruitful to consider linear transformations
whose matrices \(M\) satisfy the conditions (\ref{Doub-Stoch}),
without regard to the relations (\ref{OSDS}).

\noindent%
\textsf{DEFINITION 1.}
\begin{slshape}
 A matrix \(M=[m_{jk}]_{1\leq j,k \leq n}\) is said to be
\textsf{doubly-stochastic} if the conditions
\textup{(\ref{Doub-Stoch})} are fulfilled.
\end{slshape}

\noindent%
\textsf{DEFINITION 2}
\begin{slshape}
A matrix \(M=[m_{jk}]_{1\leq j,k \leq n}\) is said to be
\textsf{ortho-stochastic} if there exists an orthogonal matrix
\(U=[u_{jk}]_{1\leq j,k \leq n}\) such that the matrix entries
\(m_{jk}\) are representable in the form (\ref{OSDS}), i.e., if
$M$ is the Schur product of an orthogonal matrix $U$ with itself.
\end{slshape}%

\noindent%
\textsf{REMARK 1.}
\begin{slshape}
It is clear that every ortho-stochastic matrix is a
doubly-stochastic. However, not every doubly-stochastic matrix is
an ortho-stochastic. For example\,\footnote{This
example is adopted from \textup{\cite{Sch18a}}.}, the matrix %
\(P=\dfrac{1}{6}\left[\begin{array}{ccc}0&3&3\\3&1&2\\3&2&1\end{array}\right]\) %
is doubly-stochastic, but not ortho-stochastic.
\end{slshape}%

Many well known elementary inequalities can be put in the form
\begin{equation}%
\label{SSI}
\Phi(\overline{x},\,\dots\,,\overline{x})\leq\Phi(x_1,\,\dots\,,x_n),%
\end{equation}%
where \(\overline{x}=(x_1+\,\cdots\,+x_n)\)/n and
\(x_1,\,\dots\,x_n\) lie in a
specified set. For example, the  inequality  %
\begin{equation}%
\label{SJI}%
\varphi(\overline{x})\leq\big(\varphi(x_1)+\,\cdots\,+\varphi(x_n)\big)/n
\end{equation}%
for a \textsl{convex} function \(\varphi\) of one variable can be
written in the form (\ref{SSI}), with
\(\Phi(\xi_1,\,\dots\,,\xi_n)=\varphi(\xi_1)+\,\cdots\,+\varphi(\xi_n).\)

We recall, that \textsl{a real valued function \(\varphi\), defined on 
a subinterval %
\((\alpha,\beta)\) of the real axis, is 
said to be \textsf{convex} if \(\varphi\) is continuous there and the inequality %
\(\varphi\big((x_1+x_2)/2\big)\leq\big(\varphi(x_1)+\varphi(x_2)\big)\)/2
holds for every \(x_1,x_2\in(\alpha,\beta)\).} The inequality
(\ref{SJI}) is a special case of the so-called

\textsf{JENSEN INEQUALITY}. %
\begin{slshape}
Let \(\varphi\) be a convex function on an interval
\((\alpha,\beta)\), let \(x_1,\,\dots\,,x_n\) be points in the
interval \((\alpha,\beta)\), and let the numbers
\(\lambda_1,\,\dots\,,\lambda_n\) satisfy the conditions
\begin{equation}%
\label{Jens-Weigh}%
\begin{array}{rll}%
\textup{i.}  &\ \ \lambda_{k}\geq 0,\quad &1\leq k \leq n;\\[5pt]
\textup{ii.} &\sum\limits_{1\leq k\leq n}\lambda_{k}=1\,.
\end{array}%
\end{equation}%
Then
\begin{equation}%
\label{GJI}%
\varphi(\lambda_1x_1+\,\cdots\, +\lambda_nx_n)\leq\lambda_1\varphi(x_1)
+\,\cdots\,+\lambda_n\varphi(x_n)\,.
\end{equation}%
\end{slshape}

 The value \(\overline{x}=\,(x_1+\,\dots\,+x_n)\)/n that  appears
in (\ref{SJI}), the so called \textsl{arithmetic mean} of the
values \(x_1,\,\dots\,,x_n\), is the most commonly used average
value for \(x_1,\,\dots\,,x_n.\) The value
\(\lambda_1x_1+\,\dots\,+\lambda_nx_n\) that appears in
(\ref{GJI}), the so-called \textsl{weighted arithmetic mean}, is a
more general average value for
\(x_1,\,\dots\,,x_n\,.\)\\
\hspace*{2.0ex}In \cite{Sch18a}, doubly-stochastic matrices
\(M=[m_{jk}]_{1\leq j,k\leq n}\) are used to construct an
\textsl{average sequence} \(y_1,\,\dots\,,y_n\) from a given
sequence of real or complex numbers  \(x_1,\,\dots\,,x_n\) by the
\textsl{averaging rule}
\begin{equation}%
\label{GenAv}%
{\bf y} =
\left[%
\begin{array}{c}%
y_{1}\\ \vdots\\
 y_{n}
\end{array}%
\right]%
= \left[%
\begin{array}{ccc}%
m_{11}&\dots&m_{1n}\\[3pt]
\dots&\dots&\dots\\[3pt]
m_{n1}&\dots&m_{nn}\\
\end{array}%
\right]%
\left[%
\begin{array}{c}%
x_1\\
\vdots\\
x_n
\end{array}%
\right]%
= M{\bf x}.
\end{equation}%
It is intuitively clear that the sequence of ``averaged" values
\(\{y_k\}\) is ``less spread out" than the original sequence
\(\{x_k\}\). In \cite{Sch18a}, inequalities of the form
\begin{equation}%
\label{GSI}%
\Phi(y_1,\,\dots\,,y_n)\leq\Phi(x_1,\,\dots\,,x_n),
\end{equation}%
are considered for points \((x_1,\,\dots\,,x_n)\) and
\((y_1,\,\dots\,,y_n)\) in the domain of definition of the
function \( \Phi \) that are related by a doubly stochastic matrix
\(M=[m_{jk}]_{1\leq j,k \leq n}\) by  means of the averaging
procedure {\bf y}=M{\bf x} given in (\ref{GenAv}). 
In particular, the inequality (\ref{GSI}) is
established there for functions
 \(\Phi\) of the
form
\(\Phi(\xi_1,\,\dots\,,\xi_n)=\varphi(\xi_1)+\,\cdots\,+\varphi(\xi_n)\):

\noindent \textsf{THEOREM I}. \begin{slshape}%
Let \(\varphi\) be a convex function defined on a
subinterval \((\alpha,\beta)\) of the real axis, %
let \(x_1,\,\dots\,,x_n\) be arbitrary numbers from
\((\alpha,\beta)\), let \(M=\big[m_{jk}\big]_{1\leq j,k\leq n}\)
be a doubly stochastic matrix and let the numbers
\(y_1,\,\dots\,,y_n\) be obtained from the averaging procedure
${\bf y}=M{\bf x}$. Then
\begin{equation}%
\label{SGSI}%
\varphi(y_1)+\,\cdots\,+\varphi(y_n)\leq\varphi(x_1)+
\,\cdots\,+\varphi(x_n)\,.%
\end{equation}%
\end{slshape}

\noindent%
PROOF of Theorem I. In view of the conditions (\ref{Doub-Stoch}.i)
and (\ref{Doub-Stoch}.ii), Jensen's inequality is applicable with
\(\lambda_k=m_{jk},k=1,\,\dots\,,n\), and implies that 
\[
m_{j1}\varphi(x_1)+\,\cdots\,+m_{jn}\varphi(x_n)\leq%
\varphi(m_{j1}x_1+\,\cdots\,+m_{jn}x_n) =  %
\varphi(y_j).
\]
The desired conclusion is now obtained by summing the last
inequality over \(j\) from  \(1,\,\dots\,,n\) and invoking the
condition (\ref{Doub-Stoch}.iii). \hspace*{1pt}\hfill
\framebox[0.45em]{}

The preceding theorem appears as Theorem V in \cite{Sch18a} and is
used there to  derive the (Hadamard) inequality
\[\prod\limits_{1\leq k\leq n}{\omega}_k
\leq\prod\limits_{1\leq k\leq n}h_{kk}\]%
for a positive definite Hermitian matrix
\(H=\big[h_{jk}\big]_{1\leq j,k\leq n}\) with eigenvalues
\({\omega}_1,\,\dots\,,{\omega}_n\). The latter is equivalent to
the inequality
\begin{equation}%
\label{HadInS}%
\sum\limits_{1\leq k\leq n}(-\log{{\,h_{kk})\leq\sum\limits_{1\leq
k\leq n}(-\log{\,\omega}_k}}),
\end{equation}%
which is of the form (\ref{SGSI}), with the convex function
\(\varphi(\xi)=-\log{\xi}\).  In this case, the averaging
doubly-stochastic matrix \(M=\big[m_{jk}\big]_{1\leq j,k\leq n}\)
is the ortho-stochastic one, with entries \(m_{jk}\) of the form
(\ref{OSDS}), as in (\ref{OrtAv}).

In \cite{Sch18a}, functions \(\Phi\) of several variables  for
which inequalities of the form (\ref{GSI}) hold are also
considered.

\noindent%
\textsf{DEFINITION 3.}
\begin{slshape}
 A function \(\Phi\)  of \(n\) variables
\(x_1,\,\dots\,,x_n.\)  is said to be \textsf{S-convex} (i.e.,
convex in the sense of  Schur) if for every doubly-stochastic
matrix \(M\) and every pair of points \({\bf
x}=(x_1,\,\dots\,,x_n)\) and ${\bf y} = M{\bf x}$ in the domain of
\(\Phi\), the inequality \textup{(\ref{GSI})} holds. The function
\(\Phi\) is said to be  \textsf{S-concave} if the opposite
inequality holds, i.e., if
\(\Phi(x_1,\,\dots\,,x_n)\leq\Phi(y_1,\,\dots\,,y_n)\), holds for
every pair of points \({\bf x}\) and \({\bf y}=M{\bf x}\) in the
domain of  \(\Phi\). A function \(\Phi\) is \(S\)-concave if and
only if the function \(-\Phi\) is \(S\)-convex.
\end{slshape}

Let \(\pi\) be a permutation of the set \(\{1,\,\dots\,,n\}\).
Then the corresponding operator on \(\mathbb{R}^n\) that permutes
coordinates according to the rule
\((x_1,\,\dots\,,x_n)\to(x_{\pi(1)},\,\dots\,,x_{\pi(n)})\) is
linear. Its matrix \(P_{\pi}\) with respect to the standard basis
in \(\mathbb{R}^n\) is  termed a permutation matrix and is of the
form
\begin{equation}%
\label{PerM}%
P_{\pi}=\big[(p_{\pi})_{jk}\big]_{1\leq j,k\leq n}\,,\ \
\text{where, for}\
k=1,\,\dots\,,n\,,\ (p_{\pi})_{jk}= %
\left\{%
\begin{array}{l}%
1,\ \  \text{if}\ j=\pi(k);\\ %
0,\ \  \text{if}\ j\not=\pi(k).
\end{array}
\right.
\end{equation}%
There are \(n!\) permutation matrices of size \(n\times n\). Every
permutation matrix is a doubly-stochastic one. The inverse of a
permutation matrix is a permutation matrix as well, and hence it
is also doubly-stochastic. Therefore,

\textsl{Every S-convex function \(\Phi\) of \(n\) variables is a
symmetric function:}
\begin{equation}%
\label{SymCond}%
\Phi(x_1,\,\dots\,,x_n)\equiv
\Phi(x_{\pi(1)},\,\dots\,,x_{\pi(n)}), \ \text{for every
permutation}\ \pi\,.
\end{equation}%

\noindent \textsf{THEOREM II}. \textsl{Let \(\Phi\) be a
\(S\)-convex function of \(n\) variables, \(n\geq 2\),  and let
all its partial derivatives of the first order exist and be
continuous. Then the function \(\Phi\) satisfies the condition}
\begin{equation}%
\label{NMC}%
\frac{\partial \Phi}{\partial x_1}\,(x_1,x_2,\,\dots\,,x_n)-%
 \frac{\partial \Phi}{\partial x_2}\,(x_1,x_2,\,\dots\,,x_n)\,
 \geq 0,\ \ \text{if}\ \ x_1>x_2\,.
\end{equation}%

This theorem provides a necessary condition for a symmetric
function \(\Phi\) be \(S\)-convex. It appears as Theorem I in
\cite{Sch18a}. Theorem II in  \cite{Sch18a} also contains a
sufficient condition for a symmetric function \(\Phi\) be
\(S\)-convex.

\noindent \textsf{THEOREM III}. \textsl{Let \(\Phi\) be a
symmetric function of \(n\) variables, \(n\geq 2\), that
satisfies the condition
\begin{equation}%
\label{SMC}%
\left(\frac{{\partial}^2\Phi}{\partial x_1^2}
+\frac{{\partial}^2\Phi}{\partial x_2^2}- 2
\frac{{\partial}^{\,2}\,\Phi}{\partial x_1\partial
x_2}\right)(x_1,x_2,\,\dots\,,x_n)\geq 0\ \ \ \text{for all}\ \ \
x_1,x_2,\,\dots\,,x_n\,.
\end{equation}%
Then the function \(\Phi\) is \(S\)-convex.}

However,  A.\,Ostrowski showed that condition (\ref{NMC})  is both
necessary and \textsl{sufficient} for a symmetric function
\(\Phi\) to be \(S\)-convex; see Theorem VIII in \cite{Ostr}. The
reasoning in \cite{Ostr} is based essentially on the  the
reasoning in \cite{Sch18a}, but is more precise.

In \cite{Sch18a} it is shown that \textsl{the elementary symmetric
functions} \(c_k(x_1,\,\dots\,,x_n),\ \ k=1,\,\dots\,,n,\) are
 \(S\)-concave, and that the functions %
\(\Phi_k(x_1,\,\dots\,,x_n)=
\dfrac{c_{k+1}(x_1,\,\dots\,,x_n)}{c_k(x_1,\,\dots\,,x_n)}\,,\)
\(k=1,\,\dots\,,n-1\,,\) are \(S\)-concave.

To this point, we have reviewed
 almost all the main results of the short paper \cite{Sch18a}.
The significance of this paper is not confined to these results,
important as they are, but rests primarily  on the fact that
linear transformations with doubly-stochastic matrices were
introduced there. This paper attracted the attention of
mathematicians to doubly-stochastic matrices. (In \cite{BeBe} the
term ``Schur transformation" is used for linear transformations
with such matrices; see \cite{BeBe}, Chapter I, \S\, 29.) Schur
himself did not use the term doubly-stochastic matrix. He just
referred to ``a matrix \(M\) that satisfies the conditions
(\ref{Doub-Stoch})." The term ``doubly-stochastic matrix" seems to
have appeared first in the first edition of the book  \cite{Fel}
by W.\,Feller, in 1950\,\footnotemark. \footnotetext{\,However,
the term ``stochastic matrix" was used as early as 1931 in
\cite{Rom1} (see also \cite{Rom2}) for matrices satisfying the
conditions (\ref{Doub-Stoch}.i) and (\ref{Doub-Stoch}.ii) only
(but not necessarily the condition (\ref{Doub-Stoch}.iii)). Such
matrices play a crucial role in the theory of Markov chains.}

Many results  were influenced by the paper \cite{Sch18a}. We
shall begin with  the theorems of Hardy-Littlewood-Polya  and
Birkhoff.

To formulate the Hardy-Littlewood-Polya Theorem, we have to
introduce the notion of majorization. Let
\(\xi_1,\,\xi_2,\,\,\dots\,\,,\xi_n\) be a sequence of
\textsl{real} numbers. By
\(\xi_1^*,\,\xi_2^*,\,\,\dots\,\,,\,\xi_n^*\) we denote the
reaarangement of this sequence in non-increasing order:
\[\xi_1^*\geq\xi_2^*\geq\,\,\dots\,\,\geq\xi_n^*,\ \ \xi_k^*=\xi_{\pi(k)}
\ \ \text{for some permutation}\ \pi\ \text{of the  set of
indices}\ 1,\,2\,\,\dots\,,n\,.\]

\noindent%
\textsf{DEFINITION 4.}
\begin{slshape}
 Let \(\boldsymbol{x}=(x_1,\,x_2,\ \dots\ ,x_n)\) and
\(\boldsymbol{y}=(y_1,\,y_2,\ \dots\ ,y_n)\) be two sequences of
real numbers. Then we say that the sequence {\bf y} is majorized
by the sequence {\bf x} (or that the sequence {\bf x} majorizes
the sequence {\bf y}) 
 if the following conditions are
satisfied:
\begin{equation}%
\label{DefSub}%
\begin{array}{llll}%
y_1^*+y_2^*+\,\dots\,+y_k^*&\leq
&x_1^*+x_2^*\,+\,\dots\,+x_k^*,&\hspace*{-3.0ex}
(k=1,\,2,\,\dots\,,n-1\,)\,;\\[2.0ex]
y_1^*+y_2^*+\,\dots\,+y_{n-1}^*+y_n^*&=
&x_1^*+x_2^*\,+\,\dots\,+x_{n-1}^*+x_n^*\,.&
\end{array}
\end{equation}%
A relation of the form \textup{(\ref{DefSub})} is said to be
\textsf{a majorization} relation and is denoted by the symbol
\begin{equation}%
\label{SubOrd}%
 \boldsymbol{y}\prec\boldsymbol{x},\ \ \text{or}\ \
 (y_1,\,y_2,\ \dots\ ,y_n) \prec(x_1,\,x_2,\ \dots\ ,x_n),
\end{equation}%
\end{slshape}

The relations (\ref{DefSub}) were considered by R.F.\,Muirhed
\cite{Muir} and by M.O.\,Lorenz \cite{Lor} in the beginning of
20th century. Muirhead introduced these relations (with integer
\(x_k,\,y_k\) only)  to study inequalities for homogeneous
symmetric functions (Muirhead's result is also presented in
\cite{HLPb}, Chapter II, sec. \textbf{2.18}). Lorenz used the
relations (\ref{DefSub}) to describe the non-uniformity of the
distribution of wealth in a population. However, the notation
(\ref{SubOrd}) and the term ``majorization" were introduced by
G.H.\,Hardy, J.W.\,Littlewood and G.\,Polya in 1934; see
\cite{HLPb}, Sec.\,\textbf{2.18}. Chapter II of the book \cite{HLPb}, in
which  majorization is introduced and discussed, contains a
number of references to private communications by Schur.

\noindent \textsf{THEOREM} (G.H.\,Hardy, J.W.\,Littlewood and
G.\,Polya, \cite{HLPb},
sec.\,\textbf{2.20})\\
\begin{slshape}%
\hspace*{2.0ex}\textup{I.} Let
\(\boldsymbol{x}=(x_1,\,\dots\,,x_n)\) and
\(\boldsymbol{y}=(y_1,\,\dots\,,y_n)\) be two sequences of real
numbers and let matrix \(M \) be a doubly-stochastic matrix such
that ${\bf x}=M{\bf y}$. Then
\(\boldsymbol{y}\prec \boldsymbol{x}\).\\
\hspace*{2.0ex}\textup{II.} Let
\(\boldsymbol{x}=(x_1,\,\dots\,,x_n),\) and
\(\boldsymbol{y}=(y_1,\,\dots\,,y_n)\) be two sequences of real
numbers such that \(\boldsymbol{y}\prec \boldsymbol{x}\). Then
there exists a doubly stochastic matrix \(M
\) such that \(\boldsymbol{x}=M\boldsymbol{y}\). (In general such
a matrix \(M\) is not unique.)
\end{slshape}

Part II of this theorem and the first cited theorem of Schur
(which appears as Theorem I in this section) implies the following
result:

\noindent \textsf{THEOREM \(\textup{I}^{\,\prime}\)}.%
\begin{slshape}
Let a sequence \(\boldsymbol{y}=(y_1,\,\dots\,,y_n)\) be majorized
by a sequence \(\boldsymbol{x}=(x_1,\,\dots\,,x_n)\), let
\(x_k,y_k\in(\alpha,\beta)\subset\mathbb{R}\) for
\(k=1,\,\dots\,,n)\), and let \(\varphi\) be a convex function on
the interval \((\alpha,\beta)\).
Then the inequality \textup{(\ref{SGSI})} holds.
\end{slshape}

It turns out that the converse statement is true (\cite{HLP1},
Theorem 8; \cite{HLPb},
Theorem \textbf{108}): %

\begin{slshape}%
Let \(x_k,\,y_k\in (\alpha,\beta)\) for \(k=1,\,\dots\,,n\) and
ssume that the inequality \textup{(\ref{SGSI})} holds for every
function \(\varphi\) which is convex on the interval
\((\alpha,\beta)\). Then \(\boldsymbol{x}=M\boldsymbol{y}\) for
some doubly-stochastic matrix \(M\).
\end{slshape}%

This means that Schur's result (which appears as Theorem I in this
section) is sharp in some sense.

In \cite{GoKr1a}, Chapt.\,II, Lemma 3.5, a very elementary proof
of the following fact
is presented: %
\begin{slshape}%
Let \(\Phi\) be a symmetric function of \(n\) variables which has
continuous derivatives of the first order. Assume that the
condition \textup{(\ref{NMC})} is satisfied. If a sequence
\(\boldsymbol{x}=(x_1,\,\dots\,,x_n)\) of real numbers majorizes
a sequence \(\boldsymbol{y}=(y_1,\,\dots\,,y_n)\), then the
inequality \textup{(\ref{GSI})} holds.
\end{slshape}%

The last result combined  with the Hardy-Littlewood-Polya theorem
that was discussed earlier yields   an independent proof of the
fact that \textsl{  a symmetric function \(\Phi\) that satisfies
the condition \textup{(\ref{NMC})}
is S-convex.}

The theorem by G.\,Birkhoff sheds light on geometric aspects of
majorization and Schur averaging. It is clear that the set of all
doubly-stochastic matrices is compact and convex. Therefore, it is
of interest to find the extreme points of this set. It is clear
that permutation matrices are doubly-stochastic and that they are
extreme points. It turns out that they are the \textsf{only}
extreme points.

\noindent \textsf{THEOREM} (G.\,Birkhoff). %
\begin{slshape}%
Every doubly-stochastic matrix
 \(M=\big[m_{jk}\big]_{1\leq j,k\leq n}\) is
representable as a convex combination of permutation matrices:
\begin{equation}%
\label{ExtPoiRepr}%
M=\sum\limits_{\pi\in\boldsymbol{\cal S}_n}\lambda_{\pi}P_{\pi}\,,
\end{equation}%
where \(\pi\) runs over the set \(\boldsymbol{\cal S}_n\) of all
permutations of the set \(\{1,\,\dots\,,n\}\), \(P_{\pi}\) are the
corresponding permutation matrices \textup{(\ref{PerM})}, and the
coefficients \(\lambda_{\pi}=\lambda_{\pi}(M)\) satisfy the
conditions
\begin{equation}%
\label{ConvCond}%
\lambda_{\pi}\geq 0\ \ (\forall\ \pi\in\boldsymbol{\cal S}_n)\,,\
\ \ \sum\limits_{\pi\in\boldsymbol{\cal S}_n}\lambda_{\pi}=1\,.
\end{equation}%
\end{slshape}%

\noindent%
\textsf{REMARK 2.}
\begin{slshape}
 In general, the coefficients \(\lambda_{\pi}(M)\) in the representation
\textup{(\ref{ExtPoiRepr})} are not uniquely determined from the
matrix \(M\).
\end{slshape}

This theorem was formulated and proved in 1946 in the paper
\cite{Birk1}. (This formulation also appeared in  Example \(4^*\)
in \cite{Birk2}, p.266.) The original proof due to Birkhoff is
based on a theorem by Ph.\,Hall on representatives of subsets,
\cite{HalP}. (The latter theorem can  also be found in
\cite{HalM},\,sec.\textsf{5.1}). G.B.\,Dantzig  \cite{Dan} gives
an algorithm for solving a transportation problem, the solution of
which leads to Birkhoff's theorem. An independent proof of
Birkhoff's theorem was given by J.\,von\,Neumann \cite{NeuJ1} in
the setting of game theory.  ``Combinatorial" proofs of Birkhoff's
theorem (based on Ph.Hall's theorem), are presented in the books
of M.\,Hall \cite{HalM} (see Theorem 5.1.9), and  C.\,Berge
\cite{Ber} (see Theorem 11 in Chapt.\,10). A geometric proof
(based on a direct investigation of extreme points) is presented
in \cite{HoJo}, Theorem \textbf{8.7.1}. Two different proofs of
Birkhoff's theorem are presented in \cite{MaOl}, Chapt.2, Sect.
\textbf{F}. The paper \cite{Mir} is a good survey of
doubly-stochastic matrices. In particular, it contains a proof of
Birkhoff's theorem. See also the problem book by I.M.\,Glazman and
Yu.I.\,Lyubich \cite{GlLy}, Ch.\,7, \S\,4, where Birkhoff's theorem
is presented in problem form.

Let \(\boldsymbol{x}=(x_1,\,\dots\,,x_n)\) be a sequence of real
numbers and, for a permutation \(\pi\) of the set
\(\{1,\,\dots\,,n\}\), let \(\boldsymbol{x}_{\pi}=
(x_{\pi(1)},\,\dots\,,x_{\pi(n)})\). (Thus, for given
\(\boldsymbol{x}\) there are \(n!\) sequences
\(\boldsymbol{x}_{\pi}\),  some of which can coincide.) We
consider
these sequences as vectors in \(\mathbb{R}^n\).  %
Let \(\boldsymbol{\cal C}_{\boldsymbol{x}}\) denote the convex
hull of all the vectors \(\boldsymbol{x}_{\pi}\ \text{where}\
\pi\in\boldsymbol{\cal S}_n\).

\noindent \textsf{THEOREM} (R.\,Rado) %
\begin{slshape}%
Let \(\boldsymbol{x}=(x_1,\,\dots\,,x_n)\) and
\(\boldsymbol{y}=(y_1,\,\dots\,,y_n)\) be two sequences of real
numbers. Then
\[
\boldsymbol{y}\in\boldsymbol{\cal C}_{\boldsymbol{x}} \quad
\Longleftrightarrow \quad
\boldsymbol{y}\prec\boldsymbol{x}.\]%
\end{slshape}%

PROOF. The implication \(\Rightarrow\) is easy. The converse can
be obtained by combining the cited theorems of
Hardy-Littlewood-Polya  and Birkhoff. \hfill \framebox[0.45em]{ }

This theorem seems to have been   established first by R.\,Rado
\cite{Rad}. His proof was based on a theorem on the separation of
convex sets by hyperplanes. A.\,Horn (\cite{HorA1}, Theorem 2)
observed it can also be obtained by combining the results of
Hardy-Littlewood-Polya and Birkhoff that were cited earlier. A
short proof of Rado's theorem, which does not use the Birkhoff
theorem, can be found in \cite{Mark} (see Theorem 1.1).

The circle of ideas related to Schur averaging, majorization and
Birkhoff's theorem is well represented in the literature. The
whole book \cite{MaOl} (of more than 550 pages) is dedicated to
this circle. It includes applications to  combinatorial analysis,
matrix theory, numerical analysis and statistics.  The books
\cite{ArnB} and \cite{PPT} are also relevant. There are
generalizations of Birkhoff's theorem to the infinite dimensional
case, see \cite{Mir}  and \cite{NeuA}.

One generalization of Birkhoff's theorem leads to an interpolation
theorem for linear operators. 
Let \(B\) be the linear space
\(\mathbb{R}^n\) provided with a norm \(\|\,.\,\|_B\) such that
\(\|\boldsymbol{x}\|_{\boldsymbol{B}}=
\|{\boldsymbol{x}}_{\pi}\|_{\boldsymbol{B}}\)
for every \(\boldsymbol{x}\in\mathbb{R}^n\) and for every
permutation \(\pi\in\boldsymbol{\mathcal{S}}_n\), where, as usual,
\(\boldsymbol{x}_{\pi}=(x_{\pi(1)},\,\dots\,,x_{\pi(n)})\). In
other words, this property of the norm \(\|\,.\,\|_{B}\) can be
expressed as \(\|P_{\pi}\|_{B\to B}=1\) for every permutation
\(\pi\in\boldsymbol{\mathcal{S}}_n\) where \textsl{the permutation
operator} \(P_{\pi}\) is defined by the permutation matrix
\(P_{\pi}\), (\ref{PerM}), in the natural basis of the space
\(\mathbb{R}^n\). A norm \(\|\,.\,\|_B\) with this property is
said to be \textsl{a symmetric norm}. A Banach space \(B\) with a
symmetric norm is said to be  \textsl{a symmetric Banach space}.

Let an operator \(A\) in the space \(\mathbb{R}^n\) be defined by
its matrix \(A=\big[a_{jk}\big]_{1\leq k\leq n}\) in the natural
basis of the space \(\mathbb{R}^n\) and  assume that it satisfies
the norm estimates
\begin{equation}%
\label{Contr}%
\|A\|_{l^1\to l^1}\leq 1\quad\text{and}\quad \|A\|_{l^{\infty}\to
l^{\infty}}\leq 1\,.
\end{equation}%
Then, as noted earlier in Section \ref{SchurPr},
\begin{equation}%
\label{ContrE}%
\sum\limits_{1\leq j\leq n}|a_{jk}|\leq 1, \ \ 1\leq k \leq n
\quad\text{and}\quad \sum\limits_{1\leq k\leq n}|a_{jk}|\leq 1, \
\ 1\leq j \leq n\,.
\end{equation}%
According to one generalization of Birkhoff's theorem, a matrix
\(A\) satisfying the conditions (\ref{ContrE}) admits a
representation of the form
\[A=\sum\limits_{\pi\in\boldsymbol{\mathcal{S}}_n}\lambda_{\pi}P_{\pi}\]
where the \(\lambda_{\pi}\) are \textsl{real} (not necessarily
non-negative) numbers satisfying
the conditions %
 \[\sum\limits_{\pi\in\boldsymbol{\mathcal{S}}_n}|\lambda_{\pi}|\leq 1\,.\]
Therefore, since \(\|P_{\pi}\|_{B\to B}=1\), the operator \(A\)
must be a contraction in this norm:
\begin{equation}%
\label{ContrB}%
\|A\|_{B\to B}\leq 1\,.
\end{equation}%
Thus, the following result holds:

\noindent \textsf{THEOREM} (Interpolation theorem for symmetric 
Banach spaces).%
\begin{slshape}%
\quad Let an operator \(A\) acting in the space \(\mathbb{R}^n\)
be a contraction in the \(l^1\) and  \(l^{\infty}\) norms, i.e.,
let the estimates \textup{(\ref{Contr})} hold. Then the operator
\(A\) is a contraction in every symmetric norm \(\|\,\cdot\,\|_B\) on
\(\mathbb{R}^n\), i.e., the estimate \textup{(\ref{ContrB})}
holds.
\end{slshape}%

Here we presented the simplest interpolation result for symmetric
spaces. A more advanced result can be found in \cite{Mit}. Thus,
the development of ideas initiated by  Schur leads to
interpolation theorems for Banach spaces with symmetric norms.

The last topic which we discuss here is the Schur-Horn convexity
theorem. A.\,Horn (\cite{HorA1}, Theorem 4) obtained the following
strengthening of the second part  of the Hardy-Littlewood-Polya
theorem:

\noindent \textsf{THEOREM} (A.\,Horn). %
\begin{slshape}%
Let \(\boldsymbol{x}=(x_1,\,\dots\,,x_n\) and \(\boldsymbol{y}
=(y_1,\,\dots\,,y_n)\) be any two points in \(\mathbb{R}^n\) such
that \(\boldsymbol{y}\prec\boldsymbol{x}\). Then there exists  an
\textsf{ortho}-stochastic matrix \(M\) such that
\(\boldsymbol{y}=M\boldsymbol{x}\). 
\end{slshape}

The following result is a direct consequence of the cited theorems
of Rado and A.\,Horn:

\begin{slshape}%
Given \(\boldsymbol{x}=(x_1,\,\dots\,,x_n)\in\mathbb{R}^n\), the
following two sets are coincide:\\[-3.9ex]
\begin{enumerate}%
\item[\textup{1.}]
The set
\(\boldsymbol{\mathcal{C}}_{\boldsymbol{x}}\stackrel{\mathrm{def}}{=}
\ \text{the convex hull of the family of vectors}\
\{\boldsymbol{x}_{\pi}\}_{\pi\in\boldsymbol{\mathcal{S}}_n}\)\,.\\[-3.9ex]
\item[\textup{2.}]
The set \(\{M\boldsymbol{x}\}\), where \(M\) runs over the set of
all \textsf{ortho}-stochastic matrices\,.
\end{enumerate}%
\end{slshape}

In view of the relations (\ref{DiagF}) and (\ref{OrtAv}), the last
statement can be reformulated in terms of eigenvalues and diagonal
entries: Let us associate with every \textsl{real symmetric}
\(n\times n\) matrix \(H=\big[h_{jk}\big]_{1\leq j,k\leq n}\) the
\(n\)-tuple \(\boldsymbol{h}(H)=(h_{11},\,\dots\,,h_{nn})\) of its
diagonal entries and the \(n\)-tuple
\(\boldsymbol{\omega}(H)=({\omega}_{1}(H),\,\dots\,,\omega_{n}(H))\)
of its eigenvalues arranged in non-increasing order:
\({\omega}_1(H)\geq,\,\dots\,,\geq{\omega}_n(H)\). We consider
these \(n\)-tuples as vectors in \(\mathbb{R}^n\).
 Given an \(n\)-tuple \(\boldsymbol{\omega}=({\omega}_1,\,\dots\,,{\omega}_n)\)
 of real numbers, arranged in non-increasing order:
 \({\omega}_1\geq\,\dots\,,{\omega}_{n}\), let
\[\boldsymbol{\mathcal{H}}_{\boldsymbol{\omega}}=%
\big\{H:\ H\ \text{is real symmetric and}\ %
\boldsymbol{\omega}(H)=\boldsymbol{\omega}\big\}\,.\]
\textsf{THEOREM} (Schur-Horn convexity theorem). %
\begin{slshape}
 Given an \(n\)-tuple \(\boldsymbol{\omega}=({\omega}_1,\,\dots\,,{\omega}_n)\)
 of real numbers: \({\omega}_1\geq\,\dots\,,\geq {\omega}_{n}\), the set
\(\big\{\boldsymbol{h}(H)\big\}_{H\in%
{\boldsymbol{\mathcal{H}}_{\boldsymbol{\omega}}}}\) of all
``diagonals" of matrices from
\(\boldsymbol{\mathcal{H}}_{\boldsymbol{\omega}}\) is convex.
Moreover,
\begin{equation}%
\label{SchHornEqual}%
\big\{\boldsymbol{h}(H)\big\}_{H\in%
{\boldsymbol{\mathcal{H}}_{\boldsymbol{\omega}}}}= %
{\boldsymbol{\mathcal{C}}}_{\boldsymbol{\omega}},\,
\end{equation}
where \({\boldsymbol{\mathcal{C}}}_{\boldsymbol{\omega}}\) is the
convex hull
of the family of \(n!\) vectors %
\({\boldsymbol{\omega}}_{\pi}=({\omega}_{\pi(1)},\,\dots\,,{\omega}_{\pi(n)})%
\), as \(\pi\) runs over the set \({\boldsymbol{\mathcal{C}}}_n\)
of all permutations of the set \(\{1,\,\dots\,,n\}\):
\begin{equation}%
\label{DefConvSet}%
{\boldsymbol{\mathcal{C}}}_{\boldsymbol{\omega}}=\textup{Conv}\,%
\big\{{\boldsymbol{\omega}}_{\pi}:\pi\in\boldsymbol{\mathcal{S}}_n \big\}.
\end{equation}%
\end{slshape}
 Schur himself established the formula
\[
\big\{\boldsymbol{h}(H)\big\}_{H\in%
{\boldsymbol{\mathcal{H}}_{\boldsymbol{\omega}}}}= %
\{M\boldsymbol{\omega}: M \quad \textnormal{is
ortho-stochastic}\}. \] He did not described the set
 on the right geometrically as a convex hull.
The term ``convex
 set" does not appear in the paper {\cite{Sch18a}) at all. The ``Schur-Horn
 convexity theorem" appeared only in the paper by A. Horn \cite{HorA2}
 ( which used in an essential way the cited results by
Hardy-Littlewood-Polya and Birkhoff.) However, the influence of
Issai Schur on the area was so great that the term ``Schur-Horn
convexity theorem" is now common.

In the last thirty years, the Schur-Horn convexity theorem has
been generalized significantly. In 1973 (fifty years after the
publication of \cite{Sch18a}) B.\,Kostant published a seminal
paper \cite{Kos} in which he interpreted the Schur-Horn result as
a property af adjoint orbits of the unitary group and generalized
it to arbitrary compact Lie groups. More precisely, he proved (see
especially \cite{Kos}, sect.\,8) that for an element \(x\) in a
maximal abelian subspace \(\mathfrak{t}\) in the Lie algebra
\(\mathfrak{k}\) of a
compact Lie group \(K\) one has %
\[pr_{\mathfrak{t}}(\textup{Ad}\,K_{\cdot}x)
=\textup{Conv}\,\mathcal{W}_{\cdot}x\,,\] where
\(pr_{\mathfrak{t}}:\,\mathfrak{k}\to\mathfrak{t}\) is the
orthogonal projection (with respect to the Killing form) and
\(\mathcal{W}\) is the Weyl group associated with the pair
\((\mathfrak{k}_{\mathbb{C}}, \mathfrak{t}_{\mathbb{C}})\).
Subsequently, M.F.\,Atiyah \cite{Ati} and, independently,
V.\,Guillemin and S.\,Sternberg \cite{GuSt1}, \cite{GuSt2} gave an
interpretation of Kostant's theorem as a special case of a theorem
on the image of the momentum map of a Hamiltonian torus action.
Atiyah's proofs depend on some ideas from Morse theory.
Subsequently, the results of Kostant, Atiyah, Guillemin and
Sternberg were extended to the setting of symmetric spaces. See,
for example, the paper \cite{HNP}, where more references can be
found, the paper \cite{BFR} and the book \cite{HiOl}, sections 4.3
and 5.5.

In yet another direction, the relevance of doubly-stochastic matrices
and Schur averaging to operator algebras and quantum physics is discussed
in the book \cite{AlU}.
     
Thus, once again a relatively short paper of Issai Schur is seen
to have had significant influence on the development of a number
of diverse areas of mathematics. In particular, \cite{Sch18a}
paved the way to important results in matrix theory, statistics,
the theory of Lie groups and symmetric spaces, symplectic geometry
and Hamiltonian mechanics. Many of these areas are very far from
the original setting.

\setcounter{equation}{0}
\begin{minipage}{15.0cm}
\section{\hspace{-0.5cm}.\hspace{0.2cm}%
{\large \ Inequalities between the eigenvalues and the singular
\newline
 values of a linear operator. \label{InSinEig}}}
\end{minipage}\\[0.18cm]

Let \(A=\big[a_{jk}\big]_{1\leq j,k \leq n}\) be an \(n\times n\)
matrix with eigenvalues \(\lambda_1,\,\dots\,,\lambda_n\in
\mathbb{C}\). In Theorem II of \cite{Sch2a},  Schur proved the
inequality
\begin{equation}%
\label{Stwo}%
\sum\limits_{\ell =1}^n|\lambda_{\ell}|^2\leq
\sum\limits_{j,k=1}^n |a_{jk}|^2.
\end{equation}%

Schur's proof was based on Theorem I of that paper, in which he
established the fundamental fact that every square matrix $A$ with
complex entries is unitarily equivalent to an upper triangular
matrix, i.e., there exists a unitary matrix $U$ such that
\begin{equation}%
\label{UTF}%
T = U^*AU = U^{-1}AU
\end{equation}%
is  upper triangular: \(\mathsf{t}_{jk}=0\) for \(j>k\).
Therefore, the set of eigenvalues of the matrix \(A\) is equal to
the  set of eigenvalues of the matrix \(T\), which in turn is
equal to the set of  diagonal entries of \(T\). Thus,
$$
\sum\limits_{\ell=1}^n |\lambda_{\ell}|^2 = \sum\limits_{j=1}^n
|t_{jj}|^2 \leq \sum\limits_{j,k=1}^n |t_{jk}|^2 =
\text{trace}\,T^*T = \text{trace}\,A^*A = \sum\limits_{j,k=1}^n
|a_{jk}|^2.
$$

Apart from its use in the proof of the inequality (\ref{Stwo}),
Theorem I  serves as a model for some important constructions in
operator theory that will be discussed below.

In \cite{Sch2a}, Schur used (\ref{Stwo}) to obtain simple proofs
of the estimates
\begin{equation}%
\label{Sinf}%
|\lambda_l|\leq n\cdot \max\limits_{1\leq j,k \leq
n}|a_{jk}|\qquad (1\leq l\leq n)
\end{equation}%
\begin{equation}%
\label{ImEig}%
\big|\text{Re}\,\lambda_l\big|\leq n\cdot \max\limits_{1\leq
j,k\leq n}\,|b_{jk}| \quad \mbox{and} \quad
\big|\text{Im}\,\lambda_l\big|\leq n\cdot \max\limits_{1\leq
j,k\leq n}|c_{jk}|\quad (1\leq l\leq n)\,,
\end{equation}%
for the eigenvalues \(\lambda_l\) of a general \(n\times n\)
matrix \(A=\big[a_{jk}\big]\), where,
\(B=\big[b_{jk}\big]=(A+A^*)/2 \ \text{and} \
C=\big[c_{jk}\big]=(A-A^*)/(2i)\,.\) The estimates (\ref{ImEig})
were first obtained by A.Hirsch \cite{Hir}. They were improved to
\begin{equation}%
\label{ImEigR}%
\big|\text{Im}\,\lambda_l\big|\leq \sqrt{\frac{n(n-1)}{2}}\cdot
\max\limits_{1\leq j,k\leq n}|c_{jk}|\quad
 (1\leq l\leq n)\,
\end{equation}%
for real matrices \(A\) by F.\,Bendixson \cite{Bend} and
 reproved in \cite{Sch2a}. In \S 7 of (\cite{Sch2a}), the interesting
inequality
\begin{equation}%
\label{Ford}%
\sum\limits_{j<k}|\lambda_j-\lambda_k|^2\leq
\sum\limits_{j<k}|a_{jj}-a_{kk}|^2+n\sum\limits_{j\not=k}|a_{jk}|^2
\end{equation}%
is derived and then used to  obtain the following estimate  for
the discriminant
\[d=\prod\limits_{j,k}(\lambda_j-\lambda_k)^2\,,\]
 of the characteristic equation
\(\text{det}\,(\lambda I_n -A)=0\):
\begin{equation}%
\label{ForD}%
|d|^{\frac{2}{n(n-1)}}\leq\frac{2}{n(n-1)}
\sum\limits_{j<k}|a_{jj}-a_{kk}|^2+
\frac{2}{n-1}\sum\limits_{j\not =k}|a_{jk}|^2\,.
\end{equation}%
In \S 5 of \cite{Sch2a}, the well known Hadamard bound
\begin{equation}%
\label{HIMVD}%
|\det A|\leq\big(\max_{1\leq j,k\leq n}|a_{jk}|\,\big)^n\cdot
n^{n/2}\,.
\end{equation}%
on the maximal value of the determinant of a matrix is derived
from the inequality (\ref{Stwo}) with the help of  the inequality
between the geometric and the arithmetic means:
\[|\det A|^2=|\lambda_1|^2\cdot\,\cdots\,\cdot|\lambda_n|^2\leq%
\left(\frac{|\lambda_1|^2+\,\cdots\,+|\lambda_n|^2}{n}\right)^n\leq%
\left(\frac{\sum_{j,k=1}^n|a_{jk}|^2}{n}\right)^n\,.
\]
The challenge of obtaining simple new proofs
 of various Hadamard inequalities seems to have been one
 of Issai Schur's favorite occupations.

 In \cite{Sch2a},  Schur also considers  integral operators
\(x(t)\rightarrow (\boldsymbol{K}x)(t)\) in \(L^2(a,b)\),
\begin{equation}%
\label{IO}%
(\boldsymbol{K}x)(t)=\int\limits_a^bK(t,\tau)x(\tau)\,d\tau\quad
(a\leq t \leq b)\,,
\end{equation}%
with kernels \(K(t,\tau)\) that satisfy the condition
\begin{equation}%
\label{HS}%
\int\limits_a^b\!\!\!\int\limits_a^b|K(t,\tau)|^2\,dt\,d\tau<\infty\,.
\end{equation}%
Today, such operators are commonly called  \textsl{Hilbert-Schmidt
integral operators}.  Schur extended the inequality (\ref{Stwo})
to these operators:
\begin{equation}%
\label{StwoHS}%
\sum\limits_l|\lambda_l(\boldsymbol{K})|^2\leq
\int\limits_a^b\!\!\!\int\limits_a^b|K(t,\tau)|^2\,dt\,d\tau\,,
\end{equation}%
where the summation on the left hand side is extended over the set
of all eigenvalues \(\lambda_l(\boldsymbol{K})\) of the integral
operator \(\boldsymbol{K}\). In particular, the series on the left
hand side of (\ref{StwoHS}) converges.

One of the fundamental results of the Fredholm theory of integral
equations \cite{Fred} is the identification of the nonzero
eigenvalues \(\lambda_l(\boldsymbol{K})\) of an integral operator
(\ref{IO}) with a continuous kernel as the reciprocals of the
zeros of an entire function \(D_{\boldsymbol{K}}(\lambda)\) (that
is constructed from the kernel \(K(t,\tau)\) of this operator).
This function is termed \textsl{the Fredholm denominator (or the
Fredholm determinant) of the operator (\ref{IO})}.
It is defined by the Taylor series
\begin{equation}%
\label{FD1}%
 D_{\boldsymbol{K}}(\lambda)=\sum\limits_{n=0}^{\infty}c_n\lambda^n\,,
\end{equation}%
with  coefficients
\begin{equation}%
\label{FD2}%
c_n=\frac{(-1)^n}{n!}\int\limits_{a}^{b}\dots\int\limits_{a}^{b}
\det
\left[%
\begin{matrix}%
K(t_1,t_1) &\,\cdots\,&K(t_1,t_n)\\
\cdots &\cdots&\cdots\\
K(t_n,t_1)&\cdots&K(t_n,t_n)
\end{matrix}%
\right]\,dt_1\,\cdots\,dt_n\,.
\end{equation}%
From (\ref{FD2}) and  the Hadamard inequality (\ref{HIMVD}), it
follows that if%
$$
\sigma=(b-a) (\max_{a\leq t,\tau \leq b}|K(t,\tau)|) < \infty \,,
$$
then
\begin{equation}%
\label{EstForC}%
|c_n|\leq\sigma^n \cdot n^{n/2}/n!\,.
\end{equation}
Consequently, the series (\ref{FD1}) converges for every complex
\(\lambda\), and its sum \(D_{\boldsymbol{K}}(\lambda)\) is an
entire function that is subject to the bound
\begin{equation}%
\label{EFDet}%
\ln|D_{\boldsymbol{K}}(\lambda)|\leq\sigma^2|\lambda|^2(1+o(1))
\quad (|\lambda|\to\infty).
\end{equation}%
Thus, the counting function of the zeros $\mu_1, \mu_2, \ldots $
of $D_{\boldsymbol{K}}(\lambda)$:
\[n_{\boldsymbol{K}}(r)=\#\{\mu_{\ell}(\boldsymbol{K}): %
|\mu_{\ell}(\boldsymbol{K})|\leq r\} =
\#\{\lambda_{\ell}(\boldsymbol{K}): %
|\lambda_{\ell}(\boldsymbol{K})|^{-1}\leq r\}\,,\] satisfies the
condition
\begin{equation}%
\label{UnEst}%
n_{\boldsymbol{K}}(r)=O(r^2),\quad\text{as}\quad r\to\infty\,.
\end{equation}%

The estimates (\ref{EFDet}) and its consequence (\ref{UnEst}) were
known \cite{Lal} before the Schur paper \cite{Sch2a} appeared.
However, the estimate (\ref{StwoHS}) is stronger than the estimate
(\ref{UnEst}). From the convergence of the series
\(\sum\limits_l|\lambda_l(\boldsymbol{K})|^2\) and from the
estimate (\ref{EFDet}) it follows that the Fredholm denominator
(\ref{FD1})-(\ref{FD2}) admits the multiplicative decomposition
\begin{equation}%
\label{MRFD}%
D_{\boldsymbol{K}}(\lambda)=e^{c\lambda+d\lambda^2}%
\prod\limits_{l}\big(1-\lambda\,\lambda_l(\boldsymbol{K})\big)\,%
e^{\lambda\,\lambda_l(\boldsymbol{K})}\,,
\end{equation}%
for some choice of constants \(c\) and \(d\). The fact that the
Fredholm denominator of the integral operator (\ref{IO}) with a
continuous kernel admits a representation of the form (\ref{MRFD})
was first noted by Schur in \S\,14\ of \cite{Sch2a}. (It is
important to note that the kernel \(K(t,\tau)\) is not assumed to
be symmetric or Hermitian.) This result of Schur is sharp in the
sense that there exists a continuous kernel \(K\) on a finite
interval \([a,\,b]\) whose
eigenvalues satisfies the condition  %
\[\sum_{\ell}|\lambda_{\ell}(\boldsymbol{K})|^{2-\epsilon}=\infty
 \ \mbox{for every} \
\epsilon>0.\]
To construct an example, let \(K(t,\tau)=\varphi(t-\tau)\) for
\(0\leq t,\tau\leq 1\), where \(\varphi(t+1)=\varphi(t)\) is a
continuous periodic function on \(\mathbb{R}\) with Fourier
expansion \(\varphi(t)\sim\sum_lc_{\ell}e^{2\pi i\ell t}\). Then
the functions \(e^{2\pi i\ell t}\) are eigenfunctions of the
kernel \(K\), and the Fourier coefficients \(c_{\ell}\) are
eigenvalues of this kernel. A kernel with the desired properties
is obtained by choosing a
 continuous
periodic function \(\varphi\) whose Fourier coefficients
\(c_{\ell}\) satisfy the condition
\(\sum_{\ell}|c_{\ell}|^{2-\epsilon}=\infty\) for every
\(\epsilon>0\). The first example of such a function was
constructed by T.\,Carleman \cite{Carl2}.
 Other examples can be
found in \cite{Bar}, Chapt.\,4, \S\,16, or in \cite{Zyg},
Chapt.\,5., \textbf{(4.9)}. In his first publication \cite{Carl1},
Carleman proved that in fact \(d=0\) in (\ref{MRFD}). Thus, the
scientific career of this outstanding analyst started with an
improvement of a result of Issai Schur.

The inequality (\ref{Stwo}) can also be presented in the form
\begin{equation}%
\label{EigSingTwo}%
\sum\limits_{l=1}^n|\lambda_l(A)|^2\leq\sum\limits_{l=1}^ns_l(A)^2\,,
\end{equation}%
where the \(\lambda_l(A)\) are the eigenvalues of the matrix \(A\)
and the numbers \(s_l(A)\) are the singular values of $A$.

The auxiliary inequality
\begin{equation}%
\label{EigSingOne}%
\sum\limits_{l=1}^n|\lambda_l(A)|\leq\sum\limits_{ l=1}^ns_l(A)
\end{equation}%
can also be proved\,\footnotemark %
\footnotetext{This proof is adopted from \cite{GoKr1b},
Chapt.\,IV,\,\S\,8. See Theorem 8.1, especially the footnote 7 on
p.\,128 of the Russian original or on p.\,98 of the English
translation.} in an elementary way by
 using the Schur transformation
(\ref{UTF}) to reduce the matrix \(A\) to  upper triangular form.
In fact, it suffices to prove ({\ref{EigSingOne}) for upper
triangular
  matrices \(A\),
since  the transformation (\ref{UTF}) does not change either the
eigenvalues or the singular values of the matrix. But then, if
\(\{e_l\}_{1\leq l\leq n}\) is the natural basis of the space
\(\mathbb{C}^n\),
\[ a_{ll}=\langle Ae_l,e_l\rangle=\lambda_l(A), \quad l=1,\,\ldots\,,n, \]
up to a reindexing of the eigenvalues, if need be. Now let
\(A=S\cdot V\) be the polar decomposition of the matrix \(A\):
\(S\geq 0, V^*V=VV^*=I_n\) and let   \(h_l=Ve_l\). Then the
vectors \(\{h_l\}_{1\leq l\leq n}\) form an orthonormal basis of
the space \(\mathbb{C}^n\) and, by the Cauchy-Schwarz inequality,
\[|\langle Ae_l,e_l\rangle|=|\langle Sh_l,e_l\rangle|
\leq \sqrt{\langle Sh_l,h_l\rangle}\,\cdot\sqrt{\langle
Se_l,e_l\rangle}\,.\] Therefore,
\[\sum\limits_{l=1}^n|\lambda_l(A)|=
\sum\limits_{l=1}^n|\langle Ae_l,e_l\rangle|\leq
\sqrt{\sum_{l=1}^n\langle
Sh_l,h_l\rangle}\,\cdot\,\sqrt{\sum_{l=1}^n \langle
Se_l,e_l\rangle}\,=\sum_{l=1}^n\lambda_l(S)=
\sum\limits_{l=1}^ns_l(A)\,,
\]
since \[\sum_{l=1}^n\langle Sh_l,h_l\rangle= \sum_{l=1}^n\langle
Se_l,e_l\rangle=\mbox{trace}\, S = \sum_{l=1}^n \lambda_l(S),\]
and, by the definition of
singular values, \(\{\lambda_l(S)\}_{l=1}^n = \{s_l(A)\}_{l=1}^n \). %

 The inequalities (\ref{Stwo}), \textsl{written in the form}
(\ref{EigSingTwo}), and (\ref{EigSingOne}) were  significantly
generalized by
 H. Weyl \cite{Wey} in 1949. The generalization is based on
the concept of  majorization that was discussed in the previous
section. A crucial role is played by the inequalities
\begin{equation}%
\label{MISE}%
|\,\lambda_1(A)\cdot
\lambda_2(A)\cdot\,\,\cdots\,\,\cdot\lambda_k(A)\,|\leq
s_1(A)\cdot s_2(A)\cdot\,\,\cdots\,\,\cdot s_k(A)\quad
(k=1,\,2,\,\dots\,,\,n-1)\,,
\end{equation}%
which are valid when the  eigenvalues \(\lambda_k(A)\) and the
singular values \(|s_k(A)|\) are indexed in such a way that
\(|\lambda_1(A)|\geq|\lambda_2(A)|\geq\,\cdots\geq |\lambda_n(A)|
\quad \mbox{and} \quad s_1(A)\geq s_2(A)\geq\,\cdots\,\geq
s_n(A).\) The equality
\begin{equation}%
\label{MESE}%
|\,\lambda_1(A)\cdot
\lambda_2(A)\cdot\,\,\cdots\,\,\cdot\lambda_n(A)\,|= s_1(A)\cdot
s_2(A)\cdot\,\,\cdots\,\,\cdot s_n(A)
\end{equation}%
holds because both  sides are equal to
\(\big|\text{det}\,A\big|.\) The relations (\ref{MISE}) and
(\ref{MESE}) mean that the sequence
\(\{\ln|\lambda_k(A)|\}_{k=1}^n\) is majorized by the sequence
\(\{\ln\,s_k(A)\}_{k=1}^n\):
\begin{equation}%
\label{MajES}%
\{\,\ln|\lambda_k(A)|\,\}_{k=1}^n \prec
\{\,\ln\,s_k(A)\,\}_{k=1}^n\,.
\end{equation}%
In \cite{Wey}, Weyl derived the inequalities (\ref{MISE}) and then
applied  the inequality
\begin{equation}%
\label{MajIn}%
\sum\limits_{k=1}^n\psi(y_k)\leq \sum\limits_{k=1}^n\psi(x_k)\,,
\end{equation}%
which holds for any convex function \(\psi(\,\cdot\,)\)  on
\((-\infty,\,\infty)\) and any pair of sequences \(\{y_k\}\) and
\(\{x_k\}\) such that \(\{y_k\}\prec\{x_k\}\), to the sequences
\(y_k=\ln|\lambda_k(A)|\) and \(x_k=\ln\,s_k(A)\).
 The inequality
(\ref{MajIn}) is a direct consequence of the result%
\footnotemark\addtocounter{footnote}{-1} \ by Schur (which states
that the inequality (\ref{MajIn}) holds for  sequences
\(\boldsymbol{x}\) and \(\boldsymbol{y}=M\boldsymbol{x}\) that are
related by a doubly-stochastic matrix \(M\)), and of the
result\footnotemark \ by Hardy, Littlewood and Polya, who proved
that
 \[\boldsymbol{y}\prec\boldsymbol{x} \quad
\Longrightarrow \quad
 \boldsymbol{y}=M\boldsymbol{x}\] for some doubly-stochastic
matrix \(M\). \footnotetext{These results were discussed in the
previous section; see Theorems I, II and \(\text{I}^{\prime}\).}
However, Weyl was not aware of these results and gave an
independent proof of the implication
\begin{equation}
\label{weylimp} \boldsymbol{y}\prec\boldsymbol{x} \quad
\Longrightarrow \quad
 (\ref{MajIn})
\end{equation}
in Lemma 1 of \cite{Wey}.
 The  inequalities (\ref{MISE}) were known before the paper
\cite{Wey} was published. (See, for example, Exercise 17 on page
110 of the book \cite{TuAi}.) However, it was Hermann Weyl who
first combined the inequalities (\ref{MISE}) with the implication
(\ref{weylimp}) to obtain
 the following

\noindent \textsf{THEOREM}. %
\begin{slshape}%
Let \(A\) be an \(n\times n\) matrix with eigenvalues
 \(\{\lambda_k(A)\}_{k=1}^n\)
and singular values \(\{s_k(A)\}_{k=1}^n\)
 (counting  multiplicities) and  let
\(\varphi(\,\cdot\,)\) be a function on \((0,\,\infty)\) such that
the function \(\psi(t)=\varphi(e^t)\) is convex on
\((-\infty,\,\infty)\). Then
\begin{equation}%
\label{WeIn}%
\sum\limits_{k=1}^n\varphi(\lambda_k(A))\leq \sum\limits_{ k=1}^
n\varphi(s_k(A))\,.
\end{equation}%
\end{slshape}%

Weyl invoked the inequality (\ref{WeIn}) with \(\varphi(t)=t^p\)
and  \(p>0\) to obtain
 the following generalization of Schur's
inequality (\ref{EigSingTwo}):
\begin{equation}%
\label{WIP}%
\sum\limits_{l=1}^n|\lambda_l(A)|^p\leq\sum\limits_{l=1}^n
s_l(A)^p\qquad (0<p<\infty).
\end{equation}%
Analogous inequalities hold for  linear operators  \(A\) in a
Hilbert space \(\mathcal{H}\) that belong to the class
\(\mathfrak{S}_p\), i.e., for which \(\sum_{l}s_l(A)^p<\infty\,\),
where the \(s_l(A)\) are the eigenvalues of the operator
\(\sqrt{A^*A}\). (Usually the singular values \(s_l(A)\) are
enumerated by the indices \(l=0,\,1,\,2,\,\dots\).) The summation
in the last inequality is then extended over all eigenvalues and
over all singular values of the operator \(A\). The resulting
inequality is very useful in the theory of integral equations. The
point  is that it is difficult to  calculate the eigenvalues
and singular values  of an integral operator in terms of its
kernel. However,  the singular values can be effectively estimated
from above by approximating the kernel \(K(t,\tau)\) by degenerate
kernels of the form
\(K_n(t,\tau)=\sum\limits_{l=1}^n\varphi_l(t)\psi_l(\tau)\) and
invoking the fact that
\[s_n(\boldsymbol{K})=\inf \|\boldsymbol{K}-\boldsymbol{K}_n\|\]
as \(K_n(t,\tau)\) runs over the set of all degenerate kernels of
the indicated form. (See \cite{GoKr1b}, Chapt.\,2, \S\,2, item
\textbf{3}.) The smoother the kernel \(K\), the more rapid  the
rate of decay of the sequence
\(\|\boldsymbol{K}-\boldsymbol{K}_n\|\) and thus, the rate of
decay of the sequence \(|s_n{\boldsymbol{K}})|\). The inequality
\begin{equation}%
\label{WeyIn}%
\sum\limits_l|\lambda_l(\boldsymbol{K})|^p\leq%
\sum\limits_{l=0}^{\infty}s_l(\boldsymbol{K})^p
\end{equation}%
is then used to
 derive the rate of decay of the eigenvalues
\(\lambda_n(\boldsymbol{K})\).
This is the ``modern" way to derive the rate of decay of the
eigenvalues \(\lambda_n(\boldsymbol{K})\) of an integral operator
from the smoothness of its kernel \(K(t,\tau)\). The theory of
spline approximation is often used to
 construct good approximating
kernels. (See, for example, the papers by M.\,Sh.\,Birman and
M.Z.\,Solomyak mentioned in Section 4.)

The ``classical" approach, which does not exploit the Weyl
inequalities, is more complicated and gives weaker results. Chang,
in a paper \cite{Chang} that appeared before the paper \cite{Wey},
proved that
\begin{equation}%
\label{Impl}%
\sum\limits_{l=0}^{\infty}s_l(\boldsymbol{K})^p<\infty \quad
\Longrightarrow \quad
\sum\limits_{l}|\lambda_l(\boldsymbol{K})^p<\infty
\end{equation}%
for  integral operators (\ref{IO}) of Hilbert-Schmidt class, i.e.,
with kernels \(K(t,\tau)\) satisfying the condition (\ref{HS}).
The ``classical" methods of the paper \cite{Chang} are involved
and rather difficult.

The Weyl inequalities are also useful  in ``abstract" operator
theory.
Taking \(\varphi(t)=\ln(1+|\lambda|t)\) %
(which is admissible, since the function
\(\psi(t)=\ln(1+|\lambda|e^t)\) is convex), one can obtain the
inequality
\begin{equation}%
\label{ForPolyn}%
\Big|\,\prod_{l}\big(1-\lambda\,\lambda_l(A)\big)\,\Big|\leq%
 \prod_{l}\big(1+|\lambda|\,s_l(A)\big)\qquad (\forall
 \lambda\in\mathbb{C})\,,
\end{equation}%
 for linear operators \(A\) from the class \(\mathfrak{S}_1\)
of trace class operators in a Hilbert space. This inequality is
useful in the  study of the so-called characteristic determinants
of trace class operators and related analytic considerations. (See
Chapter IV of \cite{GoKr1b}.) In particular, the inequality
(\ref{ForPolyn}) plays an important role in the proof of a theorem
by V.B.\,Lidski\u{\i}, which
 states that the matricial trace and the spectral
trace of a trace class operator coincide. (See \cite{Lid}, and
\cite{GoKr1b}, Chapt.\,III, \S\,8, Theorem 8.1.) This theorem is
of principal importance in operator theory.

The Weyl inequality (\ref{WeIn}) is one of the central tools in
the toolbox of modern operator theory. However, as Weyl himself
wrote \cite{Wey}, the first step was taken by Schur:

 \textit{"Long ago I.\,Schur proved \textup{(\ref{WIP})}
for\,\,\footnote{\,This reference by Weyl is not accurate. Schur
proved the inequality (\ref{WIP}) for \(p=2\), but not for
\(p=1\).} \(p=1\). Recently S.H. Chang showed in his thesis that,
in the case of integral equations, the convergence of \(\sum
s_l^p\) implies the convergence of \(\sum |\lambda_l|^p\). These
two facts led me to conjecture the relation \textup{(\ref{WIP})},
at least for \(p\leq 1\). After having conceived a simple idea for
the proof, I discussed the matter with C.L.\,Siegel and J. von
Neumann; their remarks have contributed to the final form and
generality in which the results are presented here."}

Thus, the paper \cite{Sch2a} served as source of inspiration for
both T.\,Carleman and H.\,Weyl.

\setcounter{equation}{0}
\begin{minipage}{15.0cm}
\section{\hspace{-0.5cm}.\hspace{0.2cm}%
{\large \ Triangular representations of matrices \newline%
\hspace*{9.0ex} and linear operators.\label{TriangReprLinOp}}}
\end{minipage}\\[0.18cm]

One of the  important theorems of Schur that was discussed in the
preceding section states that every square matrix is unitarily
equivalent to a triangular matrix. In the early fifties,
stimulated by this theorem of Schur,
Moshe (Mikhail Samu{\i}lovich) Liv\v{s}ic  (\,=\,Livshits )
obtained an analogue of this result for a class of
 bounded linear operators in a separable Hilbert space.

To explain his results, let us first recall that every bounded
operator \(A\) in a Hilbert space \(\mathcal{H}\) is representable
in the form
\begin{equation}%
\label{RIP}%
 A=B_A+iC_A,\quad
\end{equation}%
where
\begin{equation}%
\label{RiP}%
B_A=\text{Re}\,A =\frac{A+A^*}{2}\,=(B_A)^*  \quad \text{and}
\quad C_A=\text{Im}\,A =\frac{A-A^*}{2i}\,=(C_A)^*.
\end{equation}%

Liv\v{s}ic obtained his conclusions in the class $i\Omega$ of
bounded linear operators $A$ for which $C_A$ is of of trace class.
In the simplest case of this setting,
\begin{equation}%
\label{RankOne}%
\text{rank}(C_A)=1\,,
\end{equation}%
and hence $C_A$ must be be  definite: either \(C_A \geq 0\), or
\(C_A \leq 0\). Thus,
\begin{equation}
\label{PR}%
C_A=j\,|C_A|,\quad \text{where}\quad |C_A|=\sqrt{(C_A)^*C_A},\  \
j=+1\ \ \text{or}\ \ j=-1\,
\end{equation}%
and the imaginary parts $\beta_k$ of the eigenvalues \(\lambda_k
=\alpha_k+i\beta_k\) of the operator \(A\) are of the form
$\beta_k=j|\beta_k|$.

Without loss of generality, we may assume that the the operator
\(A\) is \textsl{completely non-selfadjoint}:
\hspace*{2.0ex}\textsl{There is no invariant subspace for the
operator \(A\) on which \(A\) induces a self-adjoint operator.}
Indeed, if the operator \(A\) is not completely non-selfadjoint,
then it splits into the orthogonal sum \(A=S\oplus
A_{\text{\,cns}}\), where \(S\) is a selfadjoint operator and
\(A_{\text{\,cns}}\)  is a completely non-selfadjoint operator.
Moreover, \(\text{Im}\,A=\text{Im}\,A_{\text{\,cns}}\).

The eigenvalues \(\lambda_k=\alpha_k+i\beta_k\) of a completely
non-selfadjoint operator \(A\) with non-negative (non-positive)
imaginary part are never real: either \(\beta_k>0\) for all \(k\)
(if \(C_A\geq 0\)), or \(\beta_k<0\) for all \(k\) (if \(C_A\leq
0\)).

The triangular model \(T\) for an operator \(A\) satisfying the
condition (\ref{RankOne}) acts in the model Hilbert space
 \(\mathcal{H}_{\text{mod}}=l^2\oplus L^2\) that is the orthogonal
sum of the space \(l^2\)  of  square summable one-sided infinite
sequences \((\xi_1,\xi_2,\,\dots\,,)\) of complex numbers of
dimension $n\leq\infty$, where $n$ is equal to the number of
eigenvalues of the operator \(A\) (counting  multiplicities), and
\(L^2\) is the space  of all square summable complex-valued
functions on a finite interval \([0,l]\), where the number \(l\)
is determined  uniquely by the operator \(A\). The spaces \(l^2\)
and \(L^2\) are equipped with the standard scalar products. The
block decomposition of the operator \(T\) that corresponds to the
decomposition \(\mathcal{H}_{\text{mod}}=l^2\oplus L^2\) of the
space \(\mathcal{H}_{\text{mod}}\), is of the form
\begin{equation}%
\label{TrForm}%
T=
\begin{bmatrix}%
T_{\,\text{dis}} & T_{\,\text{cou}}\\
0 & T_{\,\text{con}}
\end{bmatrix}\,,%
\end{equation}%
 where \(T_{\text{dis}}:l^2\to
l^2\), \(T_{\text{con}}:L^2\to L^2\)\quad\mbox{and}\quad
\(T_{\text{cou}}:L^2\to l^2\).

The operator \(T_{\text{dis}}\), the
discrete part of the operator \(T\), is  defined by its matrix  %
\(\big[t_{km}\big]\) %
in the natural basis of the space \(l^2\). This matrix is upper
triangular, i.e., with $j$ as in (\ref{PR}), 
\begin{equation}%
\label{MDP}%
t_{km}=0 \ \ \text{for} \  \ k>m,\quad t_{kk}=\lambda_k, \quad
t_{km}=i\,|\beta_k|^{1/2}\,j\,|\beta_m|^{1/2}\ \ \text{for}\ \
k<m\,.
\end{equation}%
The operator \(T_{\text{dis}}\) is bounded, since $A$ is bounded
and \(\sum\limits_k|\beta_k|\leq \text{trace}\,|C_A|<\infty\,\).
The operator \(T_{\text{con}}\), the continuous part of the
operator \(T\), is an integral operator of the form
\begin{equation}%
\label{MCP}%
(T_{\text{con}}\xi)(t)=
\lambda(t)\xi(t)+i\int\limits_t^l\!K(t,s)\,\xi(s)\,ds\quad
0\leq t\leq l\,,
\end{equation}%
where \(\lambda(t)\) is a non-decreasing bounded real-valued
function on the interval \([0,l]\) which is determined by the
operator \(A\). The kernel \(K(t,s)\) of the integral operator
(\ref{MCP}) is of the form
\begin{equation}%
\label{KTM}%
K(t,s)=0\ \ \text{for}\ \ 0\leq s<t\leq l,\quad K(t,s)=i\,j\ \
\text{for}\ \ 0\leq t<s\leq l\,,
\end{equation}%
i.e., the operator (\ref{MCP}) can be considered as upper
triangular. The summand \(\lambda(t)\xi(t)\) corresponds to the
``main diagonal" of this operator.
 For the operator
\(T_{\text{cou}}\), the so called coupling operator, an explicit
formula can be obtained. Thus, the whole\footnote{The operator
\(T\) will not contain a discrete part \(T_{\text{dis}}\) if the
operator \(A\) has no eigenvalues. It will not contain a
continuous part \(T_{\text{con}}\) if \(l=0\).} operator \(T\) can
be naturally considered as an upper triangular operator.

\noindent \textsf{DEFINITION.} \textsl{Let \(A\) be an operator
which acts in a Hilbert space \(\mathcal{H}\). An operator
\(\widetilde{A}\) acting in a larger
 Hilbert space \(\widetilde{\,\mathcal{H}\,}\),
\(\widetilde{\mathcal{H}} \supseteq \mathcal{H}\), is said to be
an \textit{inessential extension} of the operator \(A\), if
\(\widetilde{A}=A\oplus S\), where \(S\) is a selfadjoint operator
acting in the space \(\widetilde{\mathcal{H}}\ominus
\mathcal{H}\).}

\noindent \textsf{\(\text{THEOREM\ I}^{\,\prime}\)}
(M.\,Liv\v{s}ic). \textsl{Let \(A\) be a bounded completely
non-selfadjoint linear operator  in a Hilbert space
\(\mathcal{H}\) such that  \(C_A\) is one-dimensional. Then there
exists an inessential extension \(\widetilde{A}:
\widetilde{\mathcal{H}}\to\widetilde{\mathcal{H}}\) of the
operator \(A\) that is unitarily equivalent to an ``upper
triangular" model operator \(T\) of the form
\textup{(\ref{TrForm})}: There exists a unitary operator \(U\)
acting from \(\mathcal{H}_{\text{mod}}=l^2\oplus L^2\) onto
\(\widetilde{\mathcal{H}}\) such that}
\begin{equation}%
\label{UEQ}%
T=U^*\widetilde{A}U=U^{-1}\widetilde{A}U\,.
\end{equation}%

Triangular models of the same general form
(\ref{TrForm})-(\ref{MDP})-(\ref{MCP}) can also be constructed for
bounded linear operators $A$ in a separable Hilbert space
\(\mathcal{H}\) when $C_A$ is only assumed to be of trace class.
They are, however, a bit more complicated.

For an operator \(A\) in a separable Hilbert space
\(\mathcal{H}\), let us introduce the \textsl{non-hermitian
subspace} \(\mathcal{N}_A\) as the closure of the image of its
imaginary part $C_A$:
\begin{equation}
\label{NHS} \mathcal{N}_A=\overline{C_A\,\mathcal{H}}\,.
\end{equation}
The dimension \(\mathfrak{n}_A\) of the non-hermitian subspace
\(\mathcal{N}_A\)  is said to be the \textsl{non-hermitian rank}
of the operator \(A\):
\begin{equation}
\label{NHR} \mathfrak{n}_A=\text{dim}\,\mathcal{N}_A\,.
\end{equation}
The restriction \(C_A|_{\mathcal{N}_A}\) of the operator \(C_A\)
on the subspace \(\mathcal{N}_A\), considered as an operator in
the Hilbert space \(\mathcal{N}_A\), is a selfadjoint operator for
which the point \(\{0\}\) is not an eigenvalue. Therefore, the
polar decomposition of this operator is of the form
\begin{equation}%
\label{PD}%
C_A|_{\mathcal{N}_A}=J_A\cdot M_A,
\end{equation}%
where
\begin{equation}%
\label{SignN}%
J_A:\mathcal{N}_A\to\mathcal{N}_A,\ \ J_A=J_A^{\,*},\ \ J_A^{\
2}=I_{\mathcal{N}_A} \quad \mbox{and} \quad
M_A:\mathcal{N}_A\to\mathcal{N}_A,\ M_A\geq 0.
\end{equation}%
(In this polar decomposition, \(J_A\) is the unitary operator and
\(M_A\) is the operator modulus.)

To construct the triangular model of the operator \(A\), let us
choose a Hilbert space \(\mathcal{E}\) of the same dimension as
the non-hermitian subspace \(\mathcal{N}_A\):
\(\text{dim}\,\mathcal{E}=\text{dim}\,\mathcal{N}_A\). Let
\(J_{\mathcal{E}}\) be an operator in \(\mathcal{E}\),
\begin{equation}%
\label{SignE}%
J_{\mathcal{E}}:\mathcal{E}\to\mathcal{E},\ \
J_{\mathcal{E}}=J_{\mathcal{E}}^{\,*},\ \ J_{\mathcal{E}}^{\ 2}
=I_{\mathcal{E}},
\end{equation}%
of the same signature\,\footnote%
{The spectrum of an operator \(J\) which posses the properties
\(J=J^*,\,J^2=I\) can consist of the points \(\{+1\}\) and
\(\{-1\}\) only. These points are eigenvalues of \(J\). Let
\(\mathfrak{p}_J\) and \(\mathfrak{q}_J\) denote the dimensions of
the corresponding eigenspaces, \(0\leq
\mathfrak{p}_J,\,\mathfrak{q}_J\leq \infty\).
 The signature of the
operator \(J\) is the pair
\((\mathfrak{p}_{J},\mathfrak{q}_{J})\).
}\ %
as that of  the operator \(J_A\). For the sake of brevity, we
shall restrict our attention
 to the case of operators with real spectrum only. The model space
\(\mathcal{H}_{\text{mod}}\) in this case is the space
\(L^2_{\mathcal{E}}([0,\,l])\), i.e., the space of all
square-integrable functions on a finite interval \([0,\,l]\subset
\mathbb{R}\), whose values are elements of the Hilbert space
\(\mathcal{E}\), with the scalar product:
\[
{\langle\boldsymbol{\xi},\boldsymbol{\eta}\rangle}_{\mathcal{H}_{\text{mod}}}=
\int\limits_0^l{\langle\xi(t),\eta(t)\rangle}_{\mathcal{E}}\,dt\,,
\ \ \text{for}\ \
\boldsymbol{\xi}=\xi(t),\,\boldsymbol{\eta}=\eta(t)\in
L^2_{\mathcal{E}}([0,\,l].
\]
The model operator \(T\) acts in the space
\(\mathcal{H}_{\text{mod}}\) according the rule
\begin{equation}%
\label{MOA}%
(T\xi)(t)=\lambda(t)\xi(t)+i\int\limits_t^l\Pi(t)J_{{}_\mathcal{E}}
\Pi(s)^*\xi(s)\,ds,
\end{equation}%
where  \(\lambda(t)\) is a non-decreasing real-valued function on
the interval \([0,\,l]\) and
 \(\Pi(t)\) is a function on the interval
\([0,\,l]\) whose values are Hilbert-Schmidt operators in
\(\mathcal{E}\)
 that satisfy the normalization condition
\begin{equation}%
\label{NorCo}%
{\text{trace}}_{\mathcal{{}_E}}\Pi(t)^*\Pi(t)\equiv 1\,,\quad
0\leq t\leq l\,.
\end{equation}%

\noindent \textsf{\(\text{THEOREM\ I}^{\,\prime\prime}\)}
(M.\,Liv\v{s}ic). \textsl{Let \(A\) be a bounded completely
non-selfadjoint linear operator  in a Hilbert space
\(\mathcal{H}\) such that $C_A$ is of  trace class and the
spectrum of $A$ is real.
 Then there exists an  inessential extension \(\widetilde{A}:
\widetilde{\mathcal{H}}\to\widetilde{\mathcal{H}}\),
\(\widetilde{\mathcal{H}}\supseteq\mathcal{H}\), of the operator
\(A\) that is unitarily equivalent to an ``upper triangular" model
operator \(T\) of the form \textup{(\ref{MOA})} : There exists a
unitary operator \(U\) acting from
\(\mathcal{H}_{\text{mod}}=L^2_{{}_\mathcal{E}}([0,\,l])\) onto
\(\widetilde{\mathcal{H}}\) such that}
\begin{equation}%
\label{UEq}%
T=U^*\widetilde{A}U=U^{-1}\widetilde{A}U\,.
\end{equation}%

\noindent%
\textsf{REMARK.} Direct calculation shows that
\begin{equation}%
\label{}%
\big((T-T^*)\xi\big)(t)=
i\int\limits_0^l\Pi(t)J_{{}_\mathcal{E}}\Pi(s)^*\,\xi(s)\,ds
\end{equation}%
Thus, the model operator (\ref{MOA})  is of the form
\begin{equation}%
\label{BLR1}%
(T\xi)(t)=\lambda(t)\xi(t)+
2i\int\limits_0^l\chi(t,s)H(t,s)\xi(s)\,ds\,,
\end{equation}%
 where
\[\chi(t,s)=1 \ \ \text{for}\ s>t,\ \ \chi(t,s)=0\ \
\text{for}\ \  s<t\,, \]
and  the kernel \(H(t,s)=\Pi(t)J_{{}_\mathcal{E}}\Pi(s)^*\)
represents the imaginary part of the operator \(T\):
\begin{equation}%
\label{BLR2}%
\big((T-T^*)\xi\big)(t)=2i\int\limits_0^lH(t,s)\xi(s)\,ds\,.
\end{equation}%
In other words, \textsl{the kernel \(K(t,s)\) that represents the
operator \(T\) can be obtained from the kernel \(H(t,s)\) that
represents the imaginary part of \(T\), by means of `` truncation
to the upper triangle": \(K(t,s)=\chi(t,s)H(t,s)\).
 }

 For operators whose imaginary part is of trace class but
whose spectrum is not necessary real, the triangular model has a
more complicated form; see e.g., \cite{Liv5} and in \cite{BroLi}.

Moshe Liv\v{s}ic introduced the machinery of characteristic
functions of linear operators in the mid forties in order to solve
a number of problems connected with theory of extensions of linear
operators, see \cite{Liv1} and \cite{Liv2}. He then applied this
machinery to establish
 the unitary equivalence of an operator of the class \(i\Omega\)
to a triangular model in the early fifties. See
\cite{Liv3},\cite{Liv4} for the first results and \cite{Liv5} for
a detailed presentation.

The characteristic function of a non-selfadjoint linear operator
\(A\) acting in a Hilbert space \(\mathcal{H}\) is defined as
follows: Choose a Hilbert space \(\mathcal{E}\) of the same
dimension as the non-hermitian subspace
\(\overline{C_A\,\mathcal{H}}\) of the operator \(A\) and then
factor the operator \(C_A\)  in the form
\begin{equation}%
\label{FactIn}%
C_A=\Gamma J_{{}_{\mathcal{E}}} \Gamma^*\,,
\end{equation}%
where \(\Gamma,\ J_{{}_{\mathcal{E}}}\) are linear operators,
\begin{equation}%
\label{FactOp}%
\Gamma:\mathcal{E}\to\mathcal{H}\,,\ \
J_{{}_{\mathcal{E}}}:\,\mathcal{E}\to\mathcal{E},\ \
J_{{}_{\mathcal{E}}}^2=I_{{}_{\mathcal{E}}} \ \
(I_{{}_{\mathcal{E}}}\text{ denotes the identity operator in
}\mathcal{E}).
\end{equation}%
\textsl{The characteristic function} \(W_A(z)\) of the operator
\(A\) is the operator valued function of the complex variable
\(z\) that is defined for \(z\) out of the spectrum of \(A\) by
the rule
\begin{equation}%
W_A(z)=I_{{}_{\mathcal{E}}}+2i\Gamma^*(zI-A)^{-1}\Gamma
J_{{}_{\mathcal{E}}}\,.
\label{CharF}%
\end{equation}%
Notice that $W_A(z)$ acts in the Hilbert space \(\mathcal{E}\),
which, in many problems of interest, is a finite dimensional
space. In Liv\v{s}ic's terminology,  the space \(\mathcal{E}\) is
said to be \textsl{the channel space} and the operator \(\Gamma\)
is said to be \textsl{the channel operator}.

 Liv\v{s}ic showed that the characteristic function is a unitary
invariant of a completely non-selfadjoint operator: \textsl{Let
\(A_1\) and \(A_2\) be two completely non-selfadjoint operators
such that
 their characteristic functions \(W_{A_1}(z)\) and
\(W_{A_2}(z)\) (with the same channel space \(\mathcal{E}\)) are
equal: \(W_{A_1}(z)\equiv W_{A_2}(z)\). Then the operators \(A_1\)
and \(A_2\) are unitarily equivalent}. To reduce a non-selfadjoint
operator to triangular form, Liv\v{s}ic calculated its
characteristic function \(W_{A}(z)\) and then constructed
 a model operator \(T\) in such a way that  its characteristic
function \(W_T(z)\) coincides with \(W_A(z)\).  Subsequently,
triangular models of operators were partially superceded by
functional models,  see \cite{SzNFo}, \cite{Bran}, \cite{NiVa}.

There is also a very important correspondence between the
invariant subspaces of an operator \(A\) and certain divisors of
its characteristic function. However, the importance of the notion
of a characteristic function is not confined to its applications
in  operator theory. Liv\v{s}ic related \textsl{the theory of
stationary linear dynamical systems} to the theory of linear
non-selfadjoint operators and showed that the characteristic
matrix function of a linear operator that serves as an ``inner
operator'' for the dynamical system can be identified with the
scattering matrix of this system. Examples are furnished in
\footnote{For more information on  the characteristic function of
a linear operator, see also the M.S.\,Liv\v{s}ic Anniversary
Volume \cite{OTSTR}, in particular, the Preface and the paper
\cite{Kats}.} \cite{Liv6}, \cite{Liv7} and \cite{BroLi}. A
detailed presentation of the early stages\,\footnote{ A more
elaborate presentation of  scattering theory for linear stationary
dynamical systems (with emphasis on applications to the wave
equation in \(\mathbb{R}^n\)) was carried out in \cite{LaPhi}.} of
the theory of open systems (as Liv\v{s}ic termed them) can be
found in \cite{Liv8}. In particular, as he noted in the first
sentence of Section 2.2 of that source: ``The resolution of a
system into a chain of elementary systems is closely related to
the reduction of the operator ... to triangular form."

The results of Liv\v{s}ic on reducing  operators to  triangular
form are similar in form to the result of Schur.  However, the
methods that he used are absolutely different from the method of
Schur. Schur's result implies that there exists an orthonormal
basis \(e_1,\,\dots\,,e_n\) of the space \(\mathbb{C}^n\) such
that the given matrix $A$ is upper-triangular in this basis. Thus,
if
\begin{equation}%
\label{CoChain}%
\mathcal{H}_0=0,\ \ \
\mathcal{H}_k=\text{span}\,\{e_1,\,\dots\,,e_k\}\,,\ \
k=1,\,2,\,\dots\,,n\,,
\end{equation}%
 then this collection \(\big\{\mathcal{H}_k\big\}_{0\leq k\leq
n}\) of subspaces of \(\mathbb{C}^n\) possesses the following
properties:
\begin{equation}%
\label{SchCh}%
\begin{array}{rl}
\text{i}.& 0=\mathcal{H}_0\subset \mathcal{H}_1\subset
\mathcal{H}_2\subset\
\dots \ \subset\mathcal{H}_n=\mathbb{C}^n\,,\\[1.3ex]
\text{ii}.&
\text{dim}\,\big(\mathcal{H}_k\ominus\mathcal{H}_{k-1}\big)=1\,,\\[1.3ex]
\text{iii}.& \text{Every subspace }\mathcal{H}_k \text{ is
invariant for the operator }A\,.
\end{array}
\end{equation}%
Conversely, let an operator \(A\) in \(\mathbb{C}^n\) and a
collection of subspaces \(\big\{\mathcal{H}_k\big\}_{0\leq k\leq
n}\) satisfying the conditions
(\ref{SchCh}.\,i)--(\ref{SchCh}.\,iii) be given, and let
 \(e_k\in
\mathcal{H}_k\ominus\mathcal{H}_{k-1}\) be unit vectors for \(
k=1,\,2,\,\dots\,,n\,\). Then the set
 of vectors \(\{e_k\}\) forms an orthonormal basis of \(\mathbb{C}^n\)
and the matrix of the operator \(A\) in this basis is
upper-triangular. It turns out that this strategy can be adapted to
obtain analogues of Schur's theorem in infinite dimensional
Hilbert spaces. The first step in this direction was taken
 by L.A. Sakhnovich \cite{Sakh1} who noticed that although the
proof of Schur is based on the fact that every operator $A$ in a
finite dimensional linear space over \(\mathbb{C}\) has an
eigenvector, that really the proof only depended upon following
property of the operator \(A\):

\textsf{Property I\,S\,.} \textsl{For every pair of closed
invariant subspaces \(\mathcal{H}_1\) and \(\mathcal{H}_2\) of the
operator \(A\) such that \(\mathcal{H}_1\subset\mathcal{H}_2\) and
\(\text{dim}\,(\mathcal{H}_2\ominus\mathcal{H}_1)>1\), there
exists a third closed invariant subspace \(\mathcal{H}_3\) of the
operator \(A\) such that
\(\mathcal{H}_1\subset\mathcal{H}_3\subset\mathcal{H}_2\,,\
\mathcal{H}_3\not=\mathcal{H}_1,\,\mathcal{H}_3\not=\mathcal{H}_2\).}

A theorem of J.\,von\,Neumann (unpublished) and of N. Aronszajn
and K.\,Smith \cite{ArSm}, guarantees that every compact operator
in a Hilbert space possesses the property \textsf{I\,S}. In
\cite{Sakh1}, Sakhnovich proved  that if the imaginary part
\(C_A\) of the operator \(A\) is of Hilbert-Schmidt class,
i.e., if  %
\begin{equation}%
\label{HSIP}%
 \sum\limits_k\big(s_k(C_A)\big)^2<\infty\,,
\end{equation} %
then the operator \(A\) possesses the  property \textsf{I\,S}\,. From
later results of V.I.\,Matsaev it follows that this condition can
be relaxed: if
\begin{equation}%
\label{somega}%
\sum\limits_{k}\frac{s_k(C_A)}{k+1}<\infty,
\end{equation}%
then \(A\) possesses the property \textsf{I\,S}.

Let \(\mathcal{H}\) be a Hilbert space and let
\(\mathfrak{P}_{\mathcal{H}}\) denote the collection of
orthoprojectors onto all possible closed subspaces of
\(\mathcal{H}\). The set \(\mathfrak{P}_{\mathcal{H}}\) is
partially ordered: $P_1\leq P_2$ if the corresponding ranges are
ordered by inclusion, i.e., if \(\mathcal{R}_{P_1}\subseteq
\mathcal{R}_{P_2}\)\,, and
 \(P_1< P_2\) if the inclusion of the ranges is proper.
A subset \(\mathfrak{P}\) of the set
\(\mathfrak{P}_{\mathcal{H}}\) that contains  at least two
orthoprojectors is said to be a \textsl{chain} if it is fully
ordered, i.e., if the conditions
\(P_1\in\mathfrak{P},\,P_2\in\mathfrak{P},\, P_1\not=P_2\) imply
that either \(P_1<P_2\)\,, or \(P_2<P_1\).

If a chain \(\mathfrak{P}\) contains orthoprojectors \(P^-\) and
\(P^+\) (\(P^{-} <  P^{+})\) such that every orthoprojector
\(P\in\mathfrak{P}\) distinct from them satisfies either the
inequality \(P<P^-\) or the inequality \(P>P^+\), then the pair
\((P^-, P^+)\) is said to be a \textsl{jump} in the chain
\(\mathfrak{P}\), and the dimension of the subspace
\(P^+\mathcal{H}\ominus P^-\mathcal{H}\) is said to be the
\textsl{dimension of the jump}. A chain without jumps is said to
be \textsl{continuous}.

The set of all chains in \(\mathcal{H}\) can be ordered by
inclusion:
 the chain \(\mathfrak{P}_1\) is said to
\textsl{precede} the chain \(\mathfrak{P}_2\) (and we write
\(\mathfrak{P}_1\prec \mathfrak{P}_2\)) if every orthoprojector
in
 \(\mathfrak{P}_1\) also lies in \(\mathfrak{P}_2\).
A chain \(\mathfrak{P}\) is said to be \textsl{maximal} with
respect to this ordering if there is no chain
\(\mathfrak{P}^{\prime}\) satisfying the conditions
\(\mathfrak{P}\prec\mathfrak{P}^{\prime},\
\mathfrak{P}\not=\mathfrak{P}^{\prime}\).

Let \(A\) be a bounded operator in a Hilbert space \(\mathcal{H}\)
and let \(\mathfrak{P}\) be a chain of orthoprojectors in
\(\mathcal{H}\). Then \textsl{the chain \(\mathfrak{P}\) is said
to be an  \textsl{eigenchain} for the operator \(A\)  if for every
\(P\in\mathfrak{P}\) the subspace \(P\mathcal{H}\) is invariant
under  \(A\), i.e.,  if the equality \(AP=PAP\) holds for every
\(P\in\mathfrak{P}\).}

The following result is established by transfinite induction in
\cite{Brod2} and  in \cite{Brod4}, where it appears as Theorem
15.2.

\noindent%
\textsf{THEOREM II} [M.S.\,Liv\v{s}ic-M.S.\,Brodski\u{\i};
L.A.\,Sakhnovich]. \textsl{Let  \(A\)  be a bounded linear
operator in a Hilbert space \(\mathcal{H}\) that satisfies  the
condition \textup{\textsf{IS}}. Then
there exists a maximal chain of orthoprojectors that is an
eigenchain for  \(A\).}

The idea for the proof of this theorem arose in a conversation
between M.S.\,Liv\v{s}ic and M.S.\,Brodskii. (See the historical
remark in the book \cite{Brod4}, p.\,278 of the Russian original
or p.\,234 of the English translation.) L.A.\,Sakhnovich gave an
independent  proof in
 \cite{Sakh1}.
 This theorem can be considered as a first step in extending the Schur 
theorem on
reducing  a matrix to  triangular form to the setting of a more general 
class of operators in Hilbert space. Based on it, Sakhnovich
obtained  the following result in \cite{Sakh1}:

\textsl{Every bounded linear operator \(A\) in a a separable
Hilbert space that satisfies the condition \textup{\textsf{I\,S}}
has an inessential extension \(\widetilde{A}\) which is unitarily
equivalent to an integral operator of the form
\begin{equation}%
\label{IOS1}%
x(t)\to
(\boldsymbol{K}x)(t)=\frac{d}{dt}\,\int\limits_t^1K(t,s)x(s)\,ds\,,
\end{equation}%
 acting in a space \(L^2_{{}_{\mathcal{E}}}([0,\,1])\)
 of vector functions \(x(t)\)  defined on the interval \([0,\,1]\)
whose values  belong to a Hilbert space \(\mathcal{E}\),
\(\textup{dim}\,\mathcal{E}\leq \infty\), provided with the scalar
product:
\begin{equation}%
\label{NSVF}%
{\langle\boldsymbol{x},\boldsymbol{y}\rangle}_{L^2_{\mathcal{E}}}=
\int\limits_0^1\langle x(t),\,y(t)\rangle_{{}_{\mathcal{E}}}\,dt,
\quad \textup{for}\ \  \boldsymbol{x}=x(t)\ \ \textup{and}\ \
\boldsymbol{y}=y(t)\,.
\end{equation}%
The kernel \(K(t,s)\) is a function defined for \(0\leq t,\,s\leq
1\) whose values are bounded linear operators acting in
\(\mathcal{E}\).}

A limitation of this last result is that  the class of
 kernels \(K(t,\,s)\) is not described. However, starting from
this theorem, Sakhnovich obtained  the following result in
\cite{Sakh2}:

\textsl{Every bounded operator \(A\) in a Hilbert space
\(\mathcal{H}\) whose
 spectrum is real and whose imaginary part \(C_A\) is of
Hilbert-Schmidt class has an inessential extension
\(\widetilde{A}\) which is unitarily equivalent  to an operator of
the form
\begin{equation}
\label{IOS2} x(t)\to H(t)x(t)+\int\limits_t^1K(t,s)x(s)\,ds\,,
\end{equation}
acting in the space \(L^2_{{}_{\mathcal{E}}}([0,1])\) of functions
defined on the interval \([0,\,1]\) whose values belong to a
Hilbert space \(\mathcal{E}\), provided with the scalar product
\textup{(\ref{NSVF})}. \(H(t)\) is a function defined on
\([0,\,1]\) whose values are bounded self-adjoint operators in
\(\mathcal{E}\): \(H(t)=H^*(t)\) for \(t\in[0,\,1]\). The kernel
\(K(t,s)\) is a function defined for \(0\leq t,s\leq 1\) whose
values are operators acting in the Hilbert space \(\mathcal{E}\)
that  are of Hilbert-Schmidt class. The kernel \(K\) satisfies the
condition
\begin{equation}%
\label{HSMK}%
\int\limits_0^1\textup{trace}_{{}_\mathcal{E}}\{(K^*K)(t,t)\}\,dt
<\infty\,.
\end{equation}%
}%

This result of Sakhnovich is on the  one hand more general than
the corresponding result of Liv\v{s}ic (because the condition
(\ref{HSIP}) is less restrictive than requiring $C_A$ to be of
trace class),  but on the other hand it is less concrete, since it
provides less information on the form of the kernel \(K\) than the
other theorem.

Further developments in this  area are related to the theory of
the abstract triangular representation of operators in a Hilbert
space by means of an integral with respect to a chain. This
integral appeared in the papers of M.S.\,Brodski\u{\i} at the end
of the fifties, \cite{Brod1}, \cite{Brod2}, \cite{Brod3}. In a
short time the theory of this  new integral and its applications
were  developed considerably. Important contributions to this
theory were made by V.I.\,Matsaev, \cite{Mats1}, and by
I.Ts.\,Gohberg and M.G.\,Krein, \cite{GoKr3}, \cite{GoKr4},
\cite{GoKr5}. The development of this  theory stimulated  new
analytic investigations of the spectral properties of  both
selfadjoint and non-selfadjoint operators.

To explain the definition of this integral, we begin with a
finite-dimensional example. Let \(\mathcal{H}\) be a complex
\(n\)-dimensional Hilbert space, \(n<\infty\). Let \(A\) be an
operator in \(\mathcal{H}\), and let \(\{e_k\}_{1\leq k\leq n}\)
be an orthonormal basis in \(\mathcal{H}\). Then the  operator
\(A\) can be written in the form
\begin{equation}%
\label{MaFor}%
A=\sum\limits_{j,k=1}^
ne_j\,a_{jk}\,\langle\,\cdot\,,e_k\rangle\,,
\end{equation}%
where \(a_{jk}=\langle A\,e_k,e_j\rangle\) are the entries of the
matrix of the operator \(A\) in this basis. Let the subspaces
\(\mathcal{H}_k\) be defined by (\ref{CoChain}), and let \(P_k\)
be the orthoprojector onto \(\mathcal{H}_k\). The collection
\(\mathfrak{P}\) of the orthoprojectors
\begin{equation}%
\label{FDCh}%
 0=P_0<P_1<\,\cdots\,<P_{n-1}<P_n=I
\end{equation}%
forms a maximal chain in \(\mathcal{H}\). This chain is an
eigenchain for \(A\).  Let
\begin{equation}%
\label{DeltaChain}%
\Delta P_k=P_k-P_{k-1},\quad k=1,\,2,\,\dots\,,n\,.
\end{equation}%
 Then, since
\begin{equation}%
\label{ES}%
P_k-P_{k-1}=e_k\,\langle \,\cdot\,,e_k\rangle \quad \mbox{and}
\quad
e_j\,a_{jk}\,\langle\,.\,,e_k\rangle=\Delta\,P_j\,\,A\,\,\Delta
P_k\,,
\end{equation}%
formula   (\ref{MaFor}) can be written in the form
\begin{equation}%
\label{DOIn}%
A=\sum\limits_{j,k=1}^n\Delta P_j\,A\,\Delta P_k\,.
\end{equation}
Moreover,   if \(a_{jk}=0\) for some choice of \(j,k\), then, by
(\ref{ES}), \(\Delta P_j\,A\Delta\,P_k = 0\)  in (\ref{DOI}).
Thus, if the matrix \(a_{jk}\) is upper triangular, i.e., if
\(a_{jk}=0\) for \(j>k\), and if \(a_{kk}=\lambda_k\)  (an
eigenvalue of the matrix \(A\)), then the representation
(\ref{DOIn}) takes the form
\begin{equation}%
\label{TDOI}%
A=\sum\limits_{k=1}^n\!\lambda_k\,\Delta P_k+\sum\limits_{k=2}^n
\sum_{j=1}^{k-1} \Delta P_j\,A\Delta\,P_k\,.
\end{equation}%
Since \(\sum\limits_{j=1}^{k-1}\Delta P_j=P_{k-1}\), (\ref{TDOI})
can be rewritten in the form
\begin{equation}%
\label{BOI}%
A=\sum\limits_{k=1}^n\!\lambda_k\,\Delta P_k+\sum\limits_{k=1}^n
P_{k-1}\,A\,\Delta P_k\,.
\end{equation}%
The first sum on the right hand side of (\ref{BOI}) represents the
``diagonal part" of \(A\), the second sum represents the
``super-diagonal" part with respect to the Schur basis
\(\{e_k\}_{1\leq k\leq n}\). (Everything here depends on the
choice of the basis.) Since
the matrix of the adjoint operator \(A^*\) (with respect to the
same orthonormal basis) is lower triangular, i.e., \(\Delta
P_j\,A^*\,\Delta P_k=0\) for \(j<k\), and \(P_{k-1}\,A^*\,P_k=0\),
the Schur result can be expressed as follows:

\textsl{For every operator \(A\) in a finite-dimensional Hilbert
space there exists at least one maximal eigenchain
\(\mathfrak{P}=\{P_k\}_{0\leq k\leq n}\). For every such
eigenchain, the operator \(A\) admits two representations:
\textup{(\ref{BOI})} and (with appropriate indexing) the
representation
\begin{equation}%
\label{TRIP}%
A=\sum\limits_{k=1}^n\lambda_k\Delta P_k+2i\sum\limits_{k=1}^n
P_{k-1}\,C_A\,\Delta P_k\,,
\end{equation}%
where \(\Delta P_k\) is defined by \textup{(\ref{DeltaChain})} and  %
\(C_A=\dfrac{A-A^*}{2i}\,. \) }\hfill \framebox[0.45em]{ }

The sums in (\ref{TRIP}) can be considered as  ``integrals" over
the chain \(\mathfrak{P}\):
\begin{equation}%
\label{BIR}%
A=\int\limits_{\mathfrak{P}}\lambda(P)\,dP+
2i\int\limits_{\mathfrak{P}}P\,C_A\,dP\,.
\end{equation}%
In the  case of the finite-dimensional \(\mathcal{H}\) that was
just discussed,
 the ``integrals"
in (\ref{BIR}) are no more than a notation for the finite sums in
(\ref{TRIP}). It is not a problem to generalize  integrals of the
form \(\int\limits_{\mathfrak{P}}\lambda(P)\,dP\) to the
infinite-dimensional case. This is the usual integral of a scalar
function with respect to an orthogonal spectral measure. Integrals
of this kind are  well understood, because of their connection
with needs of the theory of selfadjoint operators. However,
integrals of the form
\begin{equation}%
\label{GIaC}%
\mathfrak{I}(X,\,\mathfrak{P})\stackrel{\text{def}}{=}
\int\limits_{\mathfrak{P}}P\,X\,dP\,.
\end{equation}%
 for an arbitrary chain \(\mathfrak{P}\) of orthoprojectors and a
more or less general bounded linear operator \(X\) in an infinite
dimensional Hilbert space \(\mathcal{H}\) are more difficult to handle.%
\footnote{\  Integrals of scalar valued functions with respect to
operator valued measures and integrals of operator valued
functions with respect to a scalar valued measure are usually much
easier to deal with than integrals of operator valued functions
with respect to operator valued measures. In the integral
(\ref{GIaC}),  both the function \(PX\) and ``the measure" \(dP\)
are operator valued.}. An integral of the form (\ref{GIaC}) can be
\textsl{defined}   by means of a very natural limiting process
that was introduced by M.S.\,Brodski\u{\i} \footnote{\,An integral
of the form (\ref{GIaC}) can be considered (under  appropriate
parametrisation of the chain \(\mathfrak{P}\)\,) as a special case
of a double integral operator of the form
\begin{equation*}%
\label{GenDOI}%
\int\limits_{0}^{1}\!\!\int\limits_{0}^{1}
\chi(t,\,s)\,dP(t)\,X\,dP(s)\,,\
\ \text{with}\ \ \chi(t,\,s)=1\ \ \text{for}\ \ s>t,\ \chi(t,s)=0\
\ \text{for}\ \ s<t\,,
\end{equation*}%
We already met such integrals in  Section \ref{SchurPr}. However,
here the function \(\chi\) is of a very special form, and the
results which can be obtained for double operator integrals with
\textsl{this} function  are much more precise than the results
which follow from the general theory of double integral
operators.}: as usual, certain integral sums should be constructed
and then the passage to limit should be performed. The condition
\begin{equation}%
\label{ChCond}%
(P^+-P^-)X(P^+-P^-)=0  \ \ \text{for every jump }(P-P^+)\text{ of
the chain }\mathfrak{P}
\end{equation}%
is an evident necessary condition for the existence  of the
integral (\ref{GIaC}). However, the problem of obtaining
sufficient conditions for the existence of such an  integral
 turned out to be far more difficult. The theory of such
integrals, the so called \textsl{integral of triangular
truncation}, was created mainly in the works of
M.S.\,Brodski\u{\i}, I.Ts.\,Gohberg, M.G.\,Krein and V.I.\,Matsaev
and 
served to complete a program that was initiated by 
M.S.\,Liv\v{s}ic (see the remark
to \(\text{Theorem II}^{\prime\prime}\) of this section). A
detailed exposition of this theory is presented in \cite{Brod4},
 \cite{GoKr2} and  \cite{GoGoK}. Brodskii proved that under
condition (\ref{ChCond}), the integral (\ref{GIaC}) exists, if the
operator \(X\) is of trace class \(\mathfrak{S}_1\).
\footnote{\,Recall that singular values of a compact operator
\(X\) are the eigenvalues of the operator \(\sqrt{X^*X}\) indexed
in such a way that \(s_1(X)\geq s_2(X)\geq s_3(X)\geq\,\dots\,\)
and that \(X\in\mathfrak{S}_1\) if
\(\sum\limits_{k=1}^{\infty}s_k(X)<\infty\).} V.I.\,Matsaev,
\cite{Mats1}, sharpened this result. He proved, that \textsl{under
the condition \textup{(\ref{ChCond})}, the integral
\textup{(\ref{GIaC})} exists (in the sense of the convergence of
integral sums  with respect to the uniform operator norm), if the
compact operator \(X\) belongs to the class
\(\mathfrak{S}_{\omega}\), i.e., if the condition
\(\sum\limits_{1\leq k<\infty}s_k(X)\cdot k^{-1}<\infty\) holds}.
The latter result is precise in some sense. \textsl{If a compact
operator \(X\) does not belong to the class
\(\mathfrak{S}_{\omega}\), then there exists a continuous maximal
chain \(\mathfrak{P}\) such that the integral
\textup{(\ref{GIaC})} does not exist even in the sense of  weak
convergence}; see \cite{Brod4}, Lemma 22.2. In any case,
\textsl{if the operator \(X\) is compact and if the integral
\textup{(\ref{GIaC})} exists (in the sense of the convergence of
integral sums with respect to the uniform operator norm), then
this integral represents a Volterra operator}. We recall, that a
linear operator in a Hilbert space is said to be a
\textsl{Volterra operator} if it is compact and if its spectrum
consists of only one point,  the point zero.

The representation (\ref{BIR}) of an operator \(A\) by means of
the integral of triangular truncation   can be considered as a
coordinate-free representation of \(A\) from its maximal
eigenchain and its imaginary part. On the one hand, this
representation generalizes the results of Liv\v{s}ic (see
\(\text{Theorem I}^{\prime\prime}\) and the Remark following it
that focuses  attention on the formulas (\ref{BLR1}) and
(\ref{BLR2})). On the other hand, the representation (\ref{BIR})
is ``coordinate free", i.e., it represents the operator \(A\)
\textsl{itself} in the \textsl{original} Hilbert space
\(\mathcal{H}\), rather than a ``model" operator \(T\) that acts
in the ``model" space \(L^2_{{}_{\mathcal{E}}}\) and which is only
unitarily equivalent to the original operator \(A\) (or even to an
inessential extension \(\widetilde{A}\) of \(A\) acting in a
larger space \(\widetilde{\mathcal{H}}\supset\mathcal{H}\)). In
spirit, the representation (\ref{BIR}) is much closer to the
original  work of Schur \cite{Sch2b} than the triangular model
(\ref{MOA}) of Liv\v{s}ic. The integral representation (\ref{BIR})
for a bounded linear operator \(A\) with  imaginary part
\(C_A\in\mathfrak{S}_1\)   was first obtained  by Brodskii in
\cite{Brod1} using the representation (\ref{BLR1})-(\ref{BLR2}) as
a model. Brodski\u{\i} just transformed this representation to the
coordinate free form (\ref{BIR}). This proof used  the theory of
characteristic functions. Later, in \cite{Brod2} and \cite{Brod3},
the representation  (\ref{BIR}) was obtained for
 \textsl{arbitrary} Volterra operators \(A\) in a Hilbert space 
(in which case  \(\lambda(t)\equiv 0\) in (\ref{BIR})), and
also for bounded linear operators \(A\) with real spectrum and
\(C_A\in\mathfrak{S}_{\omega}\), independently of the theory of
characteristic functions,
 by  methods based on consideration of the eigenchains of the
operators \(A\), i.e., by  generalizing the reasoning of Issai
Schur.

The study of the integral of triangular truncation has  led to
unexpected and deep connections between the spectra of the real
and imaginary components of  Volterra operators. In certain cases
the clarification of these connections  has required the
development of new  analytic tools,  see Chapter III of the book
\cite{GoKr2}.  As an example of the application of the general
results obtained in the setting of the integral of triangular
truncation, we consider the Volterra operator \(A=B+iC, B=B^*,
C=C^*\) in the Hilbert space \(L^2([0,\,1])\) that is defined by the
equality
\[Ax(t)=2i\int\limits_t^1h(t-s)x(s)\,ds\,,\]
where the function \(h(\,\cdot\,)\) is periodic: \(h(t+1)=h(t)\),
Hermitian: \(h(-t)=\overline{h(t)}\), and  summable on
\([0,\,1]\). It is easily checked that the eigenvalues
\(\{\xi_j\}_{-\infty}^{\infty}\) and
\(\{\eta_j\}_{-\infty}^{\infty}\) of the operators \(B\) and \(C\)
(appropriately indexed) are related by the discrete Hilbert
transform:
\begin{equation}%
\label{DisHT}%
\eta_k=\frac{1}{\pi}\sum\limits_{l=
-\infty}^{\infty}\frac{\xi_l}{l-k+\frac{1}{2}}\,,\quad
-\infty<k<\infty\,.
\end{equation}
Consequently, it is possible to obtain estimates for the discrete
Hilbert transform by applying some results on   the spectra of the
Hermitian components of Volterra operators (\ref{DisHT}).

Thus, the Schur paper \cite{Sch2b}, which is elementary and purely
algebraic, stimulated the creation of several deep and rich
analytic theories.

\begin{minipage}{15.0cm}
\section{\hspace{-0.5cm}.\hspace{0.2cm}%
{\large \ Sequences of multipliers that preserve\newline%
the class of polynomials with only real zeros,\newline %
 and entire functions of the Laguerre-Polya-Schur class.\label{SchurPolya}}}
\end{minipage}\\[0.18cm]
\setcounter{equation}{0}

Problems related to the distribution of the zeros of polynomials
have attracted the attention of mathematicians for a long time. In
particular, the following question has generated considerable
interest. \textsl{How many real zeros does a given a polynomial
with real coefficients have?} There are several methods for either
estimating or determining precisely the number of zeros of such a
polynomial that belong to a given interval \((a,\,b)\) of the real
axis. These include the Descartes' rule of signs, the
Budan-Fourier algorithm, the Sturm algorithm and methods based on
Hermitian forms. These methods are presented in old books on
algebra (\cite{Web}, Vol.\,I, \cite{Kur}), as well as in books
devoted to the zeros of polynomials, \cite{Obr}, \cite{Mard},
\cite{Dieu}. The article \cite{KrNe} contains a detailed survey of
the method of Hermitian forms for the separation of the zeros of
polynomials. A lot of additional material on the distribution of
roots of polynomials can be found in \cite{PoSz}, Part V.

In this section we shall focus on a different class of results
that deal with transformations that preserve the class of
polynomials \(P(t)=p_0+p_1t+p_2t^2\,\,\cdots\,\,+p_nt^n\) with
real coefficients for which
\[
\#_{nr}(P) = \ \mbox{the number of non real roots of the
polynomial} \ P(t)
\]
is equal to zero. The simplest result of this kind states that if
\(P(t)\) is a polynomial with real coefficients and
\(\alpha\in\mathbb{R}\), then
\[
\#_{nr}(P) = 0 \Longrightarrow \#_{nr}(\alpha P +P^{\prime}) = 0.
\]
If the roots of \(P(t)\) are distinct, then this follows easily
from Rolle's theorem applied to \(e^{-\alpha t}P(t)\).

\textsf{THEOREM}
\begin{slshape}
Let \(P(t)=p_0+p_1t+p_2t^2\,\,\cdots\,\,+p_nt^n\)
 and \(Q(t)=q_0+q_1t+q_2t^2+\,\,\cdots\,\,+q_mt^m\) be polynomials with real
coefficients and assume that \(\#_{nr}(Q)=0\). Then the following
conclusions hold:

(1) [Hermite] \(\#_{nr}(Q(\frac{d}{dt})P)\leq \#_{nr}(P).\)

(2) [Laguerre] If also the roots of \(Q(t)\) fall outside the
interval \([0, n]\), then
\[
\#_{nr}(Q(0)p_0+Q(1)p_1t+Q(2)a_2t^2+\,\,\cdots\,\,+\,\,Q(n)p_nt^n)\leq
\#_{nr}(P).
\]

(3) [Malo] If also the roots of $Q(t)$ are either all positive or
all negative and \(l=\min(n,m)\), then %
\[\#_{nr}(P)=0
\Longrightarrow
 \#_{nr}(p_0q_0+p_1q_1t+\,\,\cdots\,\,+p_lq_lt^l)=0.\]

(4) [Schur] If also the roots of \(Q(t)\) are either all positive
or all negative and \(l=\min(n,m)\), then %
\[\#_{nr}(P)=0
\Longrightarrow
\#_{nr}(p_0q_0+1!p_1q_1t+2!p_2q_2t^2+\,\,\cdots\,\,+l!p_lq_lt^l)=0.\]
\end{slshape}
The cited results may be found in \cite{Obr}, \cite{Lagu1}),
\cite{Malo} and \cite{Sch7a}, respectively; see also \cite{PoSz}, Part\,%
V, Chapter 1,\,\S\,5, no.\,\textbf{63} and \textbf{67} for the
first two.

The last three  three statements of the theorem deal with an
operation of the form
\begin{equation}%
\label{GammaOper}%
p_0+p_1t+\,\,\cdots\,\,+p_nt^n\,\,\longrightarrow
 \,\,\gamma_0p_0+\gamma_1p_1t+\,\,\cdots\,\,+\gamma_np_nt^n.%
\end{equation}%
In the particular case considered by Schur, \(\gamma_k=q_k\) for
\(k\leq m\), and \(\gamma_k=0\) for \(k>m\), where the \(q_k\) are
obtained from the coefficients of a polynomial \(Q(t)\) with only
negative roots that we now write as
\(Q(t)=q_0+\frac{q_1}{1!}t+\frac{q_2}{2!}t^2+\,\,\cdots\,\,
+\frac{q_m}{m!}t^m\). The importance of this result is that it
admits a converse: Every sequences \(\{\gamma_k\}_{0\leq
k<\infty}\) for which  the operation (\ref{GammaOper}) preserves
the class of polynomials with real coefficients and $\#_{nr}(P)=0$
is either generated by a polynomial
\(Q(t)=q_0+\frac{q_1}{1!}t+\frac{q_2}{2!}t^2+\,\,\cdots\,\,+\frac{q_m}{m!}t^m\)
with only negative zeros, or belongs to the \textsl{closure} of
sequences generated by such polynomials. A full description of
this class of sequences \(\{\gamma_k\}_{0\leq k<\infty}\) is
presented in the paper \cite{PS} by G.\,Polya and I.\,Schur and
will now be described briefly below.

\begin{slshape}%
 Given an infinite unilateral sequence
\begin{equation}%
\label{MS}%
\Gamma=\{\gamma_0,\,\gamma_1,\,\gamma_2\,\,\,\dots\,\,\gamma_k,\,\,\dots\}%
\end{equation}%
of real numbers, let \(\Phi(t)\) denote the (formal) power series
\begin{equation}%
\label{FS}%
\Phi(t)=\sum\limits_{k=0}^{\infty}\frac{\gamma_k}{k!}\,t^k\,
\end{equation}%
and, for any polynomial 
\begin{equation}%
\label{EWPRL}%
P(t)=p_0+p_1t+\,\,\cdots\,\,+p_nt^n\,,
\end{equation}
let \(\Gamma[{P}(t)]\) denote the new polynomial:
\begin{equation}%
\label{TE}%
\Gamma[{P}(t)]=
\gamma_0a_0+\gamma_1a_1t+\gamma_2a_2t^2+\,\,\cdots\,\,+\gamma_na_nt^n\,.
\end{equation}%
\end{slshape}
\noindent%
\textsf{DEFINITION I} (\cite{PS}).\ \ \ %
\begin{slshape}%
\textup{\textsf{I.}}\,The sequence \textup{(\ref{MS})} is said to
be a \textsf{sequence of multipliers of the first type} if for
every polynomial 
\(P(t)\) with real coefficients (of
 arbitrary degree \(n\)),  $$\#_{nr}(P)=0 \Longrightarrow
\#_{nr}(\Gamma[P])=0.$$

\textup{\textsf{II.}} The sequence \textup{(\ref{MS})} is said to
be a \textsf{sequence of multipliers of the second type} if for
every polynomial 
\(P(t)\) with real coefficients (of
 arbitrary degree \(n\)) and only negative zeros,
\(\#_{nr}(\Gamma[P])=0.\)
\end{slshape}%

\noindent%
\textsf{DEFINITION II} (\cite{PS}). \ \  %
\begin{slshape}%
\textup{\textsf{I}.}\,  An entire function \(\Phi(t)\not\equiv 0\)
is an \textsf{entire function of the first type}, if it admits a
multiplicative representation of the form
\begin{equation}%
\label{EFFT}%
\Phi(t)=c\,t^l\,e^{\alpha
t}\prod\limits_k\left(1+t\,\delta_k\right)\,,
\end{equation}%
where \(c\not=0\) is a real number, \(l\) is a non-negative
integer, \(\alpha\) is a non-negative real number, and the
\(\delta_k\) are non-negative numbers that satisfy the condition
\(\sum\limits_k \delta_k<\infty\).

\textup{\textsf{II}.} An entire function \(\Phi(t)\not\equiv 0\)
is an \textsf{entire function of the second type}, if it admits a
multiplicative representation of the form
\begin{equation}%
\label{EFST}%
\Phi(t)=c\,t^l\,e^{-\beta t^2+\alpha
t}\prod\limits_k\left(1+t\,\delta_k\right)\,e^{-\delta_kt}\,,
\end{equation}%
where \(c\not=0\) is a real number, \(l\) is a non-negative integer, %
\(\beta\) is a non-negative number, \(\alpha\) is a real number,
and the \(\delta_k\) are real numbers that satisfy the condition
\(\sum\limits_k (\delta_k)^2<\infty\).
\end{slshape}%

\noindent \textsf{THEOREM} (G.\,Polya and I.\,Schur, \cite{PS}).\ \ \ %
\begin{slshape}%
 \textup{\textsf{I}.} If the sequence
  \(\gamma_0,\,\gamma_1,\,\,\dots\,\,,\gamma_k,\,\,\dots\,\,\)
is a sequence of multipliers of the first (respectively the
second) type, then the series \textup{(\ref{FS})} converges in the
whole complex plane, and the entire function \(\Phi(t)\) which is
represented by this series is an entire function of the first
(respectively  the second) type.

\textup{\textsf{II}.} If \(\Phi(t)\) is an entire function of the
first (respectively the second) type, and \textup{(\ref{FS})} is
its Taylor expansion, then the sequence
\(\gamma_0,\,\gamma_1,\,\,\dots\,\,,\,\gamma_k,\,\,\dots\,\,\) is
a sequence of multipliers of the first (respectively the second)
type.
\end{slshape}%

This theorem gives a full description of the sequences of
multipliers of both  the first  and  second type. The appearance
of two types of multipliers (and  two types of entire functions)
corresponds to the fact that in the  Schur theorem from
\cite{Sch7a} that was stated above, the polynomials \(P(t)\) and
\(Q(t)\) appear in a \textsl{symmetric way}: if \textsl{one} of
the polynomials \(P(t)=p_0+p_1t+p_2t^2+\,\,\cdots\ \ \) or
\(Q(t)=q_0+q_1+q_2t^2+\,\,\cdots\ \ \) has only real zeros, and
the \textsl{other} has only negative zeros, then all the zeros of
the polynomial \(p_0q_0+1!p_1q_1t+2!p_2q_2t^2+\,\,\cdots\) are
real. Thus, roughly speaking, sequences of the first (respectively
second) type act on polynomials that are entire functions of the
second (respectively first) type. The two types of entire
functions arise as limits of the two classes of polynomials:

\noindent%
\textsf{THEOREM} (E.\,Laguerre,\,\cite{Lagu2}; G.P\'{o}lya,\,%
\cite{Pol1}).\ \ \\[1.3ex]%
\begin{slshape}
 \textup{I.}\, Let \(\Phi(t)\) be an entire function of the first
type (respectively the second type). Then there exists a sequence
\(\{\Phi_n(t)\}_{n=1,\,2,\,\,\dots}\) of polynomials such that the
zeros of \(\{\Phi_n(t)\}\) lie in the negative half-axis
(respectively the real axis)
 and \(\{\Phi_n(t)\}\)
converges to \(\Phi(t)\) locally uniformly in the whole complex
plane.
\\[1.3ex]
\textup{II.} If a sequence of  polynomials
\(\{\Phi_n(t)\}_{n=1,\,2,\,\,\dots}\) converges uniformly in a
neighborhood of the origin to a function that is not identically
equal to zero  and if all the zeros of every polynomial \(\Phi_n\)
lie in the negative half-axis (respectively the real axis), then
the sequence \(\{\Phi_n\}_{n=1,\,2,\,\,\dots}\) converges locally
uniformly in the whole complex plane and the limit function
\(\Phi(t)\) is an entire function of the first  (respectively
second) type.
\end{slshape}

Part  I of this theorem was obtained by Laguerre, \cite{Lagu2};
part II was obtained by Polya, \cite{Pol1}. Laguerre obtained a
weak version of part II. Namely, he assumed that the sequences of
polynomials \(\{\Phi_n\}\) considered above
 converge in the whole complex plane, not just
in a neighborhood of the origin, and deduced  the same properties
of the limiting function that are stated in  part II of the
preceding theorem. This result of Laguerre is not strong enough to
obtain a description of the multiplier sequences, the  stronger
result by P\'olya is needed. The theorems of Laguerre and Polya,
and some  generalizations, can be found  in \cite{HiWi}, Chapter
III, \S\,3, and in \cite{Levi}, Chapter VIII.

The paper \cite{PS} by P{\'o}lya and Schur served as a source of
inspiration for the investigations of I.J.\,Schoenberg on the
representation of totally positive functions and sequences. The
notion of total positivity was introduced by Schoenberg in
\cite{Scho1}.

\noindent%
\textsf{DEFINITION III.} \ \ %
\begin{slshape}%
 A real function (or, in other terms,  kernel) \(K(t,s)\) of two
 variables ranging over linearly ordered sets \(\mathcal{T}\)
 and \(\mathcal{S}\), respectively, is said to be \textsf{totally
 positive} if for every\footnote{If both sets \(\mathcal{T}\) and
\(\mathcal{S}\)
 are infinite, then \(m\) can be an arbitrary natural number; if at
 least one of the sets \(\mathcal{T}\) or \(\mathcal{S}\) is finite
then \(m\) can be  an arbitrary natural number satisfying the
restriction \(m\leq\min\{|\mathcal{T}|, |\mathcal{S}|\}\),
 where \(|\mathcal{M}|\) denotes the \textsf{cardinality} of the set
\(\mathcal{M}\).}%
 \ \  \(m\) and for every
\begin{equation}%
 \label{Rang}%
 t_1<t_2<\,\,\cdots\,\,<t_m,\ \ \ s_1<s_2<\,\,\cdots\,\,<s_m\quad
 t_i\in\mathcal{T},\,t_j\in\mathcal{S},
 \end{equation}%
 the inequalities
 \begin{equation}%
 \label{PosDet}%
 K%
\left(%
 \begin{array}{c}%
 t_1,\,t_2,\,\,\dots\,\,,t_m\\[1.0ex]
 s_1,\,s_2,\,\,\dots\,\,,s_m
 \end{array}%
\right)%
\geq 0
 \end{equation}%
 hold, where
 \begin{equation}%
 \label{NotMin}%
K%
\left(%
 \begin{array}{c}%
 t_1,\,t_2,\,\,\dots\,\,,t_m\\[1.0ex]
 s_1,\,s_2,\,\,\dots\,\,,s_m
 \end{array}%
\right)%
=\det%
\left[%
\begin{array}{cccc}%
K(t_1,\,s_1)&K(t_1,\,s_2)&\,\,\cdots\,\,&K(t_1,\,s_m)\\[0.8ex]
K(t_2,\,s_1)&K(t_2,\,s_2)&\,\,\cdots\,\,&K(t_2,\,s_m)\\ [0.8ex]
\vdots      &    \vdots  &              &   \vdots   \\ [0.8ex]
K(t_m,\,s_1)&K(t_m,\,s_2)&\,\,\cdots\,\,&K(t_m,\,s_m)\\
\end{array}%
\right].%
\end{equation}%
\end{slshape}%

Usually \(\mathcal{T}\) and \(\mathcal{S}\) are either
subintervals of the real axis (that may coincide with the full
axis), or countable sets of real numbers such as the set of all
integers or the set of all non-negative integers, or even finite
sets of integers. If \(\mathcal{T}\) and \(\mathcal{S}\) are sets
of integers, then \(K\) can be viewed as a matrix; in this case,
\(K\) is referred to as a \textsl{totally positive matrix}.

A concept that is more general than total positivity is
\textsl{sign regularity}.

\begin{slshape}%
A function \(K(t,s)\) is said to be \textsf{sign regular} if there
exists a sequence of numbers \(\varepsilon_m\), each of which is
equal to either \(+1\) or \(-1\), such that in the setting of
(\ref{Rang}),  the inequalities
\begin{equation}%
\label{SReg}%
\varepsilon_m\,\,K\!%
\left(%
 \begin{array}{c}%
 t_1,\,t_2,\,\,\dots\,\,,t_m\\[1.0ex]
 s_1,\,s_2,\,\,\dots\,\,,s_m
 \end{array}%
\right)%
\geq 0
\end{equation}%
hold.
\end{slshape}

Totally positive matrices (and kernels) have very interesting
spectral properties that were  discovered  by F.R.\,Gantmacher and
M.G.\,Krein, \cite{GaKr1}, \cite{GaKr2}, \cite{GaKr3}. All the
eigenvalues of a totally positive matrix are positive and
distinct\footnote{Under the assumption that \textsl{all} its
minors are strictly positive}. Moreover,  its eigenvectors posses
oscillatory properties that are analogous to the oscillatory
properties of the eigenfunctions of  Sturm-Liouville differential
equations. A presentation of the spectral properties of totally
positive matrices and kernels can also be found in the survey article
\cite{Pink} by A.\, Pinkus. However, the notion of total positivity was
introduced by Schoenberg  \cite{Scho1}  in his study of
\textsl{variation-diminishing} kernels. Strictly speaking, in
\cite{Scho1}, the definitions of total positivity and sign
regularity  were formulated for the case of \textsl{finite
matrices}; generalizations to wider settings were developed  later
by Schoenberg himself and by S.\,Karlin. (See the book \cite{Karl}
for the references and for the history.)

Let \({\mathcal{V}}[z_1,\,z_2,\,\,\dots\,\,,\,z_l]\) denote the
number of sign changes of a given sequence
\([z_1,\,z_2,\,\,\dots\,\,,\,z_l]\) of real numbers, when the zero
terms are discarded. For example,
\({\mathcal{V}}[1,\,0,\,1,\,0,\,-1]=1\) and
\({\mathcal{V}}[1,\,-1,\,1,\,1,\,-1,\,1]=3\).

\noindent%
\textsf{DEFINITION IV.}\ \
\begin{slshape}%
Let \(K=\big[k_{ij}\big]_{}\) be a \(p\times q\) matrix with real
entries \(k_{ij},\,1\leq i\leq p,\,1\leq j\leq q;\,
p,\,q<\infty.\) The matrix \(K\) is said to be
\textsf{variation-diminishing}, if for every sequence
\(x=[x_1,\,x_2,\,\,\dots\,\,,\,x_q]\) of real numbers, the
sequence
\begin{equation}%
\label{MVDT}%
 y_i=\sum\limits_{1\leq j\leq q}k_{ij}x_j\,,\quad (1\leq i\leq
 p)\,,
\end{equation}%
enjoys the property
\begin{equation}%
\label{VDP}%
\mathcal{V}[y_1,\,y_2,\,\,\dots\,\,,\,y_p]\leq
\mathcal{V}[x_1,\,x_2,\,\,\dots\,\,,\,x_q]\,.
\end{equation}%
\end{slshape}

\noindent%
\textsf{THEOREM} (I.J.\,Schoenberg, \cite{Scho1}).\ \
\begin{slshape}%
Let \(K\) be a \(p\times q\) matrix with real entries.\\[-3.5ex]

\textup{\,I.} If the matrix \(K\) is sign-regular (in particular,
if \(K\) is totally positive), then \(K\) is
variation-diminishing.

\textup{II.} If the matrix \(K\) is variation-diminishing, and
\begin{equation}%
\label{ROR}%
\text{\textup{rank}}\, K=q\,,
\end{equation}%
then \(K\) is sign-regular.
\end{slshape}

Under the additional restriction (\ref{ROR}), this theorem gives
necessary and sufficient conditions for a \(p\times q\) matrix
\(K\) with real entries to be variation-diminishing.
A characterization of variation-diminishing matrices without any
restrictions was obtained by Th.\,Motzkin in his PhD Thesis
(Basel, 1934). His thesis was published in 1936, \cite{Motz1}; see
also \cite{Motz2}. Additional  characterizations of matrix
variation diminishing transforms can  also be found in  Chapter 4,
\S\,8 of \cite{HiWi}, and  Chapter 5, \S\S\, 1,\,2 of \cite{Karl}.
The latter  is a storehouse of wisdom on  total positivity,
variation diminishing transformations and related issues and
applications.

For a function \(x(t)\) which is defined on a linearly ordered set
\(\mathcal{T}\), the number of sign changes
\({\mathcal{V}}[x(t)]\) is defined as \({\mathcal{V}}[x(t)]=\sup
{\mathcal{V}}[x(t_1), \,x(t_2),\,\,\dots\,\,,\,x(t_l)]\), where
the \textit{supremum} is taken over all
\(t_1,\,t_2,\,\,\dots\,\,,\,t_l\) from \(\mathcal{T}\) such that
\(t_1<t_2<\,\,\dots\,\,<t_l\). (It is possible that
\({\mathcal{V}}[x(t)]=\infty\)).

For a real-valued kernel \(K(t,\,s)\) defined for
\(t\in{\mathcal{T},\,s\in\mathcal{S}}\), where \(\mathcal{T}\) and
\(\mathcal{S}\) are subintervals (finite or infinite) of the real
axis, the variation diminishing property also can be formulated in
the form
\begin{equation}%
\label{FVDP}%
{\mathcal{V}}[y(t)]\leq {\mathcal{V}}[x(s)]\quad
(t\in{\mathcal{T}},\, s\in{\mathcal{S}}),
\end{equation}%
where
\begin{equation}%
\label{IntTr}%
y(t)=\int\limits_{c}^{d}K(t,s)\,x(s)\,ds\,,\quad a\leq t\leq b\,.
\end{equation}%
Of course some restrictions have to be imposed on the class of
functions \(x(s)\) and on the kernel \(K(t,\,s)\) to ensure the
existence of the transformation (\ref{IntTr}) and the possibility
of counting uniquely\footnote{\,Since the number of  changes of
sign of a function is defined \textit{pointwise}, the functions
\(x(s)\) and \(y(t)\) must be defined \textit{everywhere},  not
just \textit{almost everywhere} on the appropriate intervals.} the
number of changes of sign of the functions \(x(s)\) and \(y(t)\).

The  preceding theorem of Schoenberg  that characterizes matrices
with variation-diminishing properties in terms of their
sign-regularity, can be extended to continuous kernels
\(K(t,s)\),.

It should be mentioned that as early as  1912, in order to
estimate the number of real zeros of polynomials with real
coefficients, M.\,Fekete considered formal power series
with real coefficients%
\(%
\label{Fek}%
\sum\limits_{0\leq i<\infty}c_it^i
\) %
that
possess the following property: %
\textsl{for a given natural number \(r\), all the determinants
(\ref{TPDM}) with \(m=1,\,2,\,\,\dots\,\,,r\) are non-negative.}
Power series with this property are called %
\textsf{\textsl{r}-time positive}. Multiplication by such a power series %
\(\sum\limits_{0\leq i<\infty}c_it^i\), %
transforms the power series \(\sum\limits_{0\leq k\leq r}x_kt^k\)
into the power series \(\sum\limits_{0\leq j<\infty}y_jt^j\)
according to the rule
\[\sum\limits_{0\leq j<\infty}y_jt^j=
\bigg(\sum\limits_{0\leq
i<\infty}c_it^i\bigg)\cdot%
\bigg(\sum\limits_{0\leq k<\infty}x_kt^k\bigg)\,,
\]
or, equivalently,
\[y_j=\sum\limits_{0\leq k\leq j}c_{j-k}x_k\,,\quad 0\leq j<\infty.\]
Fekete formulated the following statement (see footnote number six
in \cite{Fek}):

\begin{slshape}
Let \(\sum\limits_{0\leq i<\infty}c_it^i\) be an \(r\)-time
positive (formal) power series, and let \(\sum\limits_{0\leq k
<\infty}x_kt^k\) be a polynomial of  degree \(r\) (i.e., \(x_k=0\)
for \(k>r\)) with real coefficients. Then
\[{\mathcal{V}}[y_0,\,y_1,\,y_2,\,\,\dots\,\,,y_j,\,\,\dots]\leq
{\mathcal{V}}[x_0,\,x_1,\,\,\dots\,\,,\,x_r].
\]
\end{slshape}

Totally positive matrices and kernels that depend on the
difference of their arguments are of special interest. Let
\(\{c_i\}_{-\infty<i<\infty}\) be  an infinite bilateral sequence
and let   the infinite Toeplitz matrix \(K\) be defined in terms
of this sequence by the rule
\begin{equation}%
\label{ITM}%
K=\big[k_{p,q}\big]_{0\leq p,q<\infty}\,,\quad k_{p,q}%
\stackrel{\textup{def}}{=}c_{p-q}.
\end{equation}%
If the sequence \(c_i\) is unilateral: \(c_i\) are defined only
for \(i\geq 0\), we first extend the original sequence to the set
of all  integers by  setting \(c_i\stackrel{\textup{def}}{=}0\)
for \(i<0\) and then define the Toeplitz matrix \(K\) by  rule
(\ref{ITM}) applied to  the extended sequence.

\noindent%
\textsf{DEFINITION V} (\,I.J.\,Schoenberg, \cite{Scho2}\,).
\begin{slshape}%
The real valued infinite sequence \(\{c_i\}\), bilateral or
unilateral, is said to be \textsf{totally positive} if the matrix
\textup{(\ref{ITM})} is totally positive, i.e., if for every
natural \(m\) and for every choice of integers
\(p_1<p_2<\,\,\dots\,\,<p_m\), \(q_1<q_2<\,\,\dots\,\,<q_m\) the
inequality
\begin{equation}%
\label{TPDM}%
\det%
\left[%
\begin{array}{cccc}
c_{p_1-q_1}&c_{p_1-q_2}&\,\cdots\,& c_{p_1-q_m}\\[0.5ex]
c_{p_2-q_1}&c_{p_2-q_2}&\,\cdots\,& c_{p_2-q_m}\\[0.5ex]
\cdot&\cdot&\cdot&\cdot\\[0.5ex]
c_{p_m-q_1}&c_{p_m-q_2}&\,\cdots\,& c_{p_m-q_m}
\end{array}%
\right]%
\geq 0
\end{equation}%
holds.
\end{slshape}%

\noindent%
\textsf{DEFINITION VI} (\,I.J.\,Schoenberg, \cite{Scho2}\,).
\begin{slshape}
A real valued function \(\Lambda(t)\) that is  defined for all \ %
\(t\in(-\infty,\,\infty)\) is said to be \textsf{totally positive}
if it satisfies the following three conditions:
\begin{enumerate}
\item[\textup{i.}]
The kernel \(K(t,\,s)\stackrel{\textup{def}}{=}\Lambda(t-s),\ \
-\infty<t,\,s<\infty\), is totally positive, i.e., for every
natural number \(m\) and for every \(t_1<t_2<\,\,\cdots\,\,<t_m\),
\(s_1<s_2<\,\,\cdots\,\,<s_m\) the inequality
\begin{equation}%
\label{TPDF}%
\det%
\left[%
\begin{array}{cccc}
\Lambda(t_1-s_1)&\Lambda(t_1-s_2)&\,\cdots\,& \Lambda(t_1-s_m)\\[0.5ex]
\Lambda(t_2-s_1)&\Lambda(t_2-s_2)&\,\cdots\,& \Lambda(t_2-s_m)\\[0.5ex]
\cdot&\cdot&\cdot&\cdot\\[0.5ex]
\Lambda(t_m-s_1)&\Lambda(t_m-s_2)&\,\cdots\,& \Lambda(t_m-s_m)
\end{array}%
\right]%
\geq 0
\end{equation}%
holds.
\item[\textup{ii.}] The function \(\Lambda(t)\) is measurable.
\item[\textup{iii.}] The function \(\Lambda(t)\) is positive
for at least two distinct values of \(t\).
\end{enumerate}
\vspace*{-4.0ex}
\hfill \framebox[0.45em]{ }%
\end{slshape}

It is not difficult to prove that if a function \(\Lambda(t)\) is
defined on \(\mathbb{R}\) and is nonnegative there, and if the
inequalities (\ref{TPDF}) hold for \(m=2\) and for all
\(t_1<t_2,\,s_1<s_2\), then the function
\(\psi(t)=-\ln\Lambda(t)\) is convex (in the wide
sense\,%
\footnote{\,A function defined on \(\mathbb{R}\)is said to be
convex \textsl{in the wide sense} if it is convex in the usual
sense on some subinterval of \(\mathbb{R}\), which can coincide
with \(\mathbb{R}\), can be finite or semi-infinite. and is equal
to \(+\infty\) on the complement of this interval. A  non-negative
function \(\Lambda(t)\) on \(\mathbb{R}\) is convex in the wide
sense if and only if the inequalities
(\ref{TPDF}) hold for \(m=2\) and for all \(t_1<t_2,\,s_1<s_2\).}%
) on \(\mathbb{R}\). In particular, \textsl{if a function
\(\Lambda(t)\) is totally positive, then the function
\(-\ln\Lambda(t)\) is convex (in the wide sense) on
\(\mathbb{R}\)}. Therefore, \textsl{for every totally positive
function \(\Lambda(t)\)  the limits
\begin{equation}%
\label{Lim}%
\alpha=\lim\limits_{t\to-\infty}\frac{-\ln\Lambda(t)}{t}\,,\ \ %
\beta=\lim\limits_{t\to+\infty}\frac{-\ln\Lambda(t)}{t}\,,
\end{equation}%
exist and \(\ -\infty\leq\alpha\leq\beta\leq\infty\)\,.} The
equality \(\alpha=\beta\) holds if and only if \(\Lambda(t)\) is
of the form
\begin{equation}%
\label{TrExTPF}%
\Lambda(t)=e^{kt+l}\,,\quad -\infty<t<\infty\,, \text{for  some
real constants}\quad k\text{\ and\ }l.
\end{equation}%
(A function of the form (\ref{TrExTPF}) is easily seen to be
totally positive since all the determinants (\ref{TPDF}) vanish,
if \(m\geq 2\).) Thus, \textsl{if a function \(\Lambda(t)\) is
totally positive, but not of the form (\ref{TrExTPF}), then
\(\alpha<\beta\)}\,
and hence the two-sided Laplace transform %
\(\int\limits_{\mathbb{R}}\Lambda(t)e^{zt}dt\) %
exists for all  points \(z\) in the open strip %
\(\,\alpha<\text{Re}\,z<\beta\,\) %
of the complex \(z\)-plane and represents a holomorphic function
there. Moreover, the function, represented by this Laplace
transform, takes
\textsf{strictly positive} values for %
\(z\in(\alpha,\,\beta)\subset\mathbb{R}.\) %
Hence, the reciprocal function \(\Psi(z)\):
\begin{equation}%
\label{LapTr}%
\Psi(z)\stackrel{\textup{def}}{=}%
\bigg(\int\limits_{\mathbb{R}}\Lambda(t)\,e^{zt}dt\bigg)^{-1}\,,\quad
\alpha<\textup{Re}\,z<\beta\,,
\end{equation}%
is meromorphic in  the strip \(\alpha<\text{Re}\,z<\beta\),
holomorphic in all points of the interval
\((\alpha,\,\beta)\in\mathbb{R}\)\,, and takes strictly positive
finite values in this interval: \ \(0<\Psi(x)<\infty\)\,,\
\(\alpha<x<\beta\,.\)

\noindent%
\textsf{THEOREM I} (I.J.\,Schoenberg,
 Theorem 1 in \cite{Scho2}).\\[0.5ex]
\begin{slshape}%
\hspace*{2.0ex}\textup{I.} Let \(\Lambda(t)\) be a totally
positive function that is not of the form \textup{(\ref{TrExTPF})}
and let the function \(\Psi(z)\) be defined by means of
\textup{(\ref{LapTr})} as a meromorphic function in the strip
\(\,\,\alpha<\textup{Re}\,z<\beta\,\) \textup{(}see
\textup{(\ref{Lim})}\textup{)}.

Then  \(\Psi(z)\) is holomorphic in this strip and  admits an
analytic continuation to the whole complex plane \(\mathbb{C}\).
The continued function (denoted by \(\Psi(z)\) as well) is an
entire function of the second type,\,%
\footnote{\,In the sense of the paper \textup{\cite{PS}}, see
Definition II above in this section.}%
which is not of the form%
\begin{equation}%
\label{ExepEF}%
\Psi(z)=ce^{az}, \quad \textup{\small where} \, a \,
\textup{\small and} \,c  \textup{ \small are  real constants, }
c\not=0\,.
\end{equation}%
The function \(\Lambda(t)\) can be recovered from \(\Psi(z)\) by
mean of the inversion formula
\begin{equation}%
\label{InvLapl}%
\Lambda(t)=\frac{1}{2\pi
i}\int\limits_{\gamma-i\infty}^{\gamma+i\infty}%
\frac{1}{\Psi(z)}\,e^{-zt}dz\,, \quad -\infty<t<\infty\,,
\end{equation}%
where \(\gamma\) is an arbitrary\,%
\footnote{\,The value of the integral in \textup{(\ref{InvLapl})}
does not depend on the choice of
\(\gamma\in(\alpha,\,\beta)\).} %
real number from \((\alpha,\,\beta)\).

\textup{II.} Let \(\Psi(z)\) be an entire function of the second
type that is not of the form \textup{(\ref{ExepEF})}, and let
\(\Psi(x)\) be strictly positive on an interval
\((\alpha,\,\beta)\in\mathbb{R}\) (so that the reciprocal function
\(\dfrac{1}{\Psi(z)}\) is holomorphic in the vertical strip
\,\(\alpha<\textup{Re}\,z<\beta\)). Let the function
\(\Lambda(t)\) be \textup{defined} from this \(\Psi(z)\) by means
of the integral.\,%
\footnote{The integral (\ref{InvLapl}) converges absolutely for
every entire function \(\Psi(z)\) of the second type, except when
\(\Psi(z)\) is of the form \(ce^{at}(1+\delta z)\), where \(c\not=
0,\,\delta\not= 0\), \(c,\,\delta, a\) are real, in which case the
integral \textup{(\ref{InvLapl})} converges in the sense of
principal values.} %
Then the function \(\Lambda(t)\) is totally positive, and if the
interval \((\alpha,\,\beta)\) is the maximal
interval on which the function \(\Psi(x)\) is positive,\,%
\footnote{\,That is, if either \(\Psi(\alpha)=0\), or
\(\alpha=-\infty\), and if either \(\Psi(\beta)=0\), or
\(\beta=\infty\)\,.}%
then the endpoints \(\alpha\) and \(\beta\) of this interval
coincide with the limits \(\alpha\) and \(\beta\) in
\textup{(\ref{Lim})}.
\end{slshape}%

The proof of this theorem is based essentially on methods and
results from the paper \cite{PS}. Indeed, the names  of Polya and
Schur (and Laguerre) appear in the title of \cite{Scho2}, and as
Schoenberg himself writes \textsl{``A proof of Theorem 1 is
essentially based on the results and methods developed by Polya
and Schur. The only additional element required is a set of
sufficient conditions insuring that a linear transformation be
variation diminishing."}

For the sake of added perspective, we shall sketch the proof of
part II, which is not difficult (once the theorem has been
formulated), but shall omit the proof of part I, which is not so
simple and straightforward. If the function \(\Psi(z)\) is ``a
linear factor", i.e.,  if \(\Psi(z)=(1+\delta z),\, \delta\not=0\)
when  \(1+\delta\gamma>0\) and \(\Psi(z)=-(1+\delta z),\,
\delta\not=0\) when  \(1+\delta\gamma<0\),
then 
\begin{equation}%
\label{LinFact}%
\Lambda(t)=\left\{%
\begin{array}{ll}%
\delta^{-1}e^{t/\delta}, & t<0\\
0\,,&t>0
\end{array}%
\right.,\, \text{if\ }\gamma>-1/\delta\,;
\quad%
\Lambda(t)=\left\{%
\begin{array}{ll}%
0\,, & t<0\\
\delta^{-1}e^{t/\delta},&t>0
\end{array}%
\right.,\,  \text{if\ }\gamma<-1/\delta\,.
\end{equation}%
This function \(\Lambda(t)\) is totally positive. If  the formula
(\ref{InvLapl}) is used to construct the function \(\Lambda_j(t)\)
from  \(\Psi_j(z)\), \(j=1,\,2\), and the function \(\Lambda(t)\)
from the product \(\Psi(z)=\Psi_1(z)\Psi_2(z)\), and if the same
\(\gamma\) is used for all three constructions, then
\begin{equation}%
\label{Convol}%
\Lambda(t)=\int\limits_{-\infty}^{\infty}%
\Lambda_1(t-\xi)\Lambda_2(\xi)\,d\xi\,.
\end{equation}%
Moreover, if \(\Lambda_1\) and \(\Lambda_2\) are totally positive,
then the function \(\Lambda\) is totally positive as well.
Therefore, if \(\Psi(z)=\pm\prod\limits_{1\leq k\leq
n}(1+\delta_kz)\) is a polynomial with real roots, then   the
function \(\Lambda(t)\) defined by (\ref{InvLapl}) is totally
positive. Finally, if \(\Psi(z)\) is an entire function of the
second type, then there exists a sequence of polynomials
\(\Psi_n(z)\) with real roots such that \(\Psi_n(z)\to\Psi(z)\)
and, correspondingly, \(\Lambda_n(t)\to\Lambda(t)\). Therefore, if
\(\Psi(z)\) is an entire function of the second type, then the
corresponding function \(\Lambda(t)\) that is defined by formula
(\ref{InvLapl}) is totally positive. Thus,  part II of the theorem
is proved.

The statement that \textsl{the difference kernel
\(K(t,\,\tau)=\Lambda(y-\tau)\) is variation diminishing if and
\textsf{only if} the function \(\Lambda(t)\) is of the form
(\ref{InvLapl}), where \(\Psi(z)\) is an entire function of the
second type}, was formulated explicitly in \cite{Scho3}. In
particular, a difference kernel  \(K\) is variation diminishing if
and only if either the kernel \(K\) or the kernel  \(-K\)
is totally positive.\,%
\footnote{\,This agrees with the results of Schoenberg and Motzkin
on general variation diminishing transforms since a difference
kernels \(K(t,\,\tau)=\Lambda(t-\tau)\) is sign regular if and
only if either the kernel \(K(t,\,\tau)\) or
the kernel \(-K(t,\,\tau)\) is totally positive.} %
Many results related to totally positive and variation diminishing
difference kernels can be found in \cite{HiWi}, Chapter IV, and
especially in \cite{Karl}, Chapter 7.

Discrete totally positive difference kernels were first considered
in \cite{AESW}, \cite{ASW} and \cite{Edr}. The formulations are
analogous to the formulations for continuous difference kernels,
but the proofs are  more difficult and use tools from  value
distribution theory for meromorphic functions.

\noindent%
\textsf{THEOREM}(A.\,Aissen,\,A.\,Edrei,\,I.J.\,Schoenberg,\,A.\,Whitney,
\cite{AESW}, \cite{ASW}, \cite{Edr}).

\begin{slshape}
\textup{I}. Let \(\{s_k\}_{0\leq k<\infty}\) be a totally positive
(unilateral) sequence with \(s_0=1\). Then the series
\begin{equation}%
\label{TEUTPF}%
F(z)=\sum\limits_{0\leq k<\infty}s_kz^k
\end{equation}%
converges in a neighborhood of the origin to a function of the
form
\begin{equation}%
\label{GFTPS}%
F(z)=e^{\gamma z}\frac{\prod\limits_{
k}(1+\alpha_kz)}{\prod\limits_{k}(1-\beta_kz)}%
\quad (\alpha_k\geq 0,\,\beta_k\geq 0,\,\gamma\geq 0,\,
\sum (\alpha_k +\beta_k)<\infty).%
\end{equation}%

\textup{II}. Let \(F(z)\) be a function of the form
\textup{(\ref{GFTPS})} and let \textup{(\ref{TEUTPF})} be its
Taylor expansion in the vicinity of the origin. Then the sequence
\(\{s_k\}_{0\leq k<\infty}\) is totally positive.
\end{slshape}

This theorem provides a parametrization of the set of all totally
positive unilateral sequences (under the normalizing condition
\(s_0=1\)). The sequences \(\alpha_k,\,\beta_k\) and the number
\(\gamma\) serve as independent parameters. Various results
related to unilateral and bilateral totally positive sequences can
be found in \cite{Karl}, Chapter 8.

This parametrization of totally positive sequences plays an
essential role in the theory of representations of the infinite
symmetric group. It appears in the description of non-decomposable
positive definite  functions. This was discovered by Elmar Thoma
in \cite{Thom}, where some earlier results on totally positive
sequences were rediscovered. The generating function of the
sequences that  appear there are of the form (\ref{GFTPS}), with
\(\gamma=0\), and
\(\sum\limits_k\alpha_k+\sum\limits_k\beta_k\leq1\). The theory of
totally positive functions and sequences is used in approximation
theory, mathematical statistics and in other fields. References
can be found in \cite{Karl} and in \cite{GaMic}. Recently, a
surprising connection between  total positivity and canonical
bases for quantum groups was discovered by G.\,Lusztig;  see
\cite{FoZe}.

The methods and especially the ideology of the paper \cite{PS}
underly some of the work of B.Ya.\,Levin that is considered in
Chapter IX of his monograph \cite{Levi}. The famous
S.N.\,Bernstein inequality can be formulated in the following
form: \textsl{Let \(f(z)\) be an entire function of  exponential
type \(\sigma_f\). If the inequality \(|f(x)|\leq|e^{\sigma x}|\)
holds for all real \(x,\,\,-\infty<x<\infty\) and if
\(\sigma_f\leq\sigma\), then the derivative \(f^{\,\prime}(x)\)
satisfies the inequality
\(\big|f^{\,\prime}(x)\big|\leq\big|(e^{\sigma
x})^{\,\prime}\big|\,\ (-\infty<x<\infty).\)}\\
In other words, the operator \(\frac{d}{dx}\) preserves
inequalities  on the real axis for some classes  of entire
functions. Levin has investigated the general form of linear
operators which  preserve inequalities of this sort. In this
investigation, the linear operators that preserve the class of
entire functions that is obtained as the closure of polynomials
with zeros in the open right half plane play a crucial role. The
operators of the form (\ref{GammaOper}), which were introduced and
investigated in the paper \cite{PS}, are precisely those that
commute with the operator \(z\frac{d}{dz}\).


\vspace{1.0cm}


\setcounter{section}{8}
\begin{minipage}{15.0cm}
\section{\hspace{-0.5cm}.\hspace{0.2cm}%
\large The Schur class of holomorphic functions,\newline
 and the Schur algorithm.%
\label{SchAlg}}
\end{minipage}\\[0.18cm]
\setcounter{equation}{0}

The papers \cite{Sch9k}-\cite{Sch10k} are probably the best known contributions
of Issai Schur to analysis. In these papers Schur introduced a new
parametrization of functions that are holomorphic and bounded by one in the
open unit disk  \({\mathbb D}\)
and an algorithm for calculating these parameters.
These ideas and their subsequent generalizations to matrix and operator
valued functions are widely used in a variety of applications that range
from signal processing to the study of Pisot and Salem numbers.

\noindent%
\textsf{DEFINITION 1}. %
\begin{slshape}  A function \(s(z)\) that is holomorphic in the open unit
disk \({\mathbb D}\) and satisfies the inequality
\begin{equation}%
\label{ContrIneq}%
|s(z)|\leq 1 \qquad \mbox{for all points} \quad  z\in{\mathbb D}
\end{equation}%
is
said to belong to the \textsf{ Schur class} $\mathcal{S}$. A function
$s\in\mathcal{S}$ will be referred to as a  \textsf{ Schur function}.
A Schur function $s\in\mathcal{S}$ is said to be \textsf{inner} if
the absolute value of the radial limit
\begin{equation}
s(t)\stackrel{\mbox{\tiny def}}{=}\lim\limits_{r\rightarrow
1-0}s(rt) \label{DefBV}
\end{equation}
is equal to one a.e. with respect to Lebesgue measure. The set of inner
functions will be denoted by the symbol $\mathcal{S}_{in}$. The set of
rational inner functions will be denoted $\mathcal{S}_{rin}$.
\end{slshape}

The simplest inner functions are finite \textsf{Blaschke products}:
\begin{equation}%
\label{BRESF}%
s(z)=cz^{\kappa}\prod\limits_{k=1}^{d}%
\frac{\alpha_k-z}{1-\overline{\alpha_k}z\,},
\end{equation}%
where \(c\) is a constant of modulus one, \(\kappa\) is a non-negative
integer and  $\alpha_k\in\mathbb D$.

The  \textsl{\textsf{Schur algorithm}} was
 introduced by I.\,Schur in Section 1 of \cite{Sch9k}. It exploits
the fact that if \(\gamma\in\mathbb{D}\), then
the linear fractional transformation
\begin{equation}
\zeta\rightarrow \frac{\zeta - \gamma}{1-\zeta\,\overline{\gamma}}
\label{Transf}
\end{equation}
is a one to one mapping of the open unit
disk  \(\mathbb{D}\) onto itself and  a one to one  mapping of the unit circle
\(\mathbb{T}\)\,, i.e., the boundary of
 \(\mathbb{D}\)\,,  onto itself.
If \(|\gamma|=1\) the transformation (\ref{Transf}) maps the set
\(\mathbb{C}\setminus\{\gamma\}\) into the point \(\{-\gamma\}\)
and is \textit{not defined} at the point \(\gamma\).

Let \(f\in\mathcal{S}\) be a Schur function \textit{that is not a
constant of modulus one}. Then \(|f(0)|<1\) and hence, in view of the properties of
(\ref{Transf}),
the transformation
\begin{equation}
f(z)\rightarrow \frac{f(z)-f(0)}{1-f(z)\,\overline{f(0)}}
\label{FunctTransfnoz}
\end{equation}
maps $\mathcal{S}$ into $\{s\in\mathcal{S}:s(0)=0\}$. Therefore, by the
Schwarz Lemma, the transformation
\begin{equation}
f(z)\rightarrow \frac{f(z)-f(0)}{1-f(z)\,\overline{f(0)}}\cdot
\frac{1}{z} \label{FunctTransf}
\end{equation}
maps $\{s\in\mathcal{S}:|s(0)|\not=1\}$ onto the class
$\mathcal{S}$ of all Schur functions.
In particular if
\(f\in\mathcal{S}\setminus\mathcal{S}_{rin}\), then \(|f(0)|<1\) and
the transformation (\ref{FunctTransf}) is well defined. It is easy
to see that:

\noindent%
\textsf{PROPOSITION 1.} \begin{slshape}%
The
transformation (\ref{FunctTransf}) maps \(f\in\mathcal{S}\setminus\mathcal{S}_{rin}\) into itself.
\end{slshape}

\noindent%
\textsf{PROPOSITION 2.}
\begin{slshape}%
The transformation (\ref{FunctTransf} maps rational inner functions
$f\in\mathcal{S}_{rin}$
of degree
\(n,\,n\geq 1\),
into rational inner functions of degree \(n-1\).
\end{slshape}

\noindent \textsf{DESCRIPTION OF THE SCHUR ALGORITHM.} The Schur
algorithm defines a sequence of Schur
functions \(\{s_k(z)\}_{0\leq k<\infty}\) starting from a given Schur
function \(s(z)\) that is assigned the
index zero:
\begin{equation}
\label{SInitial}
s_0(z)\stackrel{\mbox{\rm \tiny def}}{=}s(z), \quad
\end{equation}
\begin{equation}
\label{SchurSeq}
 s_k(z)\stackrel{\mbox{\rm \tiny def}}{=}
\frac{s_{k-1}(z)-s_{k-1}(0)}{1-s_{k-1}(z)\,\overline{s_{k-1}(0)}}\cdot
\frac{1}{z} \quad(k=1,\,2,\,3\,,\dots)\,.
\end{equation}

\noindent \textsf{SCHUR PARAMETERS.} \textsl{
 Let \(s\in\mathcal{S}\)
and let \(\{s_k\}\) be the sequence (finite or infinite)
of functions generated by  the Schur
algorithm with  $s_0(z)=s(z)$. The numbers
\begin{equation}
\label{SchPar}
 {\gamma}_k\stackrel{\mbox{\rm\tiny
def}}{=}s_k(0)
\end{equation}
are termed
\textsf{the Schur parameters of the function} \(s\). }

If the starting function \(s\notin\mathcal{S}_{rin}\), then, by
Proposition 1, the
algorithm continues indefinitely
 and produces infinitely many Schur functions
\(s_k(z),\,k=0,\,1,\,2,\,3,\,\dots \)\, and generates
an infinite sequence of  Schur parameters
 \(\{\gamma_k\}_{0\leq k<\infty}.\)
  In this case
\begin{equation}
\label{SchParForNE} |\gamma_k|<1\,, \quad k=0,\,1,\,2,\,\dots\, .
\end{equation}
If the starting function \(s\in\mathcal{S}_{rin}\) is
is a rational inner function of degree \(n\), then, by Proposition 2,
the algorithm terminates after \(n\) steps.
In this case, it generates a finite sequence of
 Schur parameters
\begin{equation}
\label{SchParForE} |\gamma_k|<1\,, \quad
k=0,\,1,\,2,\,\dots\,,\,n-1 \quad \mbox{and} \quad |{\gamma}_n|=1.
\end{equation}\\[-4ex]

The Schur parameter \(\gamma_k(s)\) of a Schur function
\begin{equation}
\label{TayExp}
 s(z)=\sum\limits_{k=0}^{\infty}c_k(s)z^k
\end{equation}
 depends only on the Taylor coefficients
\(c_0(s)\),\,\(c_1(s),\)\,\(\dots,\,\) \(c_k(s)\) of the function
\(s\):
\begin{equation}
\label{TaylFromSch}
 {\gamma}_k(s)={\Phi}_k(c_0(s),\,c_1(s),\,\dots,\,c_k(s))\,,
\end{equation}
where \({\Phi}_k(c_0,\,c_1,\,\dots,\,c_k)\) is a rational function
of the variables
\(c_0,\,\overline{c_0},\,c_1,\,\overline{c_1},\,%
\dots,\,c_{k-1},\,\overline{c_{k-1}},\,c_{k}\,.\)

Conversely, the Taylor coefficient \(c_k(s)\) of  a Schur
function $s$ depends only on the Schur parameters
\(\gamma_0(s),\,\gamma_1(s),\,\dots,\,\gamma_k(s)\) of this
function:
\begin{equation}
\label{SchFromTayl}
 c_k(s)={\Psi}_k(\gamma_0(s),\,\gamma_1(s),\,\dots,\,\gamma_k(s))\,,
\end{equation}
where \({\Psi}_k(\gamma_0,\,\gamma_1,\,\dots,\,\gamma_k)\) is a
polynomial in
\(\gamma_0,\,\overline{\gamma_0},\,\gamma_1,\,\overline{\gamma_1},\,%
\dots,\,\gamma_{k-1},\,\overline{\gamma_{k-1}},\,\gamma_{k}\,.\)

Explicit expressions for \(\Phi_k\) and
\(\Psi_k\) are given in  \cite{Sch9k}.

\noindent%
\textsf{DEFINITION 4}.
\begin{slshape}
A sequence \(\{{\gamma}_k\}_{0\leq k<\infty}\) of complex numbers
is said to be \textsf{strictly contractive} if
\(|{\gamma}_k|<1\) for every \(k\).
\end{slshape}

Thus, the sequence of Schur parameters of a Schur function
$s\in\mathcal{S}\setminus\mathcal{S}_{rin}$ is strictly
contractive. Moreover,
 every preassigned strictly contractive sequence
\(\{\gamma_0,\,\gamma_1,\,\gamma_2,\,\dots,\,\gamma_k,\,\dots\}\)
is the sequence of  Schur parameters for some unique
 Schur function $s\in\mathcal{S}\setminus\mathcal{S}_{rin}$.
Such a function can be constructed
by means of a continued fraction algorithm.

\noindent \textsf{ SCHUR \ CONTINUED \ FRACTIONS.}\ \ \  Given an
arbitrary strictly contractive sequence  \(\{\gamma_0,\,\gamma_1,\,
\dots\,,\gamma_k,\,\dots\}\) of complex numbers, one can construct
a sequence of rational Schur functions which converge to a
Schur function \(s\) with Schur parameters
\(\{\gamma_0(s),\,\gamma_1(s),\, \dots\,,\\
\gamma_k(s),\,\dots\}\) that coincide with the preassigned sequence.
The construction is based on the inverse of the transformation
\begin{equation}
\label{ForvIterStep} f(z)\rightarrow
\frac{f(z)-\gamma}{1-f(z)\,\overline{\gamma}}\cdot
\frac{1}{z}\,,
\end{equation}
i.e., on the transformation
\begin{equation}
\label{BackStep} f(z)\rightarrow  \frac{ \gamma + zf(z) }{
1+\overline{\gamma} zf(z)}\,,
\end{equation}
which also maps $\mathcal{S}$ into $\mathcal{S}$.
We use the `\textsf{inverse Schur algorithm}'  recursively
to construct \textsl{the \(n\)-th Schur approximant}, which
 (following Schur) we will \textsl{denote} by
\([z;\,\gamma_0,\,\gamma_1,\,\dots\,,\gamma_n]\). Namely, we write
\begin{equation}
\label{RecurForSchAppr}
\begin{array}{l}
[z;\,\gamma_n]=\gamma_n; \\[2pt]
[z;\gamma_k,\,\gamma_{k+1},\,\gamma_{k+2},\,\dots\,,\gamma_n]=\displaystyle
\frac{ \gamma_k+z
\cdot[z;\,\gamma_{k+1},\,\gamma_{k+2},\,\dots\,,\gamma_n]} {1 +
\overline{\gamma_k} \cdot z \cdot [z;\,\gamma_{k+1},\,
\gamma_{k+2},\,\dots\,,\gamma_n]}\,,
\\[4pt]
k=n-1,\,n-2,\,\dots\,,1,\,0\,.
\end{array}
\end{equation}
The function \( [z;\,\gamma_0,\,\gamma_1,\,\dots\,,\gamma_n]\) is
a rational Schur function whose Schur parameters \\
\({\gamma}_k\big(\,[z;\,\gamma_0,\,\gamma_1,\,\dots\,,\gamma_n]\,\big)\)
are equal to
\[
\begin{array}{lll}
{\gamma}_k\big(\,[z;\,\gamma_0,\,\gamma_1,\,\dots\,,\gamma_n]\,\big)=
\gamma_k\quad &\mathrm{for}&\quad k=0,\,1,\,\dots\,,n;      \\[5pt]
{\gamma}_k\big(\,[z;\,\gamma_0,\,\gamma_1,\,\dots\,,\gamma_n]\,\big)=0
\quad &\mathrm{for}&\quad k>n\,.
\end{array}
\]

Let \(n_1\) and \(n_2\) be two nonnegative integers. Since the
Schur parameters with  index \(k: 0\leq k\leq
\mbox{min}(n_1,\,n_2)\) for the functions
\([z;\,\gamma_0,\,\gamma_1,\,\dots\,,\gamma_{n_1}]\) and
\([z;\,\gamma_0,\,\gamma_1,\,\dots\,,\gamma_{n_2}]\) coincide, the
Taylor coefficients \(c_k^1\) and \(c_k^2\) \((0\leq k\leq
\mbox{min}(n_1,\,n_2))\) for these two functions coincide as well.
Hence,
\[
[z;\,\gamma_0,\,\gamma_1,\,\dots\,,\gamma_{n_1}]-
[z;\,\gamma_0,\,\gamma_1,\,\dots\,,\gamma_{n_2}]=
 \sum\limits_{\mathrm{min}(n_1,n_2)< k<\infty}(c_k^1 -c_k^2)\,z^k\,.
\]
Using the estimates
\(\big|[z;\,\gamma_0,\,\gamma_1,\,\dots\,,\gamma_{n_1}]\big|\leq
1\),
\(\big|[z;\,\gamma_0,\,\gamma_1,\,\dots\,,\gamma_{n_2}]\big|\leq
1\) for \(z\in \mathbb{D}\) and the Schwarz Lemma, we obtain  the
inequality
\begin{equation}
\label{CauchySeq}
\big|[z;\,\gamma_0,\,\gamma_1,\,\dots\,,\gamma_{n_1}]-
[z;\,\gamma_0,\,\gamma_1,\,\dots\,,\gamma_{n_2}]\big|\leq
2\,|z|^{\,1+\mathrm{min}(\,n_1,\,n_2)} \quad \mathrm{for}\
z\in\mathbb{D}\,.
\end{equation}
From (\ref{CauchySeq}) it follows that  the limit
\begin{equation}
\label{ContFrac}
[z;\,\gamma_0,\,\gamma_1,\,\dots\,,\gamma_k,\,\dots\,]
\stackrel{\mbox{\rm\tiny def}}{=} \lim_{n\rightarrow\infty}
[z;\,\gamma_0,\,\gamma_1,\,\dots\,,\gamma_n]
\end{equation}
exists in $\mathbb{D}$. \textsl{The function
\([z;\,\gamma_0,\,\gamma_1,\,\dots\,,\gamma_k,\,\dots\,]\) is said
to be \textsf{the Schur continued fraction constructed from the
sequence \(\{\gamma_0,\,\gamma_1,\,\dots\,,\gamma_k,\,\dots\}\).}}

The function
\([z;\,\gamma_0,\,\gamma_1,\,\dots\,,\gamma_k,\,\dots\,]\in\mathcal{S}\setminus\mathcal{S}_{rin}\).  Its Schur parameters
\(\gamma_k\big([z;\,\gamma_0,\,\gamma_1,\,\dots\,,\gamma_k,\,\dots\,]\big)\)
coincide with the numbers \(\gamma_k\):
\[\gamma_k\big([z;\,\gamma_0,\,\gamma_1,\,\dots\,,\gamma_k,\,\dots\,]\big)=
\gamma_k,\ \ 0\leq k<\infty .\]

Given \(s\in\mathcal{S}\setminus\mathcal{S}_{rin}\), we can form the
sequence \(\{\gamma_0(s),\,\gamma_1(s),\,\)
\(\dots\,,\,\gamma_k(s),\,\dots\,\,\}\) of its Schur parameters
and then construct the Schur continued fraction
\([z;\,\gamma_0(s),\,\gamma_1(s),\,\dots\,,\,\gamma_k(s),\,\dots\,\,]\).
The function represented by this fraction is a Schur function
whose   Schur parameters coincide with the sequence of
Schur parameters of the original function \(s\). Hence,the  Taylor
coefficients of these two functions coincide as well. Thus, we are
led to following result:

\noindent%
\textsf{THEOREM} (I.\,Schur, \cite{Sch9k}).
\begin{slshape}

\textup{I.}\ \ Every  \(s\in\mathcal{S}\setminus\mathcal{S}_{rin}\) admits
\textsf{the continued fraction expansion}
\begin{equation}
\label{ContFractExpan}
s(z)=[z;\,\gamma_0(s),\,\gamma_1(s),\,\dots\,,\,\gamma_k(s),\,\dots\,\,].
\end{equation}

\textup{II.}\ \ A Schur function \(s\in\mathcal{S}_{rin}\) of degree
\(n\) admits the representation
\begin{equation}
\label{RESFSP}%
s(z)=[z;\,\gamma_0(s),\,\gamma_1(s),\,\dots\,,\,\gamma_n(s)]\,.
\end{equation}
\end{slshape}

\noindent%
\textsf{DEFINITION 5}.
\begin{slshape}
\label{DifSchAppr} \label{SchAppr} \ \
Let \(s\in\mathcal{S}\setminus\mathcal{S}_{rin}\),
let  \(n\) be a non-negative
integer and let (\ref{ContFractExpan}) be the Schur continued
fraction expansion of the function \(s\). Then the  function
\begin{equation}
\label{DefSchurAppr}%
\mathrm{p}_{n}(s;\,z) \stackrel{\mbox{\rm\tiny def}}{=}%
[z; \gamma_0(s),\,\gamma_1(s),\,\dots\,,\,\gamma_n(s)\,]
\end{equation}
is said to be \textsf{the \(n\)-th Schur approximant of the
function \(s\).}
\end{slshape}

\noindent \textsf{REMARK 1.}%
\begin{slshape}
The \(n\)-th Schur approximant is a rational function of \(z\)
whose numerator and denominator are polynomials of  degree not
greater than \(n\). In fact the \(n\)-th Schur approximant of a
 \(s\in\mathcal{S}\setminus\mathcal{S}_{rin}\) is the \(n\)-th convergent of
its Schur continued fraction expansion (\ref{ContFractExpan}).
\end{slshape}\\

\noindent The estimate
\begin{equation}
\label{EstAppr}%
\big|s(z)-\mathrm{p}_{n}(s;\,z) \big|\leq 2\,|z|^{n+1},
\end{equation}
which follows from the Schwarz Lemma,
holds for every  \(s\in\mathcal{S}\setminus\mathcal{S}_{rin}\)
and implies that
\textsl{the sequence of the Schur approximants of such an \(s\)
converges to it
locally uniformly in the open unit disk \(\mathbb{D}\).} This
result on the locally uniform convergence of the approximants
\(\mathrm{p}_{n}(s;\,z)\) to \(s(z)\) in the open unit disc \(\mathbb{D}\)
appears in \cite{Sch9k} (with  the rougher estimate
\(\big|s(z)-\mathrm{p}_{n}(s;\,z)\big|\leq 2|z|^{n+1}(1-|z|)^{-1}\)).
 The problem of convergence of Schur approximants to \(s\) on the
unit circle \(\mathbb{T}\) is much more difficult. This problem
was studied in \cite{Nja} and  \cite{Khru2}.

The preceding results imply
 that \textsl{ the correspondence
\begin{equation}
\label{FreeParam}
\{{\gamma}_0,\,{\gamma}_1,\,\dots\,,{\gamma}_k,\,\dots\}\longleftrightarrow
[z;\,\gamma_0,\,\gamma_1,\,\dots\,,\gamma_k,\,\dots\,]
\end{equation}
is a free parametrization of the class of all
$s\in\mathcal{S}\setminus\mathcal{S}_{rin}$
 by means of the set of all strictly contractive sequences
\(\{\gamma_0,\,\gamma_1,\,\gamma_2,\,\dots,\,\gamma_k,\,\dots\}\).
where sequences serve as free parameters of this class,}
This is important
because  the
geometry of the set of all Taylor coefficients of functions of
this class is rather complicated, whereas the geometry
of the set of all Schur parameters is very simple: it is just
the direct product of the open unit disks. This geometry is
compatible with probabilistic structures and is well suited for
probabilistic study. Some results on Schur functions with
random Schur parameters are obtained in \cite{Katsb}.

\noindent%
\textsf{REMARK 2.}
\begin{slshape}
It is not  easy to express the properties of a concrete Schur function
\(s\) in terms of its Schur parameters\,\footnote{\,However, to
express the properties of a Schur function in terms of its Taylor
coefficients is, as a rule, even
more difficult.} %
\({\gamma}_k(s)\). In particular, it is not  easy to recognize
whether the function \(s\) is inner or not. Not  much is known
about this.

If \(\sum\limits_{0\leq k<\infty}|{\gamma}_k(s)|<\infty\), then
the function \(s\) is continuous in the closed unit disk
\(\overline{\mathbb{D}}\), and
\(\max\limits_{z\in\mathbb{D}}|s(z)|<1\). (Of course, \(s\) is not
inner.) This result was obtained by I.\,Schur,
\textup{\cite{Sch10k}}, \S 15, Theorem XVIII.

If \(\sum\limits_{0\leq k<\infty}|{\gamma}_k(s)|^2<\infty\), then
again the function \(s\)  is not inner, as  follows from the
identity
\[
\prod\limits_{0\leq k<\infty}(1-|{\gamma}_k(s)|^2)=
\exp\Big\{\int\limits_{\mathbb{T}}\ln\big(1-|\,s(t)|^{\,2}\big)\,m(dt)\Big\}\,.
\]
\textup{(}See \textup{\cite{Boy}} and also   formula
\textup{(8.14)} in \textup{\cite{Gers3}}, which expresses a
similar result for polynomials that are orthogonal on
\(\mathbb{T}\)\textup{.)}

If \(\overline{\lim\limits_{k\to\infty}}|{\gamma}_k(s)|=1\), then
the function \(s\) is inner. An equivalent result was obtained by
E.A.\,Rakhmanov in the setting of orthogonal polynomials on the
unit circle, \textup{\cite{Rakh}}. It is known as Rakhmanov's
Lemma. A simple function-theoretic proof of the Rakhmanov lemma in
the setting of Schur functions can be found in
\textup{\cite{Katsb}}.

If the sequence of Schur parameters \(\{{\gamma}_k(s)\}_{0\leq
k<\infty}\) satisfies the \textit{Mat\'{e}--Nevai condition}
\(\lim\limits_{k\to\infty}{\gamma}_{k}{\gamma}_{k+n}=0\) for
\(n=1,\,2,\,3,\,\dots\, \), but
\(\overline{\lim\limits_{k\to\infty}}|{\gamma}_k|>0\), then \(s\)
is an inner function. This is \textup{Theorem 5} and
\textup{Corollary 9.1} in \textup{\cite{Khru2}}.

It is also known that there exists infinite Blaschke product \(s\)
such that \(\sum\limits_{0\leq
k<\infty}|{\gamma}_{k}(s)|^{p}<\infty\) for every \(p>2\).
\textup{(}This is shown in \textup{\cite{Khru3}.)}
\end{slshape}

The following question is both natural and important:\\[1.0ex]
\textsf{QUESTION 1}: \textsl{Given a sequence ${\bf c}=\{c_0, c_1, \ldots \}$,
 does there exist a function $s\in\mathcal{S}$ such that
$$
\frac{s^{(j)}(0)}{j!}=c_j  \quad \mbox{for}\quad j=0, 1, \ldots ?
$$}
Schur obtained an answer to this question by using the  algorithm
(\ref{SInitial})\,-\,(\ref{SchurSeq})
starting with
\begin{equation}%
\label{FPS}%
s_0(z) = f(z)
= \sum\limits_{j=0}^{\infty}c_jz^j\,,
\end{equation}%
to calculate the parameters $\gamma_j$. All the series are formal. However,
since the Schur parameters $\gamma_0, \ldots , \gamma_k$ only depend upon
$c_0, \ldots , c_k$, this does not present a problem. Proceeding this way,
Schur obtained the following answer to Question 1.

\textsl{In order for the series (\ref{FPS}) to be the Taylor series of
a Schur function, it is necessary and sufficient that
either
\(|\gamma_k(f)|<1\) for every integer \(k:\,\,0\leq
k<\infty\), or
\(|\gamma_k(f)|<1\) for \(k:\,\,0\leq k<n,\) and
\(|\gamma_n(f)|=1\).
 In the second case the coefficients
\(c_k\) of the series \textup{(\ref{FPS})} coincide with the
\(k\)-th Taylor coefficients of the function
\(\big[z;\,\gamma_0(f),\,\gamma_1(f),\,%
\dots\,\gamma_n(f)\big]\) for
every\,%
\footnote{\label{Fors}%
For \(k=0,\,1,\,\dots,\,n\) this coincidence
holds automatically since the Schur parameters
\(\gamma_0(f),\,\gamma_1(f),\,%
\dots\,,\,\gamma_n(f)\) are built from
\(c_0,\,c_1,\,\dots\,,\,c_n\); the remaining coefficients $c_k$  for \(k>n\) 
are determined by  \(c_0,\,c_1,\,\dots\,,\,c_n\).} %
\(k:\,0\leq k<\infty\).}

Schur also considered the following related question:

\noindent%
\textsf{QUESTION 2}. \textsl{Given a finite set of complex numbers
$\{c_0, \ldots ,c_m\}$, does there exist a  function $s\in\mathcal{S}$ such that
\begin{equation}
\label{SIC}
\frac{s^{(j)}(0)}{j!}=c_j\quad \mbox{for}\quad j=0, 1, \dots , m?
\end{equation}
Moreover, if such functions exist, how can one describe them?}

Schur answered  Question 2 in terms of the Schur parameters
generated by  the  algorithm
(\ref{SInitial})\,-\,(\ref{SchurSeq})
starting with
\begin{equation}%
\label{HNEW}%
s_0(z) = g(z)
= \sum\limits_{j=0}^{m}c_jz^j\,.
\end{equation}%

\noindent%
\textsf{THEOREM 1} (I.\,Schur,\,\cite{Sch9k}).\ \
 \textsl{There exists a function \(s\in\mathcal{S}\)
that meets the interpolation condition
\textup{(\ref{SIC})} if and only if:  either
\(|\gamma_k(g)|<1\) for every integer \(k:\,\,0\leq
k\leq m\), or
\(|\gamma_k(g)|<1\) for \(k:\,\,0\leq k<n,\) and
\(|\gamma_n(g)|=1\) for some \(n,\,n\leq m\).
In the second case, the interpolating
Schur function \(s(z)\) is unique, namely,
\(s(z)=\big[z;\,\gamma_0(g)),\,\gamma_1(g),\,%
\dots\,\gamma_n(g)\big]\). In the first case, there
are infinitely many interpolating functions \(s(z)\). Moreover,
the first $m+1$ Schur parameters \(\gamma_k(s)\),
\(k=0,\,1,\,\dots\,,\,m\), of every interpolant $s$ coincide with the
Schur parameters
\(\gamma_k(g)\) of the given polynomial. The remaining parameters
\(\gamma_k(s)\) with \(k>m\)  are either  an infinite sequence of
 arbitrary strictly contractive complex
numbers if $s\in\mathcal{S}\setminus\mathcal{S}_{rin}$, or a finite sequence
of strictly contractive numbers that terminates  with \(|\gamma_n(s)|=1\)
for some
\(m<n\). Moreover, all interpolants are obtained this way.}

Schur  also formulated another criterion for the solvability of the
interpolation problem (\ref{SIC}) in terms of an
\((m+1)\times
(m+1)\) upper triangular Toeplitz matrix based on the
the coefficients of the given polynomial:

\begin{equation}%
\label{TMID}%
\textsl{\textsf{C}}_m= %
\left[%
\begin{array}{ccccc}
c_0&c_1&c_2&
\ \ \ \cdots \ \ &c_m\\[1.2ex]
0&c_0&c_1&%
\cdots&c_{m-1}\\[1.2ex]
0&0&c_0&\cdots&c_{m-2}\\[1.2ex]
\cdot&\cdot&\cdot&\cdots&\cdot\\[1.2ex]
0&0&0&\cdots&c_0
\end{array}
\right]%
\end{equation}%

\noindent \textsf{THEOREM 2} (I.\,Schur,\,\cite{Sch9k}).
\begin{slshape}%
The given polynomial \(g(z)=c_0+c_1z+\cdots+c_mz^m\) can be
interpolated by a function \(s\in\mathcal{S}\) if and only the
Hermitian form based on the matrix
\(\textsl{\textsl{\textsf{I}}\,
-\,\textsf{C}}_m^{\,\ast}\textsl{\textsf{C}}_m\)
is non-negative:
\begin{equation}%
\label{NHF}%
\langle(\textsl{\textsl{\textsf{I}}\,-
\,\textsf{C}}_m^{\,\ast}\textsl{\textsf{C}}_m )\,x,\,x\rangle\geq
0,\quad \forall x\in\mathbb{C}^{m+1}\,,
\end{equation}%
where \(\textsl{\textsl{\textsf{I}}}\) is the identity matrix in 
\(\mathbb{C}^{m+1}\) and \(\langle\ \,,\ \rangle\) is the
standard scalar product in \(\mathbb{C}^{m+1}\). The Hermitian
form is strictly positive if and only if there exist more than one
interpolating function \(s\in\mathcal{S}\).
\end{slshape}%

We remark that the matrix \(\textsl{\textsf{C}}_m\) and an analogous 
matrix  \(\textsl{\textsf{D}}_m\) 
based on the coefficients of the reflected 
polynomial $z^m\overline{g(1/\overline{z})}=\overline{c_m}+\overline{c_{m-1}}z 
+ \cdots + \overline{c_0}z^m$ figure in the well known Schur-Cohn test:

\textsl{ The roots of the polynomial $g(z)$ lie in $\mathbb{D}$ if and 
only if 
\(\textsl{\textsf{D}}_m^*\textsl{\textsf{D}}_m 
- \textsl{\textsf{C}}_m^*\textsl{\textsf{C}}_m>0.\)} 

A nice proof of this result based on Schur parameters may be found in the 
first chapter of \cite{FoFr}.

Schur derived the criterion (\ref{NHF}) for the solvability of the
interpolation problem (\ref{SIC})
 from the
criterion  for solvability in terms of  Schur parameters (that was
formulated as Theorem 1).  As a by product of this 
derivation, he 
obtained a formula for the 
factorization of a \(2\times 2\) square block matrix 
\[M=\left[
\begin{array}{cc}
A&B \\
C& D
\end{array}
\right]\] with square block diagonal entries   \(A\)
and \(C\)  (not necessarily of the same size) when 
 the matrix \(A\) is invertible, which in turn leads easily to 
 the identity
\begin{equation}%
\label{SFF}%
\left[
\begin{array}{cc}
A&B \\
C& D
\end{array}
\right]= \left[
\begin{array}{cc}
I&0 \\
CA^{-1}& I
\end{array}
\right]\cdot \left[
\begin{array}{cc}
A&0 \\
0& D-CA^{-1}B
\end{array}
\right]\cdot \left[
\begin{array}{cc}
I&A^{-1}B \\
0& I
\end{array}
\right].
\end{equation}%
The  matrix \(D-CA^{-1}B\) is termed \textsl{the Schur complement of 
the block entry
\(A\) with respect to \(M\).} If it is also invertible, then 
$M$ is invertible and the last formula leads easily to a formula 
 for the inverse matrix
\(M^{-1}\), since each of the three factors on the right hand side of 
 (\ref{SFF}) are easily inverted. The same formula implies that 
\[\det M=\det A \cdot \det (D-CA^{-1}B),\] as was also noted  
in \cite{Sch9k}. 
  Schur
complements are  widely used in the applied linear algebra and operator 
theory; see e.g., the long survey (of more then hundred pages)
 \cite{Oue} and \cite{Sm}, respectively.

The \(k\)-th  step (\ref{SchurSeq})
of the Schur algorithm can also be presented in the form
\[
\left[\begin{array}{c}s_k(z)\\[5pt]1\end{array}\right]=
\left[\begin{array}{cc}
z^{-1}&-\gamma_{k-1}\,z^{-1}\\[5pt]
-\overline{\gamma_{k-1}}&1
\end{array}\right]\cdot
\left[\begin{array}{c}s_{k-1}(z)\\[5pt]1\end{array}\right]\cdot
\frac{1}{-\overline{\gamma_{k-1}}\,s_{k-1}(z)+1}\,,
\]
where \(\gamma_{k-1}=s_ {k-1}(0)\). Thus, it is natural to associate
the matrix \(\left[\begin{array}{cc}
z^{-1}&-\gamma_{k-1}\,z^{-1}\\[5pt]
-\overline{\gamma_{k-1}}&1
\end{array}\right]\)
with the \(k\)-th  step of  Schur algorithm. However, it
turns out to be more fruitful to deal with 
the matrix \(m_{\gamma_{k-1}}\), where 
\begin{equation}
\label{ElemStepSchurAlg} m_{\gamma}(z)= \left[
\begin{array}{cc}
\hspace{5pt} z^{-1}\hspace{5pt} &
\hspace{5pt}- {\gamma}\cdot {z}^{-1}\hspace{5pt} \\[10pt]
\hspace{5pt}-\overline{\gamma}\hspace{5pt} &\hspace{5pt}
1\hspace{5pt}
\end{array}
\right]\ \cdot \     \frac{1}{\sqrt{1-|\gamma|^{2}}}\,, \quad \mbox{when}
\quad |\gamma |<1.
\end{equation}
The matrix \(m_{\gamma}\) is also a matrix of coefficients of the 
linear fractional
transformation (\ref{ForvIterStep}), which is the basic step
(\ref{FunctTransf}) of the Schur algorithm, whereas the matrix 
\begin{equation}
\label{ESISA} {m_{\gamma}(z)}^{-1}= \left[
\begin{array}{cc}
\hspace{5pt} z\hspace{5pt} &
\hspace{5pt} \gamma\hspace{5pt} \\[10pt]
\hspace{5pt}z\cdot \overline{\gamma}\hspace{5pt} &\hspace{5pt}
1\hspace{5pt}
\end{array}
\right]\ \cdot \     \frac{1}{\sqrt{1-|\gamma|^{2}}}\,
\end{equation}%
is the coefficient matrix of the linear fractional transformation 
(\ref{BackStep}) 
corresponding to the basic step of the inverse Schur algorithm.

The coefficient matrix of a linear fractional transformation is determined up
to a nonzero scalar factor.  The matrix of the linear
fractional  transformation  (\ref{ForvIterStep}) is chosen to be of  the form
(\ref{ElemStepSchurAlg}) because then 
\(m_{\gamma}\) 
is \textsf{\textit{j}\,-\,inner} with respect to the signature matrix
\begin{equation}
\label{j} j= \left[
\begin{array}{cc}
\hspace{5pt}-1\hspace{5pt} &\hspace{5pt} 0\hspace{5pt}\\[10pt]
\hspace{5pt}\hspace{5pt} 0 &\hspace{5pt} 1\hspace{5pt}
\end{array}
\right], 
\end{equation}
i.e., 
\begin{equation}
\label{jForm}
({m_{\gamma}\,(z)}^{\ast})^{-1}\,j\,({m_{\gamma}\,(z)})^{-1}-j=
(1-|z|^{2})\cdot \left[
\begin{array}{c}
1\\
0
\end{array}
\right]\, \left[
\begin{array}{cc}
1&0
\end{array}
\right].
\end{equation}

\noindent%
\textsf{DEFINITION.}%
\begin{slshape}
\label{MainSeq}%
\ Let \ \(\omega_m=\{\gamma_k\}_{0\leq k\leq m}\) be a strictly
contractive sequence of complex numbers and let the entries of the matrix 
valued function 
\begin{equation}
\label{Mn}
\begin{array}{c}
M_{\,\omega_m}(z)\stackrel{\mbox{\tiny{def}}}{=}
m_{\,\gamma_{{}_{m}}}(z)\cdot
m_{\,\gamma_{{}_{m-1}}}(z)\cdot\,\,\dots\,\,\cdot
m_{\,\gamma_{{}_{1}}}(z)\cdot m_{\,\gamma_{{}_{0}}}(z)\,,
\quad m=0,\,1,\,2,\,\dots
\end{array}
\end{equation}
be denoted as
\begin{equation}%
\label{Entr}%
M_{\,\omega_m}(z)= \left[
\begin{array}{cc}
\hspace{3pt}a_{\,\omega_m}(z)\hspace{3pt}&
\hspace{3pt}b_{\,\omega_m}(z)\hspace{3pt}\\[5pt]
c_{\,\omega_m}(z)&d_{\,\omega_m}(z)
\end{array}
\right]\,\cdot
\end{equation}
\end{slshape}

For \(s(z)\in\mathcal{S}\setminus\mathcal{S}_{rin}\), let
\(\omega=\{\gamma_k\}_{0\leq k <\infty}\)  be the sequence of its
Schur parameters and let \(\{s_k(z)\}_{0\leq k <\infty}\) be the
sequence of Schur functions generated by the Schur algorithm (so that 
\(\gamma_k=s_k(0)\)). Then
\begin{equation}
\label{MatrSchur} M_{\,\omega_m}(z)\, \left[
\begin{array}{c}
s(z)\\[5pt]
1
\end{array}
\right] = \left[
\begin{array}{c}
s_{m+1}\,(z)\\[5pt]
1
\end{array}
\right]\cdot \,(c_{\,\omega_m}(z)\,s(z)+d_{\,\omega_m}(z))\,,
\end{equation}
and hence 
\begin{equation}
\label{SchurFrac} \frac{a_{\,\omega_m}(z)\,s(z)+b_{\,\omega_m}(z)}
{c_{\,\omega_m}(z)\,s(z)+d_{\,\omega_m}(z)}=s_{m+1}\,(z)\,.
\end{equation}%
Moreover,
\begin{equation}
\label{ResSchurFrac}%
 s(z)=%
\frac{w_{11}(z)\,s_{m+1}(z)+w_{12}(z)}
{w_{21}(z)\,s_{m+1}(z)+w_{22}(z)}\,,
\end{equation}
where the matrix \(W(z)=\left[
\begin{array}{cc}
\hspace{1pt}w_{11}(z)\hspace{1pt}&
\hspace{1pt}w_{12}(z)\hspace{1pt}\\[5pt]
w_{21}(z)&w_{22}(z)
\end{array}
\right]=M_{\,\omega_m}^{-1}(z)\) can be expressed in the term of
the entries of the matrix \(M_{\,\omega_m}(z)\):
\begin{equation}
\label{InvM}%
W(z)= \left[
\begin{array}{cc}
\hspace{3pt}d_{\,\omega_m}(z)\hspace{3pt}&
\hspace{3pt}-b_{\,\omega_m}(z)\hspace{3pt}\\[5pt]
-c_{\,\omega_m}(z)&a_{\,\omega_m}(z)
\end{array}
\right]\cdot \frac{1}{a_{\,\omega_m}(z)d_{\,\omega_m}(z)-
b_{\,\omega_m}(z)c_{\,\omega_m}(z)}\,\cdot
\end{equation}
Furthermore, since
\begin{equation}%
\label{EntrW}%
W(z)=m_{\,\gamma_{{}_{0}}}^{-1}(z)\cdot
m_{\,\gamma_{{}_{1}}}^{-1}(z)\cdot\,\,\dots\,\,\cdot
m_{\,\gamma_{{}_{m-1}}}^{-1}(z)\cdot
m_{\,\gamma_{{}_{m}}}^{-1}(z)\,,
\end{equation}
and the matrix \(m_{\gamma}^{-1}\) is  
linear with respect to \(z\),  the entries of the matrix
\(W(z)\) are polynomials with respect to  \(z\) of 
degree \(m\) (or less). In view of (\ref{jForm}), the matrix function
\(W(z)\) satisfies the condition
\begin{equation}%
\label{JExt}%
W^{\ast}(z)jW(z)-j\geq 0\quad \text{for}\quad z\in\mathbb{D}\,,
\end{equation}
\begin{equation}%
\label{JUnit}%
W^{\ast}(t)jW(t)-j= 0\quad \text{for}\quad t\in\mathbb{T}\,.
\end{equation}
Formula (\ref{ResSchurFrac}) (in other notation) appears in
\S\,14 of \cite{Sch10k}. It expresses the Schur function \(s(z)\) 
with Schur parameters \(\gamma_k(s)\) that satisfy the condition
\(|\gamma_k(s)|<1,\, k=0,\,1,\,\dots\,,\,n\), as a 
linear fractional transformation of the function \(s_{n+1}(z)\). It is 
important to note that \textsl{the matrix \(W(z)\) of this
fractional-linear transformation can be constructed from  only the
first \(n+1\) Schur parameters
\(\gamma_k(s),\,k=0,\,1,\,\dots\,,\,n\,\)}.

The Schur parameters of the functions \(s(z)\) and \(s_{m+1}(z)\)
are related:
\(\gamma_k(s_{m+1})=\gamma_{m+k}(s)\,,\ k=0,\,1,\,2,\,\dots\,. \)
Thus, the following result holds:

\noindent%
\textsf{THEOREM 3.} %
\begin{slshape}%
 The set of all functions
\(s(z)\in\mathcal{S}\), whose Schur parameters \(\gamma_k(s)\)
coincide with a  prescribed set of numbers \(\gamma_k,\ |\gamma_k|<1, \ %
k=0,\,1,\,\cdots,\,m\,,\) can be parametrized by means of the
 linear fractional transformation
\begin{equation}%
\label{ResMatr}%
 s(z)=%
\frac{w_{11}(z)\,\omega(z)+w_{12}(z)}
{w_{21}(z)\,\omega(z)+w_{22}(z)}\,,
\end{equation}%
where the coefficient matrix \(W(z)\) of this transformation can be
constructed from  only these \(\gamma_k\), and the free parameter 
\(\omega(z)\) is an arbitrary function from the class
\(\mathcal{S}\).
\end{slshape}

Schur did not mention either formula (\ref{ResSchurFrac}) 
or formula (\ref{ResMatr}) explicitly in his description  
of the solutions of the 
interpolation problem of the form
\begin{equation}%
\label{IPSP}%
\gamma_k(s)=\gamma_k,\quad k=0,\,1,\,\dots\,,\,m\,.
\end{equation}%
Nor do the properties 
(\ref{JExt}) and (\ref{JUnit}) of the matrix \(W(z)\) in
(\ref{ResMatr})  appear in Schur's work. But this
was the starting point of subsequent  research on
interpolation problems with constraints in various classes of
analytic functions, particularly that of M.\,Riesz and
R.\,Nevanlinna. The work \cite{Sch9k}-\cite{Sch10k} stimulated interest 
in  obtaining  matrices of  linear fractional 
transformations that appear in descriptions of the sets of
solutions of such problems. The methods described above are recursive and 
depend essentially upon formulas involving the Schur parameters.  
V.P.\,Potapov showed how to obtain an expression 
 for the matrix \(W(z)\) that appears in  the description
of the set of all solutions of the problem (\ref{IPSP}) in the
class \(\mathcal{S}\) directly in terms of the data
\(c_0,c_1,\,\dots\,,c_m\), without first calculating 
the Schur parameters. This method of V.P.\,Potapov, as applied to
the the interpolation problem (\ref{SIC}), is 
elaborated on in   great detail in the monograph \cite{DFK}.

Considerations related to  formula
(\ref{ResSchurFrac}) were used by Schur to obtain the following
result:

\begin{slshape}In
order that the function \(s(z)=\dfrac{p(z)}{q(z)}\in\mathbb{S}\),
where \(p(z)\) and \(q(z)\) are coprime polynomials, be
representable in the form
\(s(z)=[z;\,\gamma_0,\,\gamma_1,\,\dots\,,\,\gamma_m]\), with
\(m<\infty\), it is necessary and sufficient that the following
two conditions are satisfied:
\begin{enumerate}
\item %
The polynomial \(q(z)\) does not vanish in the closed unit
disc \(\mathbb{D}\);
\item %
The factorization identity \(|q(t)|^2-|p(t)|^2=r\), where \(r\) is
a positive number, holds for \(t\in\mathbb{T}\).
\end{enumerate}
\end{slshape}

The circle of problems related to  Schur functions and the Schur
algorithm is closely related to the theory of polynomials that are
orthogonal on the unit circle. Note that
\[\zeta\rightarrow\frac{1+\zeta}{1-\zeta}\]
is one-to-one mapping of the unit disk
 \(\{\zeta :\,|\,\zeta|<1\}\) onto  the
right half-plane \(\{\zeta:\,\text{Re}\,\zeta>0\}\). Therefore, if
\(s(z)\) is a Schur function, then the function
\begin{equation}
\label{SchurCarath} w(z)=\frac{1+zs(z)}{1-zs(z)}
\end{equation}
is a Carath\'{e}odory function, i.e., a function which is
holomorphic and has non-negative real part in the unit disk:
\begin{equation}
\label{CarathDef}
 \text{Re}\,w(z)\geq 0\quad \text{for}\ \ z\in\mathbb{D}\,.
\end{equation}
The factor \(z\)  in \((\ref{SchurCarath})\) leads to the
normalization
\begin{equation}
\label{NoCo} w(0)=1\,.
\end{equation}
Conversely, if \(w(z)\) is a Carath\'{e}odory function that satisfies
the normalization condition (\ref{NoCo}), then it can be  uniquely
represented in the form (\ref{SchurCarath}), where \(s(z)\) is a
Schur function. Every Carath\'{e}odory function \(w(z)\) which
satisfies the normalization condition (\ref{NoCo}) admits the
Herglotz representation
\begin{equation}
\label{HerRepr}
w(z)=\int\limits_{\mathbb{T}}\frac{t+z}{t-z}\,{\sigma}(dt)\,,
\end{equation}
where \(\sigma\) is a probability measure on \(\mathbb{T}\).
Conversely, if \(\sigma\) is a probability measure on
\(\mathbb{T}\), then formula (\ref{HerRepr}) defines a normalized
Carath\'{e}odory functiom \(w(z)\). Thus, the transformation
(\ref{SchurCarath}) together with the representation
(\ref{HerRepr}), establishes a one-to-one correspondence between
Schur functions and probability measures on \(\mathbb{T}\). It is
easy to see that
\[
S\in\mathcal{S}_{rin}\Longleftrightarrow\sigma \ \mbox{has finite support}
\Longleftrightarrow w(z) \ \mbox{is rational with all its poles on}
\ \mathbb{T}.
\]
Let \(\sigma\) be a probability measure  on
\(\mathbb{T}\) with infinitely many points of support and let
\(\{{\varphi}_k\}_{0\leq k<\infty}\) be a sequence
 of polynomials  that is orthonormal, with respect to  \(\sigma\).
Such a sequence
 can be obtained by applying the  Gram-Schmidt
orthogonalization procedure  to the sequence \(\{z^k\}_{0\leq
k<\infty}.\)
Let \({\varphi}^{\ast}_k(z)\stackrel{\text{\tiny
def}}{=}z^k\overline{{\varphi}_k(1/\overline{z})}\,\)
denote the reciprocal polynomial.
It turns
out that the system of polynomials
\(\{{\varphi}_k,\,{\varphi}_k^{\ast}\}_{0\leq k<\infty}\) satisfy
linear recurrence relations that  can be written in  the
form:
\begin{equation}
\label{RecRel} \left[
\begin{array}{c}
{\varphi}_{k+1}(z)\\[5pt]
{\varphi}^{\ast}_{k+1}(z)
\end{array}
\right]=%
\frac{1}{\sqrt{1-|\,a_k|^2}} \left[
\begin{array}{cc}  z \,\, \,& \,\,-\overline{a}_k\\[5pt]
  -z\,a_k \, \,         &\,\,1
\end{array}
\right] \cdot \left[
\begin{array}{c}
{\varphi}_{k}(z)\\[5pt]
{\varphi}^{\ast}_{k}(z)
\end{array}
\right],\ \  0\leq k <\infty\,,
\end{equation}
with the initial condition
\begin{equation}
\label{InitCond} \left[
\begin{array}{c}
{\varphi}_{0}(z)\\[5pt]
{\varphi}^{\ast}_{0}(z)
\end{array}
\right]= \left[
\begin{array}{c}
1\\[5pt]
1
\end{array}
\right]\,.
\end{equation}
Here \(\{a_k\}_{0\leq k<\infty}\) is a strictly contractive sequence of
complex numbers that is determined uniquely by
the probability measure \(\sigma\) that generates the
sequence \(\{{\varphi}_k\}_{0\leq k<\infty}\) of orthogonal
polynomials. The numbers \(a_k(\sigma)\) are termed the
\textsl{reflection coefficients}\, of the measure \(\sigma\) or of
the sequence \(\{{\varphi}_k\}\) of orthogonal polynomials.

\noindent %
\textsf{THEOREM}\,[Ya.L.\,Geronimus]. \textsl{Let
$s\in\mathcal{S}\setminus\mathcal{S}_{rin}$ and let the
normalized Carath\'{e}odory function
\(w(z)\) be related to $s$ by
\textup{(\ref{SchurCarath})}. Then the Schur parameters of $s$
coincide with the reflection coefficients of the measure sigma:
\[
\gamma_k(s) = a_k(\sigma) \quad \mbox{for} \quad 0\leq
k<\infty\,.
\]}
This theorem first appears  in \cite{Gers1}. It
also appears as Theorem 18.2 in \cite{Gers2}. Unfortunately, these papers
are not easily accessible. However, an English translation of the second
one is available.
A simplified
presentation of the cited theorem by Geronimus  can be found in \cite{Khru1}
and in \cite{PiNe}.

Many (but not all) properties of a Schur function \(s(z)\) can be
naturally reformulated in terms of the related  function \(w(z)\),
i.e., in terms of the related sequence of orthogonal polynomials.
In particular:  \textsl{a Schur function \(s(z)\) is inner if and
only if the related measure \(\sigma(dt)\) is singular.} Indeed,
\(|\,s(t)|=1\) if and only if \(\text{Re}\,w(t)=0\). On the other
hand, \(\text{Re}\,w(t)= {\sigma}^{\prime}(t)\) for \(m\) almost
every \(t\in\mathbb{T}\). (Here \(s(t)\) and \(w(t)\) are the
boundary values of the appropriate functions and
\({\sigma}^{\prime}(t)\) is the derivative of the measure
\(\sigma\) with respect to the normalized Lebesgue measure \(m\).)
Some other connections between Schur functions and orthogonal
polynomials can be found in \cite{Gol1}, \cite{Gol2} and
\cite{Khru2}.

There is a rich literature dedicated to the Schur algorithm and
related topics. The volume \cite{Ausg} containes a selection of
basic early papers on Schur analysis written in German (by
G.\,Herglotz, I.\,Schur, G.\,Pick, R.\,Nevanlinna, H.\,Weyl), as
well as the Afterword written by B.\,Fritzsche and
B.\,Kirstein, the editors of this volume. See also their survey
\cite{FrKi}. 
The last several years have witnessed an explosion of interest in 
generalizations 
of the Schur algorithm and the associated parametrization to matrix and 
operator valued functions. In particular,  a $p\times q$ matrix 
valued function $s(z)$ that is holomorphic in  \({\mathbb D}\)
and satisfies the inequality 
$$
|\xi^*s(z)\eta|\leq 1 \ \mbox{for every point} \ z\in\mathbb D \ 
\mbox{and every pair of unit vectors} \ \xi\in\mathbb C^p  \ 
\mbox{and} \ \eta\in\mathbb C^q
$$
is said to belong to the Schur class $\mathcal{S}^{p\times q}$. 
The Schur algorithm is developed in detail for this class of functions in
\cite{DeDy}; see also \cite{AlDy} for a reproducing kernel Hilbert space 
interpretation and additional generalizations to Pontryagin spaces and the 
references to both of these papers. There are also many papers by the team 
 P.\,Delsarte, Y.\,Genin and Y.\,Kamp that are devoted to generalizations 
of a number of the the themes discussed in this section; 
see e.g., \cite{DeGeK} for a start.
The books \cite{Alp}, \cite{Con}, \cite{DFK} and \cite{FoFr} are 
dedicated to  function theoretic questions related to the Schur
algorithm and its applications to  operator theory.  Applications of the 
Schur algorithm to  Pisot
and Salem numbers are considered in the book
\cite{BDGPS}.  
The Schur algorithm is also  useful in the setting of fast numerical
algorithms for systems of linear equations with structured
matrices\,%
\footnote{An algorithm for the fast inversion of Toeplitz matrices was
designed be N.\,Levinson, \cite{Levsn}, about fifty five years ago,
but it seems that the relationships with Schur algorithm did not
play any role there.} %
(Toeplitz, Hankel, etc.); see the book \cite{S:Meth} (and in particular the 
papers \cite{Kail} and
\cite{LevKai}).    The terminology "I.\,Schur
methods in signal processing" is now widely used in this connection. 
See the survey \cite{KaSa1} and the volume \cite{KaSa2} for
further references in this direction.


\setcounter{section}{9}
\begin{minipage}{15.0cm}
\section{\hspace{-0.5cm}.\hspace{0.2cm}%
\large Issai Schur's papers on analysis.
\label{Schpapers}}
\end{minipage}\\[0.18cm]
\setcounter{equation}{0}



\begin{minipage}{45.0ex}
Harry Dym\\
Department of Mathematics\\
The Weizmann Institute of Science\\
Rehovot, 76100, Israel\\
E-mail: dym@wisdom.weizmann.ac.il
\end{minipage}
\begin{minipage}{45.0ex}
Victor Katsnelson\\
Department of Mathematics\\
The Weizmann Institute of Science\\
Rehovot, 76100, Israel\\
E-mail: katze@wisdom.weizmann.ac.il
\end{minipage}

\begin{thebibliography}{MoMo2}
\small
\bibitem[Adl]{Adl} \textsc{Adler,\,M.}: \textsl{On a trace
functional for formal pseudo-differential operators and the
symplectic structure of the Korteveg-Devries equations.}
Inventiones Math., \textbf{50} (1979), pp.\,219 - 248.


\bibitem[Ami]{Ami} \textsc{Amitsur,\,S.A.} Commutative linear
differential operators. Pacif. Journ. of Math. \textbf{8} (1958),
pp.\,1 - 10.

\bibitem[Bak1]{Bak1} \textsc{Baker, H.F.} Note on the foregoing
paper "Commutative ordinary differential operators". Proc. Royal
Soc. London, \textbf{118} (1928), pp.\,584 - 593.

\bibitem[Bak2]{Bak2} \textsc{Baker, H.F.} \textsl{Abelian
Functions.} (Cambridge Mathematical Library). Cambridge Univ.
Press, Cambridge 1995.

\bibitem[BuCh1]{BuCh1} \textsc{Burchnall,\,J.L.} and
\textsc{T.W.\,Chaundy}. \textsl{Commutative ordinary differential
operators.} Proc. London Math. Soc. (Second Ser.), \textbf{21}
(1922), pp.\,420 - 440.

\bibitem[BuCh2]{BuCh2} \textsc{Burchnall,\,J.L.} and
\textsc{T.W.\,Chaundy}. \textsl{Commutative Ordinary Differential
Operators.} Proc. Royal Soc. London, \textbf{118} (1928), pp.\,557
- 583.

\bibitem[BuCh3]{BuCh3} \textsc{Burchnall,\,J.L.} and
\textsc{T.W.\,Chaundy}. \textsl{Commutative Ordinary Differential
Operators. II.--- The Identity \(P^n=Q^m\).} Proc. Royal Soc.
London, \textbf{134} (1932), pp.\,471 - 485.

\bibitem[GGKM]{GGKM} \textsc{Gardner,\,C., J.\,Green,\,M.\,Kruskal, R.\.Miura.}
\textsl{A method for solving the Korteweg --- de Vries equation.}
Phys. Rev. Lett., \textbf{19} (1967), pp.\,1095 - 1098.

\bibitem[GeDi]{GeDi} \textsc{Gel'fand,\,I.M.} and \textsc{L.A.\,Dikii}
(=\textsc{L.A.\,Dickey}). \textsl{Drobnye stepeni operatorov i
Hamil'tonovy sistemy.} Funkc. Anal. i ego Prilozh., \textbf{10}:4,
pp.\,13 - 29. {(\footnotesize Russian). English transl.:}
\textsl{Fractional powers of operators and Hamiltonian systems.}
Funct. Anal. and Appl. \textbf{10}:4 (1976), pp.\,259-273.

\bibitem[Kri1]{Kri1} \textsc{Krichever,\,I.M.} \textsl{Integrirovanie
neline{\u\i}nykh uravneni{\u\i} metodami algebraichesko{\u\i} geometrii.} Funk.
Analiz i ego prilozh., \textbf{11}:1  (1977), pp.\,15 - 31 ({\footnotesize
Russian}). {\footnotesize English transl.}: \textsl{Integration of nonlinear
equations by the methods of algebraic geometry. } Funct. Anal. and Appl.
\textbf{11}:1 (1977), pp.\,12-26.

\bibitem[Kri2]{Kri2} \textsc{Krichever, I.M.} \textsl{Metody algebraichesko\u{\i}
geometrii v teorii neline\u{\i}nykh uravneni\u{\i}} (Russian).
Uspeki Matem. Nauk, \textbf{32}:6 (1977), pp.\,183 - 208. English
transl.: \textsl{Methods of algebraic geometry in the theory of
non-linear equations}. Russ. Math. Surveys, \textbf{32}:6 (1977),
pp.\,185-213.

\bibitem[Kri3]{Kri3} \textsc{Krichever, I.M.} \textsl{Kommutativnye kol'tsa
obyknovennykh differentsial'nykh operatorov.} Funkts. Anal. i ego
Prilozh., \textbf{12}:3 (1978), pp.\,21 - 30. {\footnotesize
(Russian). English transl.:} \textsl{Commutative rings of
ordinary linear differential operators}.  Funct. Anal. and Appl.
\textbf{12}:3 (1977), pp.\,175 -- 185.

\bibitem[Kri4]{Kri4} \textsc{Krichever,\,I.M.} Foreword to the
monograph \cite{Bak2}, Reedition of 1995, pp. xvii - xxxi in
\cite{Bak2}.

\bibitem[Lax]{Lax} \textsc{Lax,\,P.} \textsl{Integrals of nonlinear equations
of evolution and solitary waves}. Comm. Pure Appl. Math.
\textbf{21} (1968), pp.\,467--490.

\bibitem[LebMa]{LebMa} \textsc{Lebedev\,B.M.} and
\textsc{Yu.I.\,Manin.} \textsl{Gamil'tonov operator
Gel'fanda-Dikogo i kopri\-soedinennoe predstavlenie gruppy
Vol'terra.} Funkts. Anal. i ego Prilozh. \textbf{13}:4,
pp.\,40-46 (1979) ({\footnotesize Russian}). {\footnotesize
English transl.:} \textsl{Gel'fand-Dikii Hamiltonian operator and
the coad\-joint  representation of the Volterra group},
\textbf{13}:4 (1979), pp.\,268 - 273.

\bibitem[Man]{Man} \textsc{Manin,\,Yu.I.} \textsl{Algebraicheskie
aspecty neline{\u\i}nykh differentsial'nykh uravneni{\u\i}.}
Itogi Nauki i Tekhniki (Sovremennye Problemy Matematiki,
\textbf{11}), VINITI, Moscow  1978, pp.\,5 - 152 ({\footnotesize
Russian}). {\footnotesize English transl.}:
 \textsl{Algebraic aspects of nonlinear differential
equations.}  Journ. of Soviet. Math. \textbf{11}, (1979), pp.\,1
- 122.

\bibitem[Mum]{Mum} \textsc{Mumford,\,D.} \textsl{Tata Lectures on
Theta. Vol.\,II.} (Progress in Mathematics, Vol.\,\textbf{43})
Birkh\"auser, Boston\(\cdot\)Basel\(\cdot\)Stuttgart 1984.

\bibitem[Sch:Ges]{Sch:Gesa} \textsc{Schur,\,I.}:
\textsl{Gesammelte Abhandlungen
[Collected Works]. Vol. I, II, III.} Springer-Verlag, Berlin\(\cdot\)%
Heidelberg\(\cdot\)New\,York 1973.

\bibitem[Sch1]{Sch1a} \textsc{Schur,\,I.} \textsl{\"Uber
vertauschbare lineare Differentialausdr\"ucke} [On permutable
differential expressions - in German]. Sitzungsberichte der
Berliner Mathematischen Gesellschaft, \textbf{4} (1905), pp.\, 2 -
8. \footnotesize Reprinted in: \small \cite{Sch:Gesa}, Vol.\,I,
pp.\,170 - 176.

\bibitem[Tsuj]{Tsuj} \textsc{Tsujishita, \,Toro}: \textsl{Formal
geometry of systems of differential equations}. Sugaku
Expositions, \textbf{3}:1 (1990), pp.\,25 - 73.

\end{thebibliography}

\vspace{3.0ex} \noindent
\begin{thebibliography}{MoMo2}
\small

\bibitem[Bo]{Bo} \textsc{Boos,\,J.} \textsl{Classical and modern methods in 
summability.} (Assisted by Peter Cass.) Oxford Mathematical Monographs. 
Oxford Science Publications. Oxford University Press, Oxford 2000. 

\bibitem[Coo]{Coo} \textsc{Cooke,\,R.G.} \textsl{Infinite Matrices
and Sequence Spaces.} MacMillan, London 1950.

\bibitem[Har]{Har} \textsc{Hardy,\,G.H.} \textsl{Divergent
Series}. Clarendon Press, Oxford 1949.

\bibitem[Koj]{Koj} \textsc{Kojima,\,T.} \textsl{On generalized
Toeplitz's theorems and their applications.} T\^{o}hoku
Mathematical Journal, \textbf{12} (1917), pp.\,291 - 326.

\bibitem[Pet]{Pet} \textsc{Petersen,\,G.M.} \textsl{Regular Matrix
Transformation}. McGraw-Hill Publishing Compamy,
London\(\cdot\)New\,York\(\cdot\)\,Toronto\(\cdot\)\,Sidney 1966.

\bibitem[Pey]{Pey} \textsc{Peyerimhoff,\,A.}
         \textsl{Lectures on summability.}
          (Lect. Notes in Math, \textbf{107}),
         Spinger Verlag,
         Berlin\(\cdot\)Heidelbrg\(\cdot\)New\,York
         1969.

\bibitem[Sch:\,Ges]{Sch:Gesb} \textsc{Schur,\,I.}:
\textsl{Gesammelte Abhandlungen
[Collected Works]. Vol. I, II, III.} Springer-Verlag, Berlin\(\cdot\)%
Heidelberg\(\cdot\)New\,York 1973.


\bibitem[Sch6]{Sch6a} \textsc{Schur,\,I.}:
\textsl{\"Uber die \"Aquivalenz der Ces\`aroschen und
H\"olderschen Mittelwerte.} [\textsl{On the equivalence of
Ces\`aro's and H\"older means} - in German]. Mathematische
Annalen, \textbf{74} (1913), pp.\,447 - 458. \footnotesize
Reprinted in: \small \cite{Sch:Gesb}, Vol.\,II, pp.\,44 - 55.

\bibitem[Sch16]{Sch16a} \textsc{Schur,\,I.}:
 \textsl{\"Uber lineare Transformationen
in der Theorie der unendlichen Reihen.} [\textsl{On linear
transformations in the theory of infinite series} - in German].
Journ. f\"ur die reine und angew. Math., \textbf{151} (1921),
pp.\,79 - 121. \footnotesize Reprinted in: \small \cite{Sch:Gesb},
Vol.\,II, pp.\,289 - 321.

\bibitem[Sch19]{Sch19a} \textsc{Schur,\,I.}
\textsl{Einige Bemerkungen zur Theorie der unendlichen Reihen.}
 [\textsl{A remark on the theory of infinite series} - in German].
  Sitzungsberichte der Berliner mathematisches Gesellschaft,
\textbf{29}, (1930), pp.\,3 - 13. \footnotesize Reprinted in:
\small \cite{Sch:Gesb}, Vol.\,III, pp.\,216 - 226.

\bibitem[Silv]{Silv} \textsc{Silverman,\,L.L.} \textsl{On the definition
of the sum of a divergent series.} University of Missouri Studies,
Math. Series I, (1913), pp.\,1 - 96.

\bibitem[Toep]{Toep} \textsc{Toeplitz,\,O.}:
\textsl{\"Uber allgemeine lineare Mittelbildungen.} Prace
matematyczno-fizyczne (Warszawa), \textbf{22} (1913), pp.\,113 -
119.
\bibitem[Zel]{Zel} \textsc{Zeller,\,K.} \textsl{Theorie der
Limitierungsverfahren.} Springer, Berlin 1958.

\end{thebibliography}

\vspace{3.0ex} \noindent
\begin{thebibliography}{MoMo2}
\small

\bibitem[BiSo]{BiSo} \textsc{Birman,\,M.S.} and \textsc{M.Z.\,Solomyak.}
\textsl{Spektral'naya Teoriya Samosopryazhennykh Operatorov v
Gil'bertovom Prostranstve} (Russian), Izdatel'stvo Leningradskogo
Universiteta, Leningrad 1980, 264 pp. English transl.:
\textsl{Spectral Theory of Self-Adjoint Operators in Hilbert
Space.} D.\,Reidel Publishing Company,
Dordrecht\(\boldsymbol{\cdot}\)Boston\(\boldsymbol{\cdot}\)
Lancaster\(\boldsymbol{\cdot}\)Tokyo 1987, xvi, 301 pp.

\bibitem[HLP]{HLPa} \textsc{Hardy,\,G.H., J.W.\,Littlewood,} and
\textsc{G.\,Polya.} \textsl{Inequalities}. 1st ed., 2nd ed.
Cambridge Univ. Press, London\(\cdot\)New\,York 1934, 1952.

\bibitem[Sch:\,Ges]{Sch:Gesc} \textsc{Schur,\,I.}:
\textsl{Gesammelte Abhandlungen \textup{[}Collected
Works\textup{]}. Vol. I, II, III.} Springer-Verlag,
Berlin\(\cdot\)Heidelberg\(\cdot\)New\,York 1973.

\bibitem[Sch4]{Sch4a} \textsc{Schur,\,I.}: \textsl{Bemerkungen zur
Theorie der beschr\"ankten Bilinearformen mit unendlich vielen
Ver\"anderlichen.}  [\textsl{Remarks on the theory of bounded
bilinear forms with infinitely many variables} - in German].
Journ. f\"ur reine und angew. Math., \textbf{140} (1911), pp.\,1 -
28. \footnotesize Reprinted in: \small \cite{Sch:Gesc}, Vol.\,I,
pp.\,464 - 491.

\bibitem[Sch18]{Sch18b} \textsc{Schur,\,I.}: \textsl{\"Uber eine Klasse
von Mittelbildungen mit Anwendungen auf die Determinantentheorie.}
[\textsl{On a class of averaging mappings with applications to
the theory of determinants} - in German].  Sitzungsberichte der
Berliner Mathematischen Gesellschaft, \textbf{22} (1923), pp.\,9 -
20. \footnotesize Reprinted in: \small \cite{Sch:Gesc}, Vol.\,II,
pp.\,416 - 427.

\bibitem[Tho]{Tho} \textsc{Thorin,\,G.O.} \textsl{ Convexity theorem
generalizing those of M. Riesz and Hadamard with some
applications.} Comm. Sem. Math. Univ. Lund [Medd. Lunds Univ. Mat.
Sem.] \textbf{9}, (1948), pp.\,1--58.
\end{thebibliography}

\vspace{3.0ex} \noindent
\begin{thebibliography}{MoMo2}
\small

\bibitem[ABF]{ABF} \textsc{Arazy,\,J., T.J.\,Barton} and \textsc{Y.\,Friedman.}
\textsl{Operator differentiable functions.} Integral Equat. and
Oper. Theory, \textbf{13} (1990), pp.\,461--487.

\bibitem[Ben]{Ben} \textsc{Bennett, G.} \textsl{Schur multipliers.}
Duke Math. Journ., \textbf{44}:3 (1977), pp.\,603 - 639.

\bibitem[Bie]{Bie} \textsc{Bieberbach,\,L.} \textsl{Analytische Fortsetzung}
[\textsl{Analytic Continuation}] (in German). (Ergebnisse der
Mathematik und ihrer Gerenzgebiete, \textbf{3}). Springer Verlag,
Berlin\(\boldsymbol{\cdot}\)G\"{o}ttingen\(\boldsymbol{\cdot}\)Heidelberg
1955.

\bibitem[Bi1]{Bi1} \textsc{Birman, M.Sh.} \textsl{Ob usloviyakh sushchestvovaniya
volnovykh operatorov} (Russian). Izvestiya AN SSSR (Ser. Mat.),
\textbf{27}:4 (1963), pp.\,883 - 906. English transl.: \textsl{On
conditions for the existence of wave operators},
Amer.\,Math.\,Soc.\,Transl. (Ser. \textbf{2}), \textbf{54} (1966),
pp.\,91 -- 117.

\bibitem[Bi2]{Bi2} \textsc{Birman,\,M.Sh.} \textsl{Lokal'ny\u{\i} priznak
 sushchestvovaniya volnovykh operatorov} (Russian). Dokl. Akad. Nauk SSSR,
 \textbf{159}:3 (1964),
 pp.\,485 - 488. English transl.: \textsl{A local criterion for the
existence of wave operators }. Soviet Math. Dokl. \textbf{5}
(1965), pp.\,1505 - 1509.

\bibitem[BiSo1]{BiSo1}\textsc{Birman,\,M.Sh.} and \textsc{M.Z.\,Solomyak.}
\textsl{O dvo\u{\i}nykh operatornykh integralakh Stil'tyesa}
(Russian). Dokl. AN SSSR, \textbf{165}:6 (1965), pp. 1223 - 1226.
English transl.: {Stieltjes double operator integrals.} Soviet
Math., Dokl. \textbf{6}, (1965), 1567-1571.

\bibitem[BiSo2]{BiSo2}\textsc{Birman,\,M.Sh.} and \textsc{M.Z.\,Solomyak.}
\textsl{Dvo\u{\i}nye operatornye integraly Stil'tyesa} (in
Russian). In "Problemy matematichesko\u{\i} fiziki",
No.\,\textbf{1}. Spektral'naya Teoriya i Volnovye Processy
[Spectral Theory and Wave Processes].
(M.Sh.\,Birman-editor).Izdat. Lenigradskogo Univ., Leningrad 1966,
pp.33 - 67. English transl.: \textsl{Stieltjes double-integral
operators}, in "Topics in Mathematical Physics", vol.\textbf{1},
Consultants Bureau, New York 1967, viii+114 pp.
\bibitem[BiSo3]{BiSo3}\textsc{Birman,\,M.Sh.} and \textsc{M.Z.\,Solomyak.}
\textsl{Dvo\u{\i}nye operatornye integraly Stil'tyesa. II} (in
Russian). In "Problemy matematichesko\u{\i} fiziki",
No.\,\textbf{2}. Spectral'naya teoriya. Zadachi diffraktsii
[Spectral theory. Diffraction problems]. (Edited by
M.Sh.\,Birman). Izdat. Lenigradskogo Univ., Leningrad 1967, pp.26
- 60. English transl.: \textsl{Stieltjes double-integral
operators.II}, in "Topics in Mathematical Physics",
vol.\textbf{2}, Consultants Bureau, New York 1968, vii+134 pp.

\bibitem[BiSo4]{BiSo4}\textsc{Birman,\,M.Sh.} and \textsc{M.Z.\,Solomyak.}
\textsl{Dvo\u{\i}nye operatornye integraly Stil'tyesa. III}
(Russian). In "Problemy matematichesko\u{\i} fiziki", \textbf{6},
Izdat. Lenigradskogo Univ., 1973, pp.27 - 53.

\bibitem[BiSo5]{BiSo5}\textsc{Birman,\,M.Sh.} and \textsc{M.Z.\,Solomyak.}
\textsl{O priblizhenii funktsi\u{\i} klassov \(W_p^{\alpha}\)
kusochno-polinomial'nymi funktsiyami} (Russian), Dokl. Akad. Nauk
SSSR, \textbf{171}:5 (1966), pp.\,1015 -1018. English transl.:
\textsl{Approximation of the classes \(W_p^{\alpha}\) by piceweise
polynomial functions.} Sov. Math., Dokl. \textbf{7} (1966),
pp.\,1573 - 1577.

\bibitem[BiSo6]{BiSo6} \textsc{Birman,\,M.Sh.} and \textsc{M.Z.\,Solomyak.}
\textsl{Kusochno-polinomial'nye priblizheniya funktsi\u{\i}
klassow \(W_p^{\alpha}\).} (Russian) Matem. Sbornik (N.S.)
\textbf{73} (\textbf{115}) (1967), pp.\,331--355. English transl.:
\textsl{Piecewise-polynomial approximations of functions of the
classes \(W_p^{\alpha}\)}. Math.\,USSR--Sbornik, \textbf{2}:3
(1967), pp.\,295 - 317.

\bibitem[BiSo7]{BiSo7} \textsc{Birman,\,M.Sh.} and \textsc{M.Z.\,Solomyak.}
\textsl{Otsenki singulyarnykh chisel integral'nykh operatorov}
(Russian). Uspekhi Matem. Nauk \textbf{32}:1(193), (1977),
pp.\,17-84, 271. English transl.:\textsl{Estimates of singular
numbers of integral operators}, Russ. Math. Surveys \textbf{32}:1
(1977), pp.\,15-89.

\bibitem[BiSo8]{BiSo8} \textsc{Birman,\,M.Sh.} and \textsc{M.Z.\,Solomyak}
\textsl{Kachestvenny\u{\i} analiz v teoremakh vlozheniya Soboleva
i prilozheniya k spektral'no\u{\i} teorii} (Russian). Desyataya
matematicheskays shkola. (Letnyaya shkola, Katsiveli/Nal'chik,
1972), pp. 5--189. Izdanie Inst. Matem. Akad. Nauk Ukrain. SSR,
Kiev, 1974. English transl.: \textsl{Quantitative analysis in
Sobolev imbedding theorems and applications to spectral theory.}
American Mathematical Society Translations, Ser.\,\textbf{2},
\textbf{114}. American Mathematical Society, Providence, R.I.,
1980. viii+132 pp.

\bibitem[BiSo9]{BiSo9}\textsc{Birman,\,M.Sh.} and \textsc{M.Z.\,Solomyak.}
\textsl{Dvo\u{\i}nye operatornye integraly Stil'tyesa i zadacha o
mnozhitekyakh.} Dokl. Akad. Nauk SSSR, \textbf{171}:6 (1966), 1251
- 1254. English transl.: \textsl{Double Stieltjes operator
integrals and problem of multipliers.} Soviet Math. Dokl.,
\textbf{7} (1966), 1618 - 1621.

\bibitem[BiSo10]{BiSo10} \textsc{Birman,\,M.Sh.} and \textsc{M.Z.\,Solomyak.}
\textsl{Zamechaniya o funktsii spektral'nogo sdviga} (Russian).
Zapiski Nauchnykh Seminarov LOMI, \textbf{27} (1972), pp.\,33 --
46. English transl.: \textsl{Remarks on the spectral shift
function}. Journ. Soviet Math. \textbf{3}:4 (1975), pp.\,408 -
419.

\bibitem[BiSo11]{BiSo11} \textsc{Birman,\,M.Sh.} and \textsc{M.Z.\,Solomyak.}
\textsl{Operatornoe integrirovanie, vozmushcheniya i kommutatory}
(Russian). Zapiski Nauchnykh Seminarov LOMI, vol.\,\textbf{170}
(1989). English transl.: \textsl{Operator integration,
perturbations, and commutators}. Journal of Soviet. Math.,
\textbf{63}:2 (1993), pp.\,129 -- 148.
\bibitem[BiYa]{BiYa} \textsc{Birman,\,M.Sh.} and \textsc{D.R.\,Yafaev.}
\textsl{Funktsiya spektral'nogo sdviga. Raboty M.G.\,Kre\u{\i}na i
ikh dal'ne\u{\i}shee razvitie} (Russian). Algebra i Analiz,
\textbf{4}:5 (1992), pp.\,1 -- 44. English transl.: \textsl{The
spectral shift functioin. The work of M.G.Kre\u{\i}n and its
further development}. St. Petersburg Math. J. \textbf{4}:5 (1993),
pp.\,833 - 870.

\bibitem[Da1]{Da1} \textsc{Daletskii, \,Yu.L.} \textsl{Pro otsinku
zalushkovogo chlenu u formuli Te\u{\i}lora dlya funktsi\u{\i}
ermitovykh operatoriv} [On an estimate of the remainder term in
Taylor's formula for functions of Hermitian operators]
(Ukrainian). Dopovidi Akad. Nauk Ukra\"{\i}n. RSR, \textbf{4}
(1951), pp.\, 234--238.

\bibitem[Da2]{Da2} \textsc{Daletskii, \,Yu.L.}
\textsl{Integrirovanie i differentsirovanie funktsi\u{\i}
ermitovykh operatorov, zavisyashchikh ot parametra} (in Russian).
Uspehi Mat. Nauk (N.S.) \textbf{12} (1957), no. 1(73), 182--186.
English transl. \textsl{Integration and differentiation of
functions of Hermitian operators depending on a parameter}, Amer.
Math. Soc. Transl. (Ser. \textbf{2}), \textbf{16} (1960),
pp.\,396--400.

\bibitem[DaKr1]{DaKr1} \textsc{Daletskii, \,Yu.L.} and \textsc{Krein,\,S.G.}
\textsl{Formuly differentsirovaniya po parametru funktsi\u{\i}
ermitovykh operatorov} (Russian). [Formulas of differentiation
according to a parameter of functions of Hermitian operators].
Doklady Akad. Nauk SSSR \textbf{76} (1951), pp.\,13--16.

\bibitem[DaKr2]{DaKr2} \textsc{Daletskii, \,Yu.L.} and \textsc{Krein,\,S.G.}
\textsl{Integrirovanie i differentsirovanie funktsi\u{\i}
ermitovykh operatorov i prilozhenie k teorii vozmushchini\u{\i}.}
(in Russian). Voronezh. Trudy seminara po funktsional'nomu
analizu, vol. \textbf{1} (1956), pp.\,81 - 105. English transl.
\textsl{Integration and differentiation of functions of hermitian
operators and applications to the theory of perturbations.} Amer.
Math. Soc. Transl. (Ser.\textbf{2}), \textbf{47} (1965), pp.\,1 -
30.

\bibitem[GoKr]{GoKr} \textsc{Gohberg,\,I.Ts.} and \textsc{M.G.\,Krein}.
\textsl{Teoriya Vol'terrovykh Operatorov v Gil'bertovom
Prostranstve i ee Prilozheniya} (In Russian).Mauka, Moskow 1977.
English transl.: \textsl{Theory and applications of Volterra
operators in Hilbert space.}  (Translations of Mathematical
Monographs, \textbf{24}). American Mathematical Society,
Providence, R.I. 1970 x+430 pp.

\bibitem[Had]{Had} \textsc{Hadamard,\,J.}
\textsl{Th\'{e}or\`{e}mes sur la s\'{e}ries enti\`{e}res}. Acta
Math., \textbf{22} (1899), pp.\,55 - 63.

\bibitem[Hal1]{Hal1} \textsc{Halmos,\,P.R.}
\textsl{Finite-Dimensional Vector Spaces.} (Annals of Math.
Studies, \textbf{7}), Princeton University Press, Princeton, N.J.,
1948,. Second edition: Van Norstand, Princeton, N.J., 1958.

\bibitem[HorR1]{HorR1} \textsc{Horn,\,R.A.} \textsl{The Hadamard
product.} In: Matrix theory and applications (Phoenix, AZ, 1989),
pp.\,87 - 169, Proc. Sympos. Appl. Math., \textbf{40}, Amer. Math.
Soc., Providence, RI, 1990.

\bibitem[HorR2]{HorR2} \textsc{Horn,\,R.A.}
\textsl{Topics in Matrix Analysis}.
Cambridge University Press, Cambridge 1991, i -- viii, 607 pp.

\bibitem[Pel]{Pel} \textsc{Peller,\,V.V.} \textsl{Gankelevy operatory v
teorii vozmushcheni{\u\i} samosopryazhennykh operatorov}
(Russian). Funktsional'n. Analiz i ego prilozh., \textbf{19}:2
(1985), pp.\,37-51. English transl.: \textsl{Hankel operators in
the perturbation theory of unitary and self-adjoint operators.}
Funct. Anal. and Appl., \textbf{19} (1985), pp.111-123.

\bibitem[Pie1]{Pie1} \textsc{Pietsch,\,A.} \textsl{Absolut \(p\)-summierende
Abbildungen in normierten R\"{a}umen} (In German). Studia Math.
\textbf{28} (1967), pp.\,333 - 353.

\bibitem[Pie2]{Pie2} \textsc{Pietsch,\,A.} \textsl{Operator ideals}.
 (Mathematische Monographien
[Mathematical Monographs], \textbf{16}). Deutscher Verlag der
Wissenschaften, Berlin 1978. 451 pp.;
\textsl{Operator ideals.} (North-Holland Mathematical Library,
\textbf{20}.) North-Holland Publishing Co.,
Amsterdam\(\boldsymbol{\cdot}\)New York, 1980. 451 pp.

\bibitem[Sch:\,Ges]{Sch:Gesd} \textsc{Schur,\,I.}:
\textsl{Gesammelte Abhandlungen [\textsl{Collected Works}]. Vol.
I, II, III.} Springer-Verlag, %
Berlin\(\cdot\)Heidelberg\(\cdot\)New\,York 1973.

\bibitem[Sch4]{Sch4b} \textsc{Schur,\,I.}: \textsl{Bemerkungen zur
Theorie der beschr\"ankten Bilinearformen mit unendlich vielen
Ver\"anderlichen.}  [\textsl{Remarks on the theory of bounded
bilinear forms with infinitely many variables} - in German].
Journ. f\"ur reine und angew. Math., \textbf{140} (1911), pp.\,1 -
28. \footnotesize Reprinted in: \small \cite{Sch:Gesd}, Vol.\,I,
pp.\,464 - 491.

\bibitem[Sty]{Sty} \textsc{Styan,\,G.} \textsl{Hadamard products and
multivariate statictical analysis.} Linear Algebra and its
Applications, \textbf{6} (1973), pp.\,217 - 240.

\bibitem[Tit]{Tit} \textsc{Titchmarsh,\,E.C.} \textsc{The Theory of Functions.}
The Clarendon Press, Oxford 1932.
\end{thebibliography}

\vspace{3.0ex}
\begin{thebibliography}{MoMo2}
\small

\bibitem[AlU]{AlU} \textsc{Alberti,\,P.M.} and \textsc{A.\,Uhlmann} 
\textsl{Stochasticity and partial order. Doubly stochastic maps and unitary 
mixing.} Mathematics and its Applications, {\bf 9}. D. Reidel,  
Dordrecht-Boston,
 1982; Mathematische Monographien [Mathematical Monographs], {\bf 18}. VEB 
Deutscher Verlag der Wissenschaften, Berlin, 1981.

\bibitem[ArnB]{ArnB} \textsc{Arnold,\,B.C.} \textsl{Majorization
 and Lorenz order:
a brief introduction.} (Lect. Notes in Statist., \textbf{43}).
Springer-Verlag, Berlin 1987, vi+122 pp.

\bibitem[Ati]{Ati} \textsc{Atiyah,\,M.F.} %
\textsl{Convexity and commuting Hamiltonians.}
Bull. London. Math. Soc., \textbf{14}:1 (1982), pp.\,1-15.

\bibitem[BeBe]{BeBe} \textsc{Beckenbach,\,E.F.} and
\textsc{R.\,Bellman.} \textsl{Inequalities.} (Ergebnisse der
Mathem. und ihrer Grenzgeb., Neue Folge, \textbf{30}),
Springer-Verlag,
Berlin\(\boldsymbol{\cdot}\)G\"ottingen\(\boldsymbol{\cdot}\)Heidelberg
1961.

\bibitem[Ber]{Ber} \textsc{Berge, C.} \textsl{Th\'eorie des Graphes
et ses Applications}. Dunod, Paris 1958 (French).
 English transl.:
\textsl{The Theory of Graphs and its Applications}.
Methuen\,\&\,Co, London 1962.

\bibitem[Birk1]{Birk1} \textsc{Birkhoff,\,G.} \textsl{Tres notas
sobre el algebra lineal} [\textsl{Three observations on linear
algebra} - in Spanish]. Universidad National de Tucum\'an,
Revista, Serie \textbf{A}, \textbf{5} (1946), pp.\,147-151.

\bibitem[Birk2]{Birk2} \textsc{Birkhoff,\,G.} \textsl{Lattice Theory.}
(Amer. Math. Soc. Coll. Publ. \textbf{25}.) Amer. Math. Soc.,
Providence, RI, 1948.

\bibitem[BFR]{BFR} \textsc{Bloch,\,A. M.,  H.\,Flaschka} and
\textsc{T.A.\,Ratiu.} \textsl{Schur-Horn-Kostant convexity theorem
for the diffeomorphism group of the annulus}. Invent. Math.
\textbf{113}:3 (1993), pp.\,511--529.

\bibitem[Dan]{Dan} \textsc{Dantzig,\,G.B.} \textsl{Application of the
symplex method to a transportation problem}. Chapter XXIII in:
\textsl{Activity Analysis of Production and Allocation.}
(\textsc{Koopmans,\,T.C.}\,-\,ed.), Wiley, New\,York 1951.

\bibitem[Fel]{Fel} \textsc{Feller,\,W.} \textsl{An Introduction to
Probability Theory and Its Applications,} 1st ed.,
Vol.\,\textbf{1}, Wiley, New\,York 1950.

\bibitem[GlLy]{GlLy} \textsc{Glazman,\,I.M.} and
\textsc{Yu.I.\,Lyubich.} \textsl{Konechnomerny\u{\i}
line\u{\i}ny\u{\i} analiz v zadachakh} (Russian). Nauka, Moskow
1969, 475 pp. English transl.: \textsl{Finite-dimensional linear
analysis: a systematic presentation in problem form.} The M.I.T.
Press, Cambridge, MA\(\boldsymbol{\cdot}\)London, 1974. xvi+520
pp. French transl.: \textsl{ Analyse lin\'eaire dans les espaces
de dimensions finies: manuel en probl\`emes.}  Mir, Moscow, 1972.
400 pp.

\bibitem[GoKr1]{GoKr1a} \textsc{Gohberg,\,I.Ts.} and \textsc{M.G.\,Krein}.
\textsl{Vvedenie v Teoriyu Line\u{\i}nykh Nesamosopryazhennykh
Operatorov} (Russian). Nauka, Moskow 1965, 448 pp. English
Transl.: \textsl{Introduction to the Theory of Linear
Non-Selfadjoint Operators.} (Transl. Math. Monogr. \textbf{18}).
Amer. Math. Soc., Providence, R.I., 1969.

\bibitem[GuSt1]{GuSt1} \textsc{Guillemin,\,V.} and \textsc{S.\,Sternberg.}
\textsl{Convexity properties of the moment mapping.} Invent. Math.
\textbf{67} (1982), pp.\,491-513.

\bibitem[GuSt2]{GuSt2} \textsc{Guillemin,\,V.} and \textsc{S.\,Sternberg.}
\textsl{Convexity properties of the moment mapping. II.}
Invent.\,Math. \textbf{77} 1984), pp.\,533-546.

\bibitem[HalM]{HalM} \textsc{Hall,\,M.,\,\footnotesize{Junior}}\,
\textsl{Combinatorial Theory}. Blaisdell Publishing Co., Waltham\,MA\,%
\(\boldsymbol{\cdot}\)To- ronto\,\(\boldsymbol{\cdot}\)London
1967, x+310\,pp.

\bibitem[HalP]{HalP} \textsc{Hall,\,Ph.}
\textsl{On representatives of subsets.}
Journ. Lond. Math. Soc, \textbf{10} (1935), pp.\,26-30.

\bibitem[Hal2]{Hal2} \textsc{Halmos,\,P.R.}
\textsl{Bounded Integral Operators in
\(L^2\) Spaces.} Springer Verlag, Berlin 1978.

\bibitem[Har:Col]{Har:Col} \textsc{G.H.\,Hardy.}
\textsl{Collected Papers. Vol.\,2.} Clarendon
Press, Oxford 1967.

\bibitem[HiOl]{HiOl} \textsc{Hilgert,\,J.} and
\textsc{G.\,{\`O}lafsson.} \textsl{Causal Symmetric Spaces.
Geometry and Harmonic Analysis.} (Perspectives in Mathematics,
\textbf{18}). Academic Press,
San\,Diego\(\boldsymbol{\cdot}\)London 1997, i-ivx+286 pp.

\bibitem[HLP]{HLPb} \textsc{Hardy,\,G.H., J.E.\,Littlewood,} and
\textsc{G.\,P\'{o}lya.} \textsl{Inequalities}. 1st ed., 2nd ed.
Cambridge Univ. Press, London\(\boldsymbol{\cdot}\)New\,York 1934,
1952.

\bibitem[HLP1]{HLP1} \textsc{Hardy,\,G.H., J.E.\,Littlewood} and
\textsc{G.\,P\'{o}lya}. \textsl{Some simple inequalities satisfied by
convex functions.} Messenger of Mathematics, \textbf{58},
pp.\,145-152. Reprinted in \cite{Har:Col}, pp.\,500-508.

\bibitem[HNP]{HNP} \textsc{Hilgert,\,J., K.-H.\,Neeb}, and
\textsc{W.\,Plank.} \textsl{Symplectic convexity theorems and
coadjoint orbits.} Compos. Math., \textbf{94} (1994), pp.129-180.


\bibitem[HoJo]{HoJo} \textsc{Horn,\,R.A.}, and \textsc{Ch.R.\,Johnson}.
\textsl{Matrix Analysis}. Cambrigde University Press,
Cambridge\(\boldsymbol{\cdot}\)London\(\boldsymbol{\cdot}\)New\,York
1986.

\bibitem[HorA1]{HorA1} \textsc{Horn,\,A.} \textsl{Doubly stochastic
matrices and the diagonal of the rotation matrix.} Amer. J. Math.
\textbf{76} (1954), 620 - 630.

\bibitem[HorA2]{HorA2} \textsc{Horn,\,A.} \textsl{On the eigenvalues of a matrix
with prescribed singular values.} Proc. Amer. Math. Soc.,
\textbf{5}:1 (1954), pp.\,4-7.

\bibitem[Kos]{Kos} \textsc{Kostant,\,B.} \textsl{On convexity, the Weyl group
and the Iwasawa decomposition.} Ann. Sci. \'Ecole Norm. Sup.
(Ser.\,4), \textbf{6}\,(1973), pp.\,413-455.

\bibitem[Lor]{Lor} \textsc{Lorenz,\,M.O.}
\textsl{Methods of measuring concentrations of wealth.} Journ.
Amer. Statist. Assoc., \textbf{9} (1905), pp. 209-219.

\bibitem[MaOl]{MaOl} \textsc{Marshall,\,A.W.} and
\textsc{I.\,Olkin.} \textsl{Inequalities: Majorization and Its
Applications.} Academic Press,
New\,York\(\cdot\)London\(\cdot\)Toronto 1979.

\bibitem[Mark]{Mark} \textsc{Markus,\,A.S.}
\textsl{Sobstvennye i singulyarnye chisla
summy i proizvedeniya line\u{\i}nykh operatorov} (Russian).
Uspekhi Matem. Nauk, \textbf{19}:4 (1964), pp.\,93-123.  English
transl.: \textsl{The eigen- and singular values of the sum and
product of linear operators}. Russian Math. Surveys, \textbf{19}
(1964), pp.\,91-120.

\bibitem[Mir]{Mir} \textsc{Mirsky,\,L.}
\textsl{Results and problems in the theory of
doubly-stochastic matrices.} Zeitschr. f\'ur die
Wahrscheinlichkeitstheorie und Verw. Gebiete, \textbf{1}
(1962/1963), 319-334.

\bibitem[Mit]{Mit} \textsc{Mityagin,\,B.S.}
\textsl{Interpolyatsionnaya teorema dlya
modulyarnykh prostranstv} (Russian). Matem. Sbornik (N.S.),
\textbf{66} (1965), pp.\,473-482. English transl.: \textsl{An
interpolation theorem for modular spaces.} In:
\textsl{Interpolation spaces and allied topics in analysis (Lund,
1983)}, Lecture Notes in Math., \textbf{1070}, Springer, Berlin,
1984, pp.\,10-23.

\bibitem[Muir]{Muir} \textsc{Muirhead,\,R.F.}
\textsl{Some methods applicable to identities and inequalities of
symmetric algebraic functions of \(n\) letters.} Proc. Edinburgh
Math. Soc., \textbf{21} (1903), pp.144-157.

\bibitem[NeuA]{NeuA} \textsc{Neumann,\,A.}
\textsl{An infinite dimensional version of the
Schur-Horn convexity theorem.} Journ. of Funct. Anal.,
\textbf{161} (1999), pp.\,418 - 451.

\bibitem[NeuJ1]{NeuJ1} \textsc{Neumann, John von.}
\textsl{A certain zero-sum two-person game equivalent to the
optimal assignement problem.} In: \textsl{Contributions to the
Theory of Games. Vol. \textup{II}} (\textsc{Kuhn,\,H.W.} and
\textsc{A.W.\,Tucker} - editors). (Annals Math. Study
\textbf{28}), Princeton Univ. Press, Princeton 1953,  pp.\,5-12.
Reprinted in: \cite{NeuJ2}, pp.\,44-49.

\bibitem[NeuJ2]{NeuJ2} \textsc{Neumann, John von.}
\textsl{Collected Works. Vol.\textup{VI}: Theory of games,
astrophysics, hydrodynamics and meteorology.}  General
editor:\,\textsc{A.H.\,Taub}, Pergamon Press,
Oxford\(\boldsymbol{\cdot}\)London\(\boldsymbol{\cdot}\)New
York\(\boldsymbol{\cdot}\)Paris: 1963, x+ 538 pp.

\bibitem[Ostr]{Ostr} \textsc{Ostrowski,\,A.}
\textsl{Sur quelques applications des fonctions convexes and
concaves au sens de I.\,Schur} (French). Journ. de Math\'ematiques
Pures et Appliqu\'ees, \textbf{31} (1952), pp.\,253-292.

\bibitem[PPT]{PPT} \textsc{Pe{\v{c}}ari{\'c},\,J.E.,
F.\,Proschan} and \textsc{Y.L.\,Tong}. \textsl{Convex functions,
partial ordering, and statistical applications}. (Mathematics in
Science and Engineering, \textbf{187}). Academic Press, Boston,
MA, 1992, xiv+467 pp.

\bibitem[Rad]{Rad} \textsc{Rado,\,R.}  \textsl{An inequality.}
Journ. Lond. Math. Soc., \textbf{27} (1952), pp.\,1-6.

\bibitem[Rom1]{Rom1} \textsc{Romanovsky,\,V.}
\textsl{Sur les z{\`e}ros des matrices stochastiques.} Compt. Rend.
Acad. Sci. Paris, \textbf{192} (1931), pp.\,266-269.

\bibitem[Rom2]{Rom2} \textsc{Romanovsky,\,V.}
\textsl{Recherches sur les chaines de Markoff.} Acta Math.,
\textbf{66} (1930), pp.\,147-251.

\bibitem[Sch18]{Sch18a} \textsc{Schur,\,I.}: \textsl{\"Uber eine Klasse
von Mittelbildungen mit Anwendungen auf die Determinantentheorie.}
[\textsl{On a class of averaging mappings with applications to the
theory of determinants} - in German]. Sitzungsberichte der
Berliner Mathematischen Gesellschaft, \textbf{22} (1923), pp.\,9 -
20. \footnotesize Reprinted in: \small \cite{Sch:Gese}, Vol.\,II,
pp.\,416 - 427.

\bibitem[Sch:\,Ges]{Sch:Gese} \textsc{Schur,\,I.}:
\textsl{Gesammelte Abhandlungen [\textsl{Collected Works}]. Vol.
I, II, III.} %
Springer-Verlag,
Berlin\(\boldsymbol{\cdot}\)Heidelberg\(\boldsymbol{\cdot}\)New\,York,
1973.

\end{thebibliography}

\vspace*{3.0ex}
\begin{thebibliography}{MoMo2}
\small
\bibitem[Bar]{Bar} \textsc{Bari,\,N.K.}
 \textsl{Trigonometricheskie ryady} (Russian). Fiz.-Mat.\,Giz, Moscow
  1961, 936 pp. English transl.: \textsl{A treatise on trigonometric
series. Vols. I, II.} A Pergamon Press Book The Macmillan Co., New
York 1964 , Vol. I: xxiii+553 pp. Vol. II: xix+508 pp.

\bibitem[Bend]{Bend} \textsc{Bendixson,\,F.} \textsl{Sur les racines d'une
\'{e}quation fondamentale.} Acta Math., \textbf{25} (1902),
pp.359\,-\,365.

\bibitem[Carl1]{Carl1} \textsc{Carleman,\,T.} \textsl{Sur le genre du
d\'enominateur \(D(\lambda)\) de Fredholm.} Arkiv f\"{o}r
Mathematik, Astronomi och Fysik, \textbf{12} (1917), pp.\,1\,-\,5.
Reprinted in \cite{Carl3}, pp.\,1\,-\,5.

\bibitem[Carl2]{Carl2} \textsc{Carleman,\,T.} \textsl{\"Uber die
Fourierkoeffizienten einer stetigen Funktion.} Acta Math.,
\textbf{41} (1918), pp.\,377\,-\,384. Reprinted in \cite{Carl3},
pp.\,15\,-\,22.

\bibitem[Carl3]{Carl3} \textsl{\'Edition Compl\`ete des Articles de
 Torsten Carleman.} Publ\'ee par l'Institute Math\'ematique
 Mittag-Leffler. Litos Reprotryk, Malm\"o 1960.

\bibitem[Chang]{Chang} \textsc{Chang,\,Shih-Hsun.} %
\textsl{On the distribution of the characteristic values and
singular values of linear integral equations.}
Trans.\,Amer.\,Math.\,Soc., \textbf{67}:2 (1949),
pp.\,351\,-\,367.

\bibitem[Fred]{Fred} \textsc{Fredholm,\,E.I.} \textsl{Sur une classe
d'\'equations fonctionnelles}. Acta Math., \textbf{27},
pp.\,365\,-\,396.

\bibitem[GoKr1]{GoKr1b} \textsc{Gohberg,\,I.Ts.} and \textsc{M.G.\,Krein}.
\textsl{Vvedenie v Teoriyu Line\u{\i}nykh Nesamosopryazhennykh
Operatorov} (Russian). Nauka, Moskow 1965, 448 pp. English
Transl.: \textsl{Introduction to the Theory of Linear
Non-Selfadjoint Operators.} (Transl. Math. Monogr. \textbf{18}).
Amer. Math. Soc., Providence, R.I., 1969.

\bibitem[Hir]{Hir} \textsc{Hirsch,\,A.} \textsl{Sur les racines d'une
\'{e}quation fondamentale.} Acta Math., \textbf{25} (1902),
pp.367\,-\,370.

\bibitem[Lal]{Lal} \textsc{Lalesco,\,T.}
\textsl{Sur l'ordre de la fonction
enti\`ere \(D(\lambda)\) de Fredholm}. Comp. Rend.  des Seances
de'l Academie des Science, Paris, 25 November (1907), p.\,906.

\bibitem[Lid]{Lid} \textsc{Lidski\u{\i},\,V.B.}
\textsl{Nesamosopryazhennye
operatory, imeyushchie sled} (Russian). Dokl. Akad. Nauk SSSR,
\textbf{125} (1959), pp.\,485\,-\,488. English transl.:
\textsl{Non-selfadjoint operators with a trace}. Amer. Math. Soc.
Thansl. (\textbf{2}) \textbf{47} (1965), 43\,-\,46.

\bibitem[Sch2]{Sch2a} \textsc{Schur,\,I.}: \textsl{\"Uber die
charakteristischen Wurzeln einer linearen Substitution mit einer
Anwendung auf die Theorie der Integralgleichungen.} [\textsl{On the 
characteristic roots of a linear substitution with an application
to the theory of integral equations} - in German].  Mathematische
Annalen, \textbf{66} (1909), pp.\,488 - 510. \footnotesize
Reprinted in: \small \cite{Sch:Gesf}, Vol.\,I, pp.\, 272 - 294.

\bibitem[Sch:\,Ges]{Sch:Gesf} \textsc{Schur,\,I.}:
\textsl{Gesammelte Abhandlungen}. [\textsl{Collected Works}]. Vol.
I, II, III. Springer-Verlag,
Berlin\(\cdot\)%
Heidelberg\(\cdot\)New\,York 1973.

\bibitem[TuAi]{TuAi} \textsc{Turnbbull,\,H.W.} and
\textsc{A.C.\,Aitken.}
\textsl{An Introductian to the Theory of Canonoical Matrices.}
Blackie\ \&\  %
Son, Ltd., London and Glasgow, First printed, 1932; Reprinted
1945, 1948, 1950, 1952.
\bibitem[Wey]{Wey} \textsc{Weyl,\,H.} \textsl{Inequalities between
the two kinds of
eigenvalues of a linear transformation.} Proc. Nat. Acad. Sci.
USA, \textbf{35} (1949), pp.\,408\,-\,411. Reprinted in:
\cite{Wey:Ges}, Band IV, pp.\,390\,-\,393.

\bibitem[Wey:Ges]{Wey:Ges} \textsc{Weyl,\,H.}
\textsl{Gesammelte Abhandlungen
\textup{[}Collected Works\textup{]}. Vol.\,I, II, III, IV.}
Springer\,-\,Verlag,
Berlin\(\boldsymbol{\cdot}\)Heidelberg\(\boldsymbol{\cdot}\)New\,York
1968.

\bibitem[Zyg]{Zyg} \textsc{Zygmund,\,A.} \textsl{Trigonometric series.}
 2nd ed. Vols. I, II. Cambridge University Press, New York 1959.
 Vol. I. xii+383 pp.; Vol. II. vii+354 pp.

\end{thebibliography}

\vspace{3.0ex}
\begin{thebibliography}{MoMo2}
\small
\bibitem[ArSm]{ArSm} \textsc{Aronszajn,\,N.} and \textsc{K.T.\,Smith.}
\textsl{Invariant subspaces of completely continuous operators.}
Annales of Math., \textbf{60} (1954), pp.\,345\,-\,350.

\bibitem[Bran]{Bran}  \textsc{Branges,\,Louis de.} \textsl{Hilbert Spaces of
Entire Functions}, Prentice Hall, Englewood Cliffs, NJ, 1968.

\bibitem[Brod1]{Brod1} \textsc{Brodski\u{\i},\,M.S.} \textsl{Ob integral'nom
predstavlenii ogranichennykh nesamosopryazhennykh operatorov s
veshchestvennym spektrom} (Russian). [Integral representations of
bounded non-selfadjoint operators with a real spectrum]. Dokl.
Akad. Nauk SSSR, \textbf{126} (1959), 1166--1169.

\bibitem[Brod2]{Brod2} \textsc{Brodski\u{\i},\,M.S.} \textsl{O treugol'nom
predstavlenii vpolne nepreryvnykh operatorov s odno\u{\i}
tochko\u{\i} spectra} (Russian). Uspehi Mat. Nauk \textbf{16}:1
(1961), pp.\,135--141. English transl.:{On the triangular
representation of completely continuous operators with one-point
spectra}.  Amer. Math. Soc. Transl. (\textsf{2}), \textbf{47}
(1965), pp.\,59\,-\,65.

\bibitem[Brod3]{Brod3} \textsc{Brodski\u{\i},\,M.S.} \textsl{O treugol'nom
predstavlenii nekotorykh operatorov s vpolne nepreryvno\u{\i}
mnimo\u{\i} chast'yu} (Russian). Dokl. Akad. Nauk SSSR,
\textbf{133} (1960), pp.\, 1271\,-\,1274. English transl.:
\textsl{Triangular representation of some operators with
completely continuous imaginary part}. Soviet Math. Dokl.,
\textbf{1} (1960), pp.\,952\,-\,955.
\bibitem[Brod4]{Brod4} \textsc{Brodski\u{\i},\,M.S.} \textsl{Treugol'nye
i Zhordanovy predstavleniya line\u{\i}nykh operatorov.}
(Sovremennye problemy matematiki) (Russian). Nauka, Moscow 1969,
287 pp. English translation: \textsl{Triangular and Jordan
Representations of Linear Operators.} (Transl. of Math.
Monographs, \textbf{32}). Amer. Math. Soc., Providence, RI, 1971.

\bibitem[BroBr]{BroBr} \textsc{Brodski\u{\i},\,V.M.} and
\textsc{M.S.\,Brodski\u{\i}}. \textsl{Ob abstraktnom treugol'nom
predstavlenii line\u{\i}nykh ogranichennykh operatorov i
mul'tiplikativnom razlozhenii sootvetstvuyushchikh im
kharakteristicheskikh funktsi\u{\i}} (Russian). Dokl. Akad. Nauk
SSSR, \textbf{181}:3 (1968), pp.\,511\,-\,514. English transl.:
\textsl{The abstract triangular representation for bounded linear
operators and multiplicative expansions of their respective
characteristic functions}. Soviet Math. Dokl., \textbf{9} (1968),
pp.\,846\,-\,850.

\bibitem[BroLi]{BroLi} \textsc{Brodski\u{\i},\,M.S.} and
\textsc{M.S.\,Livshits} (=\textsc{M.S.\,Liv\v{s}ic}).
\textsl{Spektral'ny\u{\i} analiz nesamosopryazhennykh operatorov i
promezhutochnye sistemy}(Russian). Uspehi Mat. Nauk (N.S.)
\textbf{13}:1 (1958), pp.\.3\,-\,85. English transl.:
\textsl{Spectral analysis of non-selfadjoint operators and
intermediate systems.} Amer. Math. Soc. Transl. (\textsf{2})
\textbf{13} (1960), pp.\,265\,-\,346.

\bibitem[Dav]{Dav} \textsc{Davidson,\,K.R.} \textsl{Nest
Algebras.}\,{\footnotesize Triangular forms for operator algebras
on Hilbert space.} (Pitman Research Notes in Mathematics Series,
\textbf{191}). Longman, Harlow \& Wiley New\,York 1988.

\bibitem[GoGoK]{GoGoK} \textsc{Gohberg,\,I.Ts., S.\,Goldberg} and
\textsc{M.\,Kaashoek.} \textsl{Classes of Linear Operators. Vol.
\textup{II}}. (Operator Theory: Advances and Applications.
\textbf{OT\,63}). Birkh\"auser Verlag, Basel 1993.

\bibitem[GoKr2]{GoKr2} \textsc{Gohberg,\,I.Ts.} and \textsc{M.G.\,Krein}.
\textsl{Theoriya Vol'terrovykh Operatorov v Gil'bertovom
Prostranstve i ee Prilozheniya} (Russian), Nauka, Moscow 1967.
English transl.: \textsl{Theory and Applications of Volterra
Operators in Hilbert Space}. (Transl. of Mathem. Monogr.,
\textbf{24}), Amer. Math. Soc., Providence 1970.

\bibitem[GoKr3]{GoKr3} \textsc{Gohberg,\,I.Ts.} and \textsc{M.G.\,Krein}.
\textsl{O vpolne nepreryvnykh operatorakh so spectrom
sosredotochennym v nule} (Russian). [Completely continuous
operators with a spectrum conzentrated at zero]. Doklady Akad.
Nauk SSSR, \textbf{128}:2 (1959), pp.\,227\,-\,230.


\bibitem[GoKr4]{GoKr4} \textsc{Gohberg,\,I.Ts.} and \textsc{M.G.\,Krein}.
\textsl{K teorii treugol'nykh predstavleni\u{\i}
nesamosopryazhennykh operatorov} (Russian). Doklady Akad. Nauk
SSSR, \textbf{137}:5 (1961), pp.\,1034\,-\,1037. English transl.:
\textsl{On the theory of triangular representations of linear
operators}. Soviet Math., Doklady, \textbf{2} (1961),
pp.\,392\,-\,396.

\bibitem[GoKr5]{GoKr5} \textsc{Gohberg,\,I.Ts.} and \textsc{M.G.\,Krein}.
\textsl{O vol'terrovykh operatorakh s mnimo\u{\i} komponento\u{\i}
togo ili inogo klassa} (Russian). Doklady Akad. Nauk SSSR,
\textbf{139}:4 (1961), pp.\,779\,-\,782. English transl.:
\textsl{Volterra operators with imaginary component in one class
ore another}. Soviet Mat., Doklady, \textbf{2} (1961),
pp.\,983\,-985.

\bibitem[Herr]{Herr} \textsc{Herrero,\,D.A.} \textsl{Triangular operators.}
Bull. Lond. Math. Soc., \textbf{23} (1991), pp.\,513 - 554.


\bibitem[Kats]{Kats} \textsc{Katsnelson,\,V.} \textsl{Right and left joint
system representation of a rational matrix function in general
position.}

\bibitem[LaPhi]{LaPhi} \textsc{Lax,\,P.} and \textsc{R.\,Phillips.}
\textsl{Scattering Theory}. Academic Press,
New\,York\(\boldsymbol{\cdot}\)London 1967.

\bibitem[Liv1]{Liv1} \textsc{Liv\v{s}ic, M.S.} (=\textsc{Livshits,\,M.S.})
\textsl{Ob odnom classe line\u{i}nykh operatorov v gil'bertovom
prostranstve} (Russian). Matem. Sbornik, \textbf{19} (1946),
pp.\,239\,-\,260. English transl.:

\bibitem[Liv2]{Liv2} \textsc{Liv\v{s}ic, M.S.} (=\textsc{Livshits,\,M.S.})
\textsl{Izometricheskie operatory s ravnymi defektnymi chislami,
kvazi-unitarnye operatory} (Russian). Matem. Sbornik, \textbf{26}
(1950), pp.\,247\,-\,264. English translation:

\bibitem[Liv3]{Liv3} \textsc{Liv\v{s}ic, M.S.} (=\textsc{Livshits,\,M.S.})
\textsl{O privedenii line\u{\i}nykf ermitovykh operatorov k
"treugol'nomu" vidu.} Doklady Akad. Nauk SSSR, \textbf{84}:5,
pp.\,873\,-\,876.

\bibitem[Liv4]{Liv4} \textsc{Liv\v{s}ic, M.S.} (=\textsc{Livshits,\,M.S.})
\textsl{O rezol'vente line\u{\i}nogo nesimmetricheskogo
operatora.} Doklady Akad. Nauk SSSR, \textbf{84}:6,
pp.\,1131\,-\,1134.

\bibitem[Liv5]{Liv5}   \textsc{Liv\v{s}ic, M.S.} (=\textsc{Livshits,\,M.S.})
\textsl{O spektral'nom razlozhenii line\u{\i}nykh
nesamosopryazhennykh operatorov.} Matem. Sbornik, \textbf{34}
(1954), pp.\,145\,-\,199 (Russian),. English transl.: \textsl{On
the spectral resolution  of linear non-selfadjoint operators.}
Amer. Math. Soc. Transl. (Ser.\textsf{2}), \textbf{5} (1957), pp.
67\,-\,114.

\bibitem[Liv6]{Liv6}   \textsc{Liv\v{s}ic, M.S.} (=\textsc{Livshits,\,M.S.})
\textsl{}Zhurnal Experimental'no\u{\i} i Theor. Phys.,
\textbf{31}:1, pp.\,121-131 (Russian). English transl.:
\textsl{The application of non-self-adjoint operators to
scattering theory.} Soviet Physics JETP, \textbf{4}:1 (1957),
pp.\,91\,-\,98.

\bibitem[Liv7]{Liv7}  \textsc{Liv\v{s}ic, M.S.} (=\textsc{Livshits,\,M.S.})
 \textsl{Metod nesamosopryazhennykh operatorov v teorii rasseyaniya} (Russian).
Uspekhi Matem. Nauk, \textbf{12}:1 (1957), pp.\,212 - 218
(Russian). Engl. transl.: \textsl{The  method of non-self-adjoint
operators in dispersion theory.} Amer. Math. Soc. Transl.
(\textsf{2}), \textbf{16} (1960), pp.427-434.

\bibitem[Liv8]{Liv8} \textsc{Liv\v{s}ic, M.S.}. (\,=\,\textsc{Livshits,\,M.S.})
\textsl{Operatory, Kolebanija, Volny (Otkrytye sistemy)}, Nauka,
Moscow, 1966 (Russian). English transl.: \textsl{Operators,
oscillations, waves \textup{(}open systems\textup{)}}.
(Translations of Mathematical Monographs, Vol. \textbf{34}.)
American Mathematical Society, Providence, R.I., 1973.

\bibitem[Mats1]{Mats1} \textsc{Matsaev,\,V.I.} (\,=\,\textsc{Macaev,\,V.I.})
\textsl{Ob odnom klasse vpolne nepreryvnykh operatprov} (Russian).
Doklady Akad. Nauk SSSR, \textbf{139}:3 (1961), pp.\,548\,-\,551.
English transl.: \textsl{On a class of completely continuous
operators.} Soviet Math., Doklady, \textbf{2} (1961),
pp.\,972\,-\,975.

\bibitem[Mats2]{Mats2} \textsc{Matsaev,\,V.I.} (\,=\,\textsc{Macaev,\,V.I.})
\textsl{O vol'terrovykh operatorakh, poluchaemykh vozmushcheniem
samosopryazhennykh} (Russian). Doklady Akad. Nauk SSSR,
\textbf{139}:4 (1961), pp.\,%
810\,-\,813. English transl.:  \textsl{Volterra operators produced
by perturbation of self-adjoint operators.} English thransl.:
Soviet Math., Doklady, \textbf{2} (1961), pp.\,1013\,-\,1016.

\bibitem[NiVa]{NiVa} \textsc{Nikolski\u{\i},\,N.K.} and
\textsc{V.I.\,Vasyunin}. \textsl{A unified approach to functional
models, and the transcription problem}. The Gohberg anniversary
collection, Vol.II.
 (Operator Theory: Advances and
Applications. \textbf{OT\,41}). Birkh{\"a}user, Basel 1989,
pp.\,403\,-\,434.

\bibitem[OTSTR]{OTSTR} \textsl{Operator Theory, System Theory and
Related Topics}.(\textsl{The Moshe Liv\v{s}ic Anniversary
Volume}). \textsc{Alpay,\,D.} and \textsc{V.\,Vinnikov}-editors.
(Operator Theory: Advances and Applications, \textbf{OT 123}),
Birkh{\"a}user Verlag, Basel, 2001.
\bibitem[Sakh1]{Sakh1} \textsc{Sakhnovich,\,L.A.}
(\,=\,\textsc{Sahnovi\v{c},\,L.A.}) \textsl{O privedenii
nesamosopryazhennykh operatorov k treugol'nomu vidu} [\textsl{The
reduction of non-selfadjoint operators to triangular form} - in
Russian]. Izvestiya Vys\v{s}. U\v{c}ebn. Zaved. Matematika
\textsf{1959}, no.\,1, pp.\,180\,-\,186.

\bibitem[Sakh2]{Sakh2} \textsc{Sakhnovich,\,L.A.}
(\,=\,\textsc{Sahnovi\v{c},\,L.A.}) \textsl{Issledovanie
treugol'no\u{\i} modeli nesamosopryazhennykh operatorov}
[\textsl{A study of the triangular form of non-selfadjoint
operators} - in Russian]. Izvestiya Vys\v{s}. U\v{c}ebn. Zaved.
Matematika \textsf{1959}, no.\,4, pp.\,141\,-\,149.

\bibitem[Sch2]{Sch2b} \textsc{Schur,\,I.}: \textsl{\"Uber die
charakteristischen Wurzeln einer linearen Substitution mit einer
Anwendung auf die Theorie der Integralgleichungen} [\textsl{On the 
characteristic roots of a linear substitution with an application
to the theory of integral equations} - in German]. Mathematische
Annalen, \textbf{66} (1909), pp.\,488 - 510. \footnotesize
Reprinted in: \small \cite{Sch:Gesg}, Vol.\,I, pp.\, 272 - 294.

\bibitem[Sch:\,Ges]{Sch:Gesg} \textsc{Schur,\,I.}:
\textsl{Gesammelte Abhandlungen \textup{[}Collected
Works\textup{]}. Vol. I, II, III.} Springer-Verlag,
Berlin\(\cdot\)%
Heidelberg\(\cdot\)New\,York 1973.

\bibitem[SzNFo]{SzNFo} \textsl{Sz.-Nagy,\,B.} and \textsc{C.\,Foias.}
\textsl{Analyse Harmonique des Op\'erateurs de l'espace  de
Hilbert} (French). Masson and Acad\'emiae Kiado, 1967. English
transl.:\textsl{Harmonic Analysis of Operators in Hilbert Space}. 
North Holland, Amsterdam 1970.

\end{thebibliography}

\begin{thebibliography}{MoMo2} \small
\bibitem[AESW]{AESW} \textsc{Aissen,\,M., A.\,Edrei, I.J.\,Schoenberg}
and \textsc{A.\,Whitney.} \textsl{On the generating functions of
totally positive sequences.} Proc.\,Nat.\.Acad.\,Sci. U.S.A.,
\textbf{37} (1951), pp.\,303\,-\,307.

\bibitem[ASW]{ASW}  \textsc{Aissen,\,M., A., I.J.\,Schoenberg}
and \textsc{A.\,Whitney.} \textsl{On the generating functions of
totally positive sequences I.} J. d'Anal. Math.,\,Jerusalem,
\textbf{2} (1952), pp.\,93\,-\,103.


\bibitem[Dieu]{Dieu} \textsc{Dieudonne,\,J.A.} \textsl{La Theorie
analytique des Polynomes d'une variable : (a coefficients
Quelconques)}. Paris, Gauthier-Villars, 1938.

\bibitem[Edr]{Edr} \textsc{Edrei,\,A.} \textsl{On the generating functions
of totally positive sequences II.} J. d'Anal. Math., Jerusalem,
\textbf{2} (1952), pp.\,104\,-\,109.

\bibitem[Fek]{Fek} \textsc{Fekete,\,M.} \textsl{\"{U}ber ein Problem von
Laguerre} (German). [On one problem of Laguerre.] Rendiconti del
Circolo Matematico di Palermo, \textbf{34} (1912), pp.\,89\,-100.
Reprinted in: \cite{PolCo}, pp.\,1\,-\,32.

\bibitem[FoZe]{FoZe} \textsc{Fomin,\,S.} and \textsc{A.\,Zelevinski.}
\textsl{Total positivity: tests and parametrizations}. Math.
Intelligencer, \textbf{22}:1 (2000), pp.\,23\,-\,33.

\bibitem[GaKr1]{GaKr1} \textsc{Gantmacher,\,F.R.} and
\textsc{M.G.\,Krein}. \textsl{Sur les matrices oscillatoires}
(French). [On oscillatory matrices]. Compt. Rend. Acad. Sci
(Pqris). \textbf{201} (1935), pp.\,577\,-\,579.

\bibitem[GaKr2]{GaKr2} \textsc{Gantmacher,\,F.R} and
\textsc{M.G.\,Krein}. \textsl{Sur les Matrices Compl\`{e}ment
Non-n\'{e}gatives et Oscillatoires} (French). [On completely
non-negative and oscillatory matrices] Compositio  Math.,
\textbf{4} (1937), pp.\,445\,-\,476.

\bibitem[GaKr3]{GaKr3} \textsc{Gantmacher,\,F.R.} and
\textsc{M.G.\,Krein.} \textsl{Ostsillyatsionnye Matritsy i Yadra i
Malye Kolyebaniya Mekhanicheskikh Sistem} (Russian) [Oscillatory
Matrices and Kernels and small Vibrations of Mechanical systems].
GITTL, Moskow 1950. English transl: German transl.:

\bibitem[GaMic]{GaMic} \textsc{Gasca,\,M.} and
\textsc{C.A.\,Micchelli} (eds). \textsl{Total Positivity and its
Applications.} (Mathematics and its Applications, \textbf{359}).
Kluver Academic Publishers, Dordrecht, 1996. x+518\,pp.

\bibitem[HiWi]{HiWi} \textsc{Hirschman,\,I.I.} and {D.V.\,Widder}.
\textsl{The Convolution Transform.} Princeton Univ. Press,
Princeton, NJ, 1955, x\,+\,268 pp. Russian transl.:
\textsl{Preobrazovaniya Tipa Svertki}. \ Izdat. Inostranno\u{\i}
Literatury, Moscow 1958, 312 pp.

\bibitem[Kar]{Karl} \textsc{Karlin,\,S.} \textsl{Total
Positivity, Vol.\,I}. Stanford Univ. Press, Stanford, CA, 1968.
xi+576 pp.

\bibitem[KrNa]{KrNe} \textsc{Krein,\,M.G.} and \textsc{M.A.\,Naimark}.
\textsl{Metod Simmetricheskikh i Ermitovykh Form v Teorii
Otdeleniya Korne\u{\i} Algebraicheskikh Uravneni\u{\i}} (Russian).
Nauchno-Issledovatel'ski\u{\i} Institut Matematiki i Mekhaniki pri
Khar'kovskom Universitete. DNTVU, Khar'kov 1936, 43 p. \ \ English
transl.: \textsl{The Method of Symmetric and Hermitian Forms in
the Theory of the Separation of the Roots of Algebraic Equations},
Linear and Multilinear Algebra, \textbf{10} (1981), pp.\,265 -
308.

\bibitem[Kur]{Kur} \textsc{Kurosh,\,A.G.} \textsl{Kurs Vysshe\u{\i} Algebry}
(Russian). Nauka, Moskow 1965. English transl.: {Higher algebra},
Mir, Moscow 1972.

\bibitem[Lag1]{Lagu1} \textsc{Laguerre,\,E.} \textsl{Sur quelques points de la
th\'eorie des \'equations num\'eriques.} Acta Math, \textbf{4}
(1864). Reprinted in \cite{LaguO}, pp.\,184\,-\,206.

\bibitem[Lag2]{Lagu2} \textsc{Laguerre,\,E.} \textsl{Ser les fonctions du
genre z\'ero et du genre un}. Compt.\,Rend.\,Acad.\,Sciences
Paris, \textbf{98} (1882), pp.828\,-\,831. Reprinted in
\cite{LaguO}, pp.\,174\,-\,177.

\bibitem[LagO]{LaguO} \textsc{Laguerre,\,E.} \textsl{Euvres. Tome I}
(French). Paris. Gauthier Villars, 1898, 1-sr ed.; Chelsea, 1972,
2-nd ed.

\bibitem[Lev]{Levi} \textsc{Levin,\,B.Ya.} \textsl{Raspredelenie
korne\u{\i} tselykh funktsi\u{\i}} (in Russian). Gostekhizdat,
Moskow 1956. English  transl.: \textsl{Distribution of Zeros of
Entire Functions.} (Ser.: Transl. of Math. Monogr., \textbf{5},
Amer. Math. Soc., Providence, RI, 1-st ed. 1964, viii\,+\,493 pp.,
2-nd ed. 1980, xii\,+\, 523 pp.

\bibitem[Mar]{Mard} \textsc{Marden,\,M.} \textsl{Geometry of polynomials}
(2-nd ed.). (Ser.: Mathematical surveys and monographs,
\textbf{3}). Amer.\,Math.\,Sos., Providence, Rhode Island 1966,
xiii, 243 p.

\bibitem[Mal]{Malo} \textsc{Malo,\,E.} \textsl{Note sur {\'e}quations
alg{\'e}briques dont toutes les racines sont r{\'e}elles}
(French). [Note on algebraic equations all roots of which are
real]. Journal de Math{\'e}matiques sp{\'e}ciales, (ser.\,4),
t.\,\textbf{4} (1895), p.\,7\,-\,10.

\bibitem[Mot1]{Motz1} \textsc{Motzkin,\,Th.}
\textsl{Beitraege zur Theorie der linearen Ungleichungen}
(German). [Contributions to the theory of linear inequalities].
(Doctoral dissertation. Basel, 1934).  Azriel, Lerusalem , (1936).
73pp. English transl. in \cite{Motz3}, pp.\,1\,-\,80.

\bibitem[Mot2]{Motz2} \textsc{Motzkin,\,Th.} \textsl{Sur les transformations
qui n'augmentent pas le nombre des variations du signe} (French).
[On transformations which don't increase the nomber of the
variations of signs] C. R. Acad. Sci., Paris \textbf{202}, (1936),
pp.\,894-897.

\bibitem[Mot3]{Motz3} \textsc{Motzkin,\,Th.} \textsl{Selected papers}.
(Contemporary Mathematicians). Birkh{\"a}user, Boston, MA,1983.
xxvi+530

\bibitem[Obr]{Obr} \textsc{Obreschkoff,\,N.} \textsl{Verteilung und Berechnung
der Nullstellen reeller Polynome}. Deutscher Verlag der
Wissenschaften, Berlin 1963, viii, 298 p.

\bibitem[Pink]{Pink} \textsc{Pinkus,\,A.} \textsl{Spectral properties of
totally positive kernels and matrices.} In \cite{GaMic},
pp.\,457\,-\,511.

\bibitem[Pol1]{Pol1} \textsc{P\'{o}lya,\,G.} %
\textsl{\"{U}ber Ann\"aherung durch Polynome mit lauter reellen
Wurzeln}. Rend.\,Circ.\,Math.\,Palermo, \textbf{36} (1913),
pp.279\,-\,295.

\bibitem[PolCo]{PolCo} \textsc{P\'{o}lya,\,G.} \textsl{Collected Papers, Vol.
II}. MIT Press, Cambridge, MA, and London, England 1974,
x\,+\,444\,pp.

\bibitem[PS]{PS} \textsc{P\'{o}lya,\,G.} and \textsc{I.\,Schur}:
\textsl{\"Uber zwei Arten von Faktorenfolgen in der Theorie der
algebraischen Gleichungen.} Journ. f\"ur reine und angew. Math.
\textbf{144} (1914), pp.\,89 - 113. \footnotesize Reprinted in:
\small \cite{Sch:Gesh}, Vol.\,II, pp.\,88 - 112. Reprinted also in
\cite{PolCo}, pp.\,100\,-\,124.

\bibitem[PoSz]{PoSz} \textsc{P\`{o}lya,\,G.} and \textsc{G.\,Szeg{\"o}}.
\textsl{Aufgaben and Lehrs{\"a}tze in der Analysis. Band
\textup{II}} (German). (Ser.: Grundlehren der mathematischen
Wissenschaften, \textbf{20}). Springer-Verlag, Berlin, 1-st ed.
1925, 2-nd ed. 1954. English transl.: \textsl{Problems and
Theorems in Analysis, Volume \textup{II}.} (Ser.: Grundlehren der
mathematischen Wissenschaften, \textbf{216}). Springer-Verlag,
Berlin\(\cdot\)Heidelberg\(\cdot\)New\,York, 1976, 301

\bibitem[Scho1]{Scho1} \textsc{Schoenberg,\,I.} \textsl{\"{U}ber
variationsvermindernde lineare Transformationen}.
Math.\,Zeitschrift, \textsf{32} (1930), pp.\,321\,-\,328.

\bibitem[Scho2]{Scho2} \textsc{Schoenberg,\,I.} \textsl{On totally
positive functions, Laplace Integral and entire functions of
Laguerre-Polya-Schur type.} Proc. Nat. Acad. Sci. U.S.A.,
\textbf{33} (1947), pp.11\,-\,17.

\bibitem[Scho3]{Scho3} \textsc{Schoenberg,\,I.} \textsl{On
variation-diminishing integral operators of the convolutions
 type}. Proc. Nat. Acad. Sci. U.S.A., \textbf{34} (1948).
 pp.\,164\,-\,169.

\bibitem[Scho4]{Scho4} \textsc{Schoenberg,\,I.J.} \textsl{On
P{\'o}lya frequency functions. I.{\small{The totally positive
functions and their Laplace transforms}}}. J. d'Analyse Math.,
Jerusalem, \textbf{1} (1951), pp.\,331\,-\,374.

\bibitem[Sch:\,Ges]{Sch:Gesh} \textsc{Schur,\,I.}:
\textsl{Gesammelte Abhandlungen
[Collected Works]. Vol. I, II, III.} Springer-Verlag, Berlin\(\cdot\)%
Heidelberg\(\cdot\)New\,York 1973.

\bibitem[Sch7]{Sch7a} \textsc{Schur,\,I.}: \textsl{Zwei S\"atze \"uber
algebraische Gleichungen mit lauter reellen Wurzeln} [\textsl{Two
theorems on algebraic equations with only real roots} - in
German].  Journ. f\"ur die reine und angew. Math., \textbf{144}
(1914), pp.\,75 - 88. \footnotesize Reprinted in: \small
\cite{Sch:Gesh}, Vol.\,II, pp.\,56 - 69.

\bibitem[Tho]{Thom} \textsc{Thoma,\,E.} \textsl{Die unzerlegbaren,
positiv-definiten Klassenfunktionen der abz{\"a}hlbar unendlichen,
symmetrischen Gruppe} [\textsl{The undecomposible
positive-definite class functions of the countably infinite
symmetric group} - in German]. Math. Zeitschr., \textbf{85}
(1964), pp.\,40\,-\,61.

\bibitem[Web]{Web} \textsc{Weber,\,H.} \textsl{Lehrbuch der Algebra}. Vieweg
Verlag, Braunschweig 1912 (2-nd ed.), New York, Chelsea 1961 (3-d %
ed).
\end{thebibliography}

\begin{thebibliography}{MoMo2}
\small
\bibitem[Alp]{Alp} \textsc{Alpay,\,D.} \textsl{Algorithme de
Schur, espaces {\`a} noyau reproduisant et th{\'e}orie des
syst{\`e}mes}. (Panorame and Synth{\`e}ses, \textsf{6}).
Soci{\`e}t{\`e} Math{\`e}matique de France, Paris, 1998.
viii+189pp. (French). Englisch transl.: \textsl{The Schur
Algorithm, Reproducing Kernel Spaces and System Theory.} (SMF/AMS
Texts and Monographs, \textsf{5}). American Mathematical Society,
Providence, RI; Soci{\`e}t{\`e} Math{\`e}matique de France, Paris,
2001. viii+150pp.

\bibitem[AlDy]{AlDy} \textsc{Alpay,\,D.} and \textsc{H. Dym}, 
\textsl{On applications
of reproducing kernel spaces to the Schur algorithm and rational $J$ unitary
factorization}. In \cite{S:Meth}, pp.\,89\,-\,159.  

\bibitem[Ausg]{Ausg}
\textsl{Ausgew\"ahlte Arbeiten zu den Urspr\"ungen der
Schur-Analysis} (German). [\textsl{Selected papers on the origin
of the Schur analysis}]. (\textsc{Fritzsche,~B.} and
\textsc{B.~Kirstein}
- editors) (German). {\rm(}Series: Teubner-Archiv zur %
Mathematik\,-\,\textsf{\textup{016}}). B.G.Teubner
Verlagsgesellschaft, Stuttgart$\cdot$Leipzig 1991, 290pp.

\bibitem[BDGPS]{BDGPS} \textsc{Bertin,\,M.J.,
A.\,Decomps-Guilloux, M.\,Grandet-Hugot, M.\.Pathiaux-Delefosse}
and \textsc{J.P.\,Schreiber}. \textsl{Pisot and Salem Numbers}.
Birkh{\"a}user-Verlag, Basel\(\boldsymbol{\cdot}\)Boston%
\(\boldsymbol{\cdot}\)Berlin, 1992, xiii+291pp.

\bibitem[Boy]{Boy} \textsc{Boyd,\,D.}: \textsl{Schur's algorithm for
bounded holomorphic functions.} Bull. London Math. Soc.,
\textbf{11} (1979), pp.\,145 - 150.

\bibitem[Con]{Con} \textsc{Constantinescu,\,T.}
\textsl{Schur Parameters, Factorization and Dilation Problems}.
(Operator Theory: Advances and Applications, \textbf{OT 82}).
Birkh{\"a}user Verlag, Basel, 1996. x+253\,pp.

\bibitem[DeGeK]{DeGeK} \textsc{PH.\,Delsarte, Y.\,Genin} and \textsc{Y.\,Kamp},
 \textsl{Schur parametrization of positive definite block-Toeplitz systems},
 SIAM J. Appl. Math., {\bf 36} (1979), 34-46.

\bibitem[DeDy]{DeDy} \textsc{Dewilde,\,P.} and \textsc{H.\,Dym}, 
\textsl{Lossless chain scattering matrices and optimum linear prediction: 
The vector case}, Circuit Theory and Appl., {\bf 9} (1981), 135-175.

\bibitem[DFK]{DFK} \textsc{Dubovoj,\,V.K., B.\,Fritzsche} and
\textsc{B.\,Kirstein}. \textsl{Matricial Version of the Classical
Schur Problem}. (Teubner\,-Texte zur Mathematik, Band
\textbf{129}).
B.G.\,Teubner Verlagsgesellschaft, %
Stuttgart\(\boldsymbol{\cdot}\)Leipzig, 1992.

\bibitem[FoFr]{FoFr} \textsc{C.\,Foias} and  \textsc{A.\,Frazho}, 
\textsl{The Commutant Lifting Approach to Interpolation Problems},  
(Operator Theory: Advances and Applications, \textbf{OT 44}).
Birkh{\"a}user Verlag, Basel, 1990.

\bibitem[FrKi]{FrKi} \textsc{Fritzsche,\,B.} and \textsc{B.\,Kirstein}.
\textsl{Schuranalysis - Umfassende Entfaltung einer mathematischen
Methode} [\textsl{The Schur analysis - the comprehensive scope of
a mathematical method} - in German.] In: Jahrbuch
\textsl{\"Uberblicke Mathematik 1992.} - \textsc{Chatterji,\,S.D.,
B.\,Fuchssteiner, U.\,Kulisch,\,R.\,Liedl} und
\textsc{W.\,Purkert} - editors, Vieweg, Braunschweig, 1992,
pp.\,117\,-\,136.

\bibitem[Ger1]{Gers1} \textsc{Geronimus,~Ya.L.}: %
\textsl{O polinomakh, ortogonal'nykh na kruge, o
trigonometricheskoi probleme momentov i ob associirovannykh s neyu
funktsiyakh tipa Carath\'{e}odory i Schur'a.  (On polynomials
orthogonal on the circle, on trigonometric moment-problem and on
allied Carath\'{e}odory and Schur functions.} Russian,
R\'{e}sum\'{e} Engl.) Matem. Sbornik. Nov.\,ser. \textbf{15}
(\textbf{42}), no.1 (1944), pp.~99-130.

\bibitem[Ger2]{Gers2} \textsc{Geronimus,~Ya.L.}: %
\textsl{Polinomy, ortogonal'nye na kruge, i ikh prilozheniya.}
Zapiski nauchno-issledovatel'skogo instituta matematiki i
mekhaniki i Khar'kovskogo matem. obshchestva, \textbf{19} (1948),
pp.~35-120 (Russian). English transl.: \textsl{Polynomials
orthogonal on a circle and their applications.} Amer. Math. Soc.
Transl. (ser. \textbf{1}), vol.\textbf{3} (1962), pp.~1-78.

\bibitem[Ger3]{Gers3} \textsc{Geronimus,~Ya.L.}: \textsl{Polinomy
Ortogonal'nye na Okruzhnosti i na Otrezke.} Fizmatgiz, Moskow,
1958 (In Russian). English transl.: \textsl{Polynomials Orthogonal
on a Circle and Interval.} Pergamon Press, New York, 1960.

\bibitem[Gol1]{Gol1} \textsc{Golinskii,~L.B}: \textsl{Schur
functions, Schur parameters and orthogonal polynomials on the unit
circle.} Zeitschrift f\"{u}r Analysis und ihre Anwendungen
\textbf{12} (1993), pp.\,457 -- 469.

\bibitem[Gol2]{Gol2} \textsc{Golinskii,~L.}: \textsl{On Schur
functions and Szeg\"{o} orthogonal polonomials.} In:
\textsl{Topics in Interpolation Theory} (Operator Theory:\,%
Advances and Applications,\,\textbf{OT\,95}). \ (\textsc{Dym,~H.,
B.~Fritzsche, V.~Katsnelson, B.~Kirstein} - editors.)
Birkh\"{a}user Verlag, Basel\(\cdot\)Boston\(\cdot\)Berlin 1997.

\bibitem[Kail]{Kail} \textsc{Kailath,\,T.} \textsl{A theorem of
I.\,Schur and its impact on modern signal processing}. In
\cite{S:Meth}, pp.\,9\,-\,30.

\bibitem[Kats]{Katsb} \textsc{Katsnelson,\,V.} \textsl{A generic Schur
 function is an inner one.} In: Interpolation Theory, System Theory and
Related Topics. The Harry Dym Anniversary Volume. (Operator
Theory: Advances and Applications). (\textsc{Alpay,\,D.,
I.\,Gohberg} and \textsc{V.\,Vinnikov}-editors). Birkh\"auser
Verlag, Basel%
\(\boldsymbol{\cdot}\)Boston\(\boldsymbol{\cdot}\)Berlin, 2002,
pp.\,249\,-\,293.

\bibitem[KaSa1]{KaSa1} \textsc{Kailath,\,T.} and
\textsc{A.H.\,Sayed}. \textsl{Displacement structure: theory and
applications}. SIAM Review, \textbf{37}:3, (1005),
pp.\,297\,-\,386.

\bibitem[KaSa2]{KaSa2} (\textsc{Kailath,\,T.} and
\textsc{A.H.\,Sayed}-aditors.) \textsl{Fast Reliable Algorithms
for Matrices with Structure.} Society of Industrial and Applied
Mathematics (SIAM), Philadelphia, PA, 1999, xvi+342pp.

\bibitem[Khru1]{Khru1} \textsc{Khrushchev,~S.}: \textsl{Parameters of
orthogonal polynomials.} %
In: \textsl{"Methods of Appriximation Theory in Complex Analysis
and Mathematical Physics. Selected papers from the international
seminar held in Leningrad, May 13-26, 1991."}
\textsc{Gonchar,~A.A} and \textsc{E.~Saff} - editors. Lecture
Notes in Math., vol.\,\textbf{1550}. Reissued by Springer-Verlag,
Berlin 1993, originally published by "Nauka", Moskow 1993,
pp.\,185-191.

\bibitem[Khru2]{Khru2} \textsc{Khrushchev,~S.V.}:
\textsl{Schur's algorithm, orthogonal polynomials and convergence
of Wall's continued fractions in \(L^2(\mathbb{T})\).} Journal of
Approximation Theory, vol.\,\textbf{108}, no.2, (2001),
pp.\,161-248.

\bibitem[Khru3]{Khru3} \textsc{Khrushchev,~S.V.}:
\textsl{A singular Riesz product in the Nevai class and inner
functions with  the Schur parameters \(\bigcap\limits_{p>2}l^p\).}
Journal of Approximation Theory, vol.\,\textbf{108}, no.2, (2001),
pp.\,249-255.

\bibitem[LevKai]{LevKai} \textsc{Lev-Ari,\,H.} and \textsc{T.\,Kailath}.
\textsl{Triangular factorization of structured hermitian
matrices}. In \cite{S:Meth},\,pp.\,301\,-\,323.

\bibitem[Levsn]{Levsn} \textsc{Levinson,\,N.} \textsl{The Wiener's RMS (root
mean square) error criterion in filter desighn and prediction.}
Journal of mathematics and physics, \textbf{25} (1946),
261\,-\,278.

\bibitem[Nja]{Nja} \textsc{Nj{\aa}stad,~O.} \textsl{Convergence of the Schur
algorithm.} Proc. Amer. Math. Soc. \textbf{110}\,,\,No.4\,,
(1990), pp.\,1003-1007.

\bibitem[Oue]{Oue} \textsc{Ouelette,\,D.V.} \textsl{Schur
complements and statistics}. Linear Algebra and its Applications,
\textbf{36}, (1981), pp.\,187\,-\,295.

\bibitem[PiNe]{PiNe}  \textsc{Pint\'{e}r,\,F.} and \textsc{P.\,Nevai}. %
\textsl{Schur functions and orthogonal polynomials on the unit
circle.} In: \textsl{Approximation Theory and Function Series.
Budapest, 1995}, pp.\,293 - 306. \textsc{V\'{e}rtesi,\,P.,
L.\,Leindler, F.\,M\'{o}ricz, Sz.\,R\'{e}v\'{e}sz,\, J.\,Szabados}
and \textsc{V.\,Totik} - editors. J\'{a}nos Bolyai Mathematical
Society, Budapest 1996.

\bibitem[Rakh]{Rakh}
\textsc{Rahmanov,~E.A.(=Rakhmanov,~E.A.)} \textsl{Ob asimptotike
otnosheniya ortogonal'nykh polinomov.\,II.}  %
Matem. Sbornik. Nov.ser. \textbf{118}\,(\textbf{160}), no.\,1
(1982), pp.104-117 (In Russian). English transl.: \textsl{On the
asymptotics of the ratio of orthogonal polynomials.\ II.}\ Math.
USSR Sbornik \textbf{46} (1983), pp.\,105-117.

\bibitem[Sch:\,Ges]{Sch:Ges} \textsc{Schur,\,I.}:
\textsl{Gesammelte Abhandlungen
\textup{[}Collected Works\,\textup{]}. Vol. I, II, III.} Springer-Verlag,
Berlin\(\boldsymbol{\cdot}\)%
Heidelberg\(\boldsymbol{\cdot}\)New\,York 1973.

\bibitem[Sch9]{Sch9k} \textsc{Schur,~I}.: \textsl{\"Uber Potenzreihen, die im
Innern des Einheitskreises beschr\"ankt \mbox{sind,} I\,.\,} (in
German) J. reine und angewandte Math. \textbf{147}(1917), 205 -
232. \footnotesize Reprinted in:
 \small \cite{Sch:Ges}, Vol.\,II,
pp.\,137 - 164. \footnotesize Reprinted also in: \small
\cite{Ausg}, pp. 22-49. {\footnotesize English translation}:
\small \textsl{On power series which are bounded in the interior
of the unit circle.\,I.}, In: \,\cite{S:Meth}, pp. 31-59.

\bibitem[Sch10]{Sch10k} \textsc{Schur,~I.}: \textsl{\"Uber Potenzreihen, die im
Innern des Einheitskreises beschr\"ankt \mbox{sind,} II\,.\,} (in
German). J. reine und angewandte Math. \textbf{148} (1918), 122 -
145. \footnotesize Reprinted in: \small \cite{Sch:Ges}, Vol.\,II,
pp.\,165 - 188.
 \footnotesize Reprinted also in: \small\cite{Ausg}, pp.
50-73. {\footnotesize English translation}: \small\textsl{On power
series which are bounded in the interior of the unit circle.\,II.}
In: \cite{S:Meth}, pp. 61 - 88.

\bibitem[Sm]{Sm} \textsc{\v{S}mul'yan,\,Yu.L.}
(=\textsc{Shmul'yan,\,Yu.L.}) \textsl{Operatorny\u{\i} integral
Hellingera} (Russian). Matem. Sbornik \textbf{47} (\textbf{91}):4
(1959), pp.\,381\,-\,430. English transl.: \textsl{A Hellinger
operator integral}, Amer. Math.Soc. Transl. (\textbf{2})
\textbf{22} (1962), 289\,-\,337.


\bibitem[TIT]{TIT} \textsl{Topics in Interpolation Theory}.
(\textsc{Dym,\,H., B.\,Fritzsche, V.\,Katsnelson, B.\,Kirstein} -
editors). (Operator Theory: Advances and Applications, \textbf{OT
95}), Birkh\"auser Verlag,
Basel\(\boldsymbol{\cdot}\)Boston\(\boldsymbol{\cdot}\)Berlin,
1997.

\bibitem[S:Meth]{S:Meth} \textsl{I. Schur Methods in Operator Theory and Signal
Processing}  (Operator Theory: Advances and Applications, {\bf OT
18}) (\textsc{Gohberg,~I.} - editor). Birkh\"auser Verlag,
Basel\(\cdot\)\,Boston\(\cdot\)\,Stuttgart, 1986.
\end{thebibliography}

\begin{thebibliography}{MoMo2}

\bibitem[Sch1]{Sch1x} \textsc{Schur,\,I.} \textsl{\"Uber
vertauschbare lineare Differentialausdr\"ucke} (German)
[\textsl{On permutable differential expressions}].
Sitzungsberichte der Berliner Mathematischen Gesellschaft,
\textbf{4} (1905), pp.\, 2 - 8. \footnotesize Reprinted in: \small
\cite{Sch:Gesx}, Vol.\,I, pp.\,170 - 176.

\bibitem[Sch2]{Sch2x} \textsc{Schur,\,I.}: \textsl{\"Uber die
charakteristischen Wurzeln einer linearen Substitution mit einer
Anwendung auf die Theorie der Integralgleichungen}
(German)[\textsl{On the characteristic roots of a linear
substitution with an application to the theory of integral
equations}]. Mathematische Annalen, \textbf{66} (1909), pp.\,488 -
510. \footnotesize Reprinted in: \small \cite{Sch:Gesx}, Vol.\,I,
pp.\, 272 - 294.

\bibitem[Sch3]{Sch3x} \textsc{Schur, I.} \textsl{Zur Theorie der linearen
homogenen Integralgleichungen}(German) [\textsl{On the theory of
linear homogeneous integral equations}]. Math. Annalen,
\textbf{67} (1909), pp.\,336 - 339. \footnotesize Reprinted in:
\small \cite{Sch:Gesx}, Vol.\,I, pp.\,312 - 345.

\bibitem[Sch4]{Sch4x} \textsc{Schur,\,I.}: \textsl{Bemerkungen zur
Theorie der beschr\"ankten Bilinearformen mit unendlich vielen
Ver\"anderlichen} (German) [\textsl{Remarks on the theory of
bounded bilinear forms with infinitely many variables}]. Journ.
f\"ur reine und angew. Math., \textbf{140} (1911), pp.\,1 - 28.
\footnotesize Reprinted in: \small \cite{Sch:Gesx}, Vol.\,I,
pp.\,464 - 491.

\bibitem[Sch5]{Sch5x} \textsc{Schur,\,I.}: \textsl{\"Uber einen Satz von
C.\,Carath\'eodory} (German) [\textsl{On a theorem of
C.\,Carath\'eodory}]. Sitzungsberichte der K\"onigl. Preuss. Akad.
der Wiss., (1912), pp.\,4 - 15. \footnotesize Reprinted in: \small
\cite{Sch:Gesx}, Vol.\,II, pp.\,12 - 23.

\bibitem[Sch6]{Sch6x} \textsc{Schur,\,I.}:
\textsl{\"Uber die \"Aquivalenz der Ces\`aroschen und
H\"olderschen Mittelwerte} (German) [\textsl{On the equivalence
of Ces\`aro's and H\"older means}]. Mathematische Annalen,
\textbf{74} (1913), pp.\,447 - 458. \footnotesize Reprinted in:
\small \cite{Sch:Gesx}, Vol.\,II, pp.\,44 - 55.

\bibitem[Sch7]{Sch7x} \textsc{Schur,\,I.}: \textsl{Zwei S\"atze \"uber
algebraische Gleichungen mit lauter reellen Wurzeln} (German) [\textsl{Two
theorems on algebraic equations with only real roots}].  
Journ. f\"ur die reine und angew. Math., \textbf{144}
(1914), pp.\,75 - 88. \footnotesize Reprinted in: \small
\cite{Sch:Gesx}, Vol.\,II, pp.\,56 - 69.

\bibitem[Sch8]{Sch8x} \textsc{Schur,\,I.} \textsl{\"Uber die
Entwicklung einer gegebenen Funktion nach den Eigenfunktionen
eines positiv definiten Kerns} (German) [\textsl{On the expansion
of a given function in eigenfunctions of a positive definite
kernel}]. Schwarz-Festschrift (1914), pp.\,392 - 409.
\footnotesize Reprinted in: \small \cite{Sch:Gesx}, Vol.\,II,
pp.\,70 - 87.

\bibitem[Sch9]{Sch9x} \textsc{Schur,~I}.: \textsl{\"Uber Potenzreihen, die im
Innern des Einheitskreises beschr\"ankt \mbox{sind,} I\,.\,} (German) 
J. reine und angewandte Math. \textbf{147}(1917), 205 -
232. \footnotesize Reprinted in:
 \small \cite{Sch:Gesx}, Vol.\,II,
pp.\,137 - 164. \footnotesize Reprinted also in: \small
\cite{Ausgx}, pp. 22-49. {\footnotesize English translation}:
\small \textsl{On power series which are bounded in the interior
of the unit circle.\,I.}, In: \,\cite{S:Methx}, pp. 31-59.

\bibitem[Sch10]{Sch10x} \textsc{Schur,~I.}: \textsl{\"Uber Potenzreihen, die im
Innern des Einheitskreises beschr\"ankt \mbox{sind,} II\,.\,} (in
German). J. reine und angewandte Math. \textbf{148} (1918), 122 -
145. \footnotesize Reprinted in: \small \cite{Sch:Gesx}, Vol.\,II,
pp.\,165 - 188.
 \footnotesize Reprinted also in: \small\cite{Ausgx}, pp.
50-73. {\footnotesize English translation}: \small\textsl{On power
series which are bounded in the interior of the unit
circle.\,II.} In: \cite{S:Methx}, pp. 61 - 88.

\bibitem[Sch11]{Sch11x} \textsc{Schur,\,I.} \textsl{\"Uber die Verteilung
der Wurzeln bei gewissen algebraischen Gleichungen mit
ganzzahligen Koeffizienten} (German) [\textsl{On the distribution
of the roots of some algebraic equations with  integer
coefficients}]. Math. Zeitschrift, \textbf{1} (1918),
pp.\,377\,-\,402. \footnotesize Reprinted in: \small
\cite{Sch:Gesx}, Vol.\,II, pp.\,213\,-\,238.
\bibitem[Sch12]{Sch12x} \textsc{Schur,~I.}: \textsl{\"Uber die
Koeffizientensummen einer Potenzreihe mit positivem reellen Teil}
(German) [\textsl{On the sum of the coefficients of a power
series with  positive real part}]. Archiv der Math. und Phys.
(3) \textbf{27} (1918), pp.\,126 - 135. \footnotesize Reprinted
in: \small \cite{Sch:Gesx}, Vol.\,II, pp.\, 249 - 258.

\bibitem[Sch13]{Sch13x} \textsc{Schur,\,I.}: \textsl{\"Uber das Maximum
des absoluten Betrages eines Polynoms in einem gegebenen
Intervall} (German) [\textsl{On the maximum of the absolute value
of a polynomial on a given interval}]. Math. Zeitschrift,
\textbf{4} (1919), pp.\,271 -287. \footnotesize Reprinted in:
\small \cite{Sch:Gesx}, Vol.\,II, pp.\,259 - 275.

\bibitem[Sch14]{Sch14x}\textsc{Schur,\,I.} \textsl{\"Uber einen von Herrn
L.\,Lichtenstein benuzten Integralsatz} (German) [\textsl{On a
theorem on integrals used by L.\,Lichtenstein}]. Math. Zeitschr.
\textbf{7} (1920), pp.\,232 - 234. \footnotesize Reprinted in:
\small \cite{Sch:Gesx}, Vol.\,II, pp.\,286  - 288.

\bibitem[Sch15]{Sch15x} \textsc{Schur,\,I.} \textsl{\"Uber algebraische
Gleichungen, die nur Wurzeln mit negativen Realteilen besitzen}
(German) [\textsl{On algebraic equations which possess only
roots with  negative real parts}]. Zeitschr. f\"ur angew. Math.
und Mech., \textbf{1} (1921), pp.\.307\,-\,311. {\footnotesize
Reprinted in:} \small \cite{Sch:Gesx}, Vol.\,II, pp.\,322\,-\,326.

\bibitem[Sch16]{Sch16x} \textsc{Schur,\,I.}:
 \textsl{\"Uber lineare Transformationen
in der Theorie der unendlichen Reien} (German) [\textsl{On
linear transformations in the theory of infinite series}]. Journ.
f\"ur die reine und angew. Math., \textbf{151} (1921), pp.\,79 -
121. \footnotesize Reprinted in: \small \cite{Sch:Gesx}, Vol.\,II,
pp.\,289 - 321.

\bibitem[Sch17]{Sch17x} \textsc{Schur,\,I.} \textsl{Ein Beitrag zur
Hilbertschen Theorie vollstetigen quadratischen Formen}
(German) [\textsl{A contribution to the theory of completely
continuous quadratic forms}]. Math. Zeitschr., \textbf{12} (1922),
287\,-\,297. \footnotesize Reprinted in: \small \cite{Sch:Gesx},
Vol.\,II, pp.\,402 - 412.

\bibitem[Sch18]{Sch18x} \textsc{Schur,\,I.}: \textsl{\"Uber %
eine Klasse von Mittelbildungen mit Anwendungen auf die
Determinantentheorie} (German) [\textsl{On a class of averaging
mappings with applications to the theory of determinants}].
Sitzungsberichte der Berliner Mathematischen Gesellschaft,
\textbf{22} (1923), pp.\,9 - 20. \footnotesize Reprinted in:
\small \cite{Sch:Gesx}, Vol.\,II, pp.\,416 - 427.

\bibitem[Sch19]{Sch19x} \textsc{Schur,\,I.}
\textsl{Einige Bemerkungen zur Theorie der unendlichen Reihen}
(German) [\textsl{A remark on the theory of infinite series}].
Sitzungsberichte der Berliner Mathematischen Gesellschaft,
\textbf{29}, (1930), pp.\,3 - 13. \footnotesize Reprinted in:
\small \cite{Sch:Gesx}, Vol.\,III, pp.\,216 - 226.

\bibitem[Sch20]{Sch20x} \textsc{Schur,\,I.} \textsl{Affektlose Gleichungen
in der Theorie der Laguerreschen und Hermiteschen Polynome}
(German) [\textsl{On affect free equations in the theory of 
Laguerre polynomials and Hermite polynomials}]. Journ. f\"ur die
reine und angewandte Mathematik, \textbf{165} (1931),
pp.\,52\,-\,58. Reprinted in: \small\cite{Sch:Gesx}, Vol.\,III,
pp.\,227 - 233.

\bibitem[Sch21]{Sch21x} \textsc{Schur,\,I.}
\textsl{On Faber Polynomials}. American Journ. of Math.
\textbf{67} (1945), pp.\,33 - 41. \footnotesize Reprinted in:
\small \cite{Sch:Gesx}, Vol.\,III, pp.\,361 - 369.

\bibitem[Sch22]{Sch22x} \textsc{Schur, I.}
\textsl{Identities in the theory of power series.} American
Journ. of Math. \textbf{69} (1947), pp.\,14 - 26. \footnotesize
Reprinted in: \small \cite{Sch:Gesx}, Vol.\,III, pp.\,379 - 391.

\textbf{Joint papers}:
\bibitem[KnSch]{KnSchx} \textsc{Knopp,\,K.} and \textsc{I.\,Schur}.
\textsl{\"Uber die
Herleitung der Gleichung %
\(\sum\limits_{n=1}^{\infty}\frac{1}{n^2}=\frac{\pi^2}{6}\)%
} (German) [\textsl{On a derivation of the equation %
\(\sum\limits_{n=1}^{\infty}\frac{1}{n^2}=\frac{\pi^2}{6}\)%
}]. Archiv der Math. und Phys. (3)\textbf{27} (1918),
pp.\,174\,-\,176. \footnotesize Reprinted in: \small
\cite{Sch:Gesx}, Vol.\,II, pp.\,246\,-\,248.


\bibitem[PS]{PSx} \textsc{Polya,\,G.} and \textsc{I.\,Schur}:
\textsl{\"Uber zwei Arten von Faktorenfolgen in der Theorie der
algebraischen Gleichungen} (German) [\textsl{On two types of
multiplier sequences in the theory of algebraic equations}].
Journ. f\"ur reine und angew. Math. \textbf{144} (1914), pp.\,89 -
113. \footnotesize Reprinted in: \small \cite{Sch:Gesx}, Vol.\,II,
pp.\,88 - 112.

\bibitem[SchSz]{SchSzx} \textsc{Schur,\,I.} and
\textsc{G.\,Szeg\"o.} \textsl{\"Uber die Abschnitte einer im
Einheitkreise beschr\"ankten Potenzreihe} (German) [\textsl{On
the truncation of a power series bounded in the unit disk}].
Sitzungsberichte der Preu\ss ischen Akademie der Wissenschaften
1925, Physikalisch-Mathematische Klasse, pp.\,545 - 560.
\footnotesize Reprinted in: \small \cite{Sch:Gesx}, Vol.\,II,
pp.\,249 - 258. \footnotesize Reprinted in: \small \cite{Sch:Gesx},
Vol.\,III, pp.\,27 - 42.

\textbf{Sources in which Schur papers were reprinted or
translated}:
\bibitem[Sch:\,Ges]{Sch:Gesx} \textsc{Schur,\,I.}:
\textsl{Gesammelte Abhandlungen \textup{[}Collected
Works\textup{]}. Vol. I, II, III.} Springer-Verlag,
Berlin\(\cdot\) Heidelberg\(\cdot\)New\,York, 1973.

\bibitem[Ausg]{Ausgx}
\textsl{Ausgew\"ahlte Arbeiten zu den Urspr\"ungen der
Schur-Analysis} (German) [\textsl{Selected papers on the origin
of  Schur analysis}]. (\textsc{Fritzsche,~B.} and
\textsc{B.~Kirstein}, 
editors) (German). {\rm(}Series: Teubner-Archiv zur %
Mathematik\,-\,\textsf{\textup{016}}). B.G.Teubner
Verlagsgesellschaft, Stuttgart$\cdot$Leipzig 1991, 290pp.

\bibitem[S:Meth]{S:Methx} \textsl{I. Schur Methods in Operator Theory
and Signal Processing}  (Operator Theory: Advances and
Applications, {\bf OT 18}) (\textsc{Gohberg,~I.}, editor).
Birkh\"auser Verlag, Basel\(\cdot\)\,Boston\(\cdot\)\,Stuttgart
1986.


\end{thebibliography}
\end{document}